\documentclass[review]{elsarticle}

\usepackage{lineno,hyperref}
\modulolinenumbers[5]

\usepackage{graphicx}
\usepackage{subfig}
\usepackage{amsmath}
\usepackage{amssymb}
\usepackage{lscape} 

\usepackage{nomencl} 
\usepackage{multicol}

\usepackage{etoolbox}
\renewcommand\nomgroup[1]{%
	\item[\bfseries
	\ifstrequal{#1}{A}{List of Symbols}{%
		\ifstrequal{#1}{B}{List of Subscripts}{}}%
	]}

\journal{Journal of Computational and Applied Mathematics, }

\begin{document}

\begin{frontmatter}

\title{Development and Preliminary Assessment of the Open-Source CFD toolkit SU2 for Rotorcraft Flows}

\author[mymainaddress,mysecondaryaddress]{Myles Morelli\corref{mycorrespondingauthor}}
\cortext[mycorrespondingauthor]{Corresponding author, email address: mylescarlo.morelli@polimi.it}

\author[mymainaddress]{Tommaso Bellosta}

\author[mymainaddress]{Alberto Guardone}

\address[mymainaddress]{Department of Aerospace Science and Technology, Politecnico di Milano, Italy}
\address[mysecondaryaddress]{CFD Laboratory, School of Engineering, University of Glasgow, United Kingdom}

\begin{abstract}
	Computational aerodynamic analyses of rotorcraft main rotor blades are performed in both hover and forward flight. The open-source SU2 code is used for rotor performance prediction. The core of the code is the set of RANS equations, which are solved for determining the flow. In hover, both steady-state and time-accurate modelling techniques of varying complexity are used and assessed. Simulation specific parameters which have a significant influence on the solution are also addressed. In forward flight, the code is developed to include the main rotor blade kinematics which is a prerequisite for modelling a trimmed rotor. Two databases are used for the validation of the rotor performance prediction. The renowned Caradonna-Tung experimental tests of a model rotor are used to evaluate the pressure distribution along the blade during hover. The extensive aerodynamic and aeroacoustic data survey of the AH-1G Cobra helicopter is used to assess the pressure distribution at different advancing and retreating azimuth angles during forward flight. The prediction capabilities of the solver in terms of rotor performance are demonstrated and are overall in good agreement with the measured data.
\end{abstract}


\end{frontmatter}

%


\section{Introduction}
\label{chap4:sec:intro}
Rotorcraft have the unique ability to be able to vertically take-off and land. This allows for their operation in highly demanding and challenging flight scenarios where conventional fixed-wing aircraft cannot pursue. Their frequent use for demanding operations however can lead to dangerous in-flight situations which can be seen from the relatively high number of accident reports from the National Transportation Safety Board (NTSB) \cite{NTSB}. Tools that allow the analysis of rotorcraft are therefore essential for supporting design and for reducing accident rates. Sustained computational development over the past 30 years \cite{strawn200630} means that many codes are now capable of modelling rotorcraft behaviour. Computational advancements have enabled CFD-based methods to simulate and further understand complex rotorcraft aerodynamics which are rich in flow physics. 

NASA have a long history of rotorcraft CFD code development with OVERFLOW \cite{jespersen1997recent} and FUN3D \cite{anderson1994implicit}. Of the two codes, OVERFLOW has been used more extensively for rotorcraft simulations. OVERFLOW uses body-fitted structured grids near solid surfaces and automatically generated Cartesian grids in the background and is a finite-difference node-based solver. FUN3D uses unstructured grids throughout the domain and is a finite volume node-based solver. Academia also plays a leading role in rotorcraft CFD code development. In-house codes HMB from the University of Glasgow \cite{steijl2006framework} and ROSITA from Politecnico di Milano \cite{biava2012simulation} are both finite-volume solvers and utilize structured Chimera multi-block grids to account for the blade motion. Each of the codes have been part of large-scale collaborations for the assessment of their predictive capabilities such as the GOAHEAD project \cite{antoniadis2012assessment}. 
Another code developed within academia is the TURNS research code from the University of Maryland \cite{srinivasan1992flowfield, srinivasan1993turns}. The TURNS code has Chimera overset grid capabilities and uses a finite difference numerical algorithm that evaluates the inviscid fluxes using an upwind-biases flux scheme. European research centers also have their own codes which have been heavily developed. French and German research institutes ONERA and DLR have their own multi-block Chimera based rotorcraft codes. ONERA have the elsA solver \cite{gazaix2002elsa} and DLR have the FLOWer solver \cite{raddatz2005block}.

Despite all of the codes showing excellent predictive capabilities, none of the mentioned codes are freely available in the open-source domain which hinders technological developments. Furthermore, Validation and Verification (VnV) is shown to be increasingly important so to identify a wide variety of physical modelling, discretization, and solution errors \cite{oberkampf1998issues}. To that end, an open-source code with an active and growing community of users provides a platform for extensive VnV and innovative new solutions. 

This work looks to provide the first open-source and validated rotorcraft CFD code by developing the well established SU2 code \cite{economon2016su2}. The main contribution of this work is the introduction of the blade motion and an approach for modelling forward flight. The contemporary nature of the code means that it can benefit from new and effective numerical techniques. The recent implementation by Gori et al. \cite{gori2017sliding} of a method for dealing with non-conformal boundary interfaces using the supermesh technique is one example of this. Another example of this is radial basis function mesh deformation which can now utilize computationally efficient algorithms for dealing with large-scale problems. The outline of this paper is as follows; the physical modelling of the RANS equations is discussed in Section~\ref{chap4:sec:rans}, the numerical implementation is described in Section~\ref{chap4:sec:numerical_implementation}, the validation work and results of an isolated rotor in hover and forward flight are discussed in Section~\ref{chap4:sec:validation}, and finally the main talking points of the work are concluded in Section~\ref{chap4:sec:conclusion}.

\section{Physical Modelling}
\label{chap4:sec:rans}
The following section on the physical modelling of the Reynolds-averaged Navier-Stokes equations inside of the SU2 code is primarily a summary of the work from Economon et al. \cite{economon2016su2} and is here to give context to the numerical implementation in Section~\ref{chap4:sec:numerical_implementation}.  

\subsection{Reynolds-Averaged Navier-Stokes Equations}
Within this framework, we are interested in time-accurate turbulent flow around rotor blades with arbitrary motion. Therefore we are concerned with compressible flow governed by the Reynolds-averaged Navier-Stokes (RANS) equations. These mass, momentum, and energy conservation equations can be expressed in arbitrary Lagrangian-Eulerian differential form as,

\begin{equation}
\label{chap4:eq:rans}
\left\{\begin{matrix}
\begin{array}{cclllc}
\mathcal{R}(\boldsymbol{U}) & 
= & 
\frac{\partial U}{\partial t} + \triangledown \cdot \boldsymbol{F}_{ale}^{c} - \triangledown \cdot ( \mu_{tot}^1 \boldsymbol{F}^{\upsilon 1} + \mu_{tot}^2 \boldsymbol{F}^{\upsilon 2}) - \boldsymbol{Q} = 0 & 
\textrm{in} & \Omega , & 
t > 0\\ 
\boldsymbol{\upsilon} & 
= & 
\boldsymbol{u}_{\Omega} & 
\textrm{on} & 
S & \\ 
\partial_n T & 
= & 
0 & 
\textrm{on} & 
S & \\ 
(W)_+ & 
= & 
W_{\infty} & 
\textrm{on} & 
\Gamma_{\infty} & 
\end{array}
\end{matrix}\right.    
\end{equation}
with the vector of conservative variables being represented by, $\boldsymbol{U} = \left \{ \rho, \, \rho \boldsymbol{\upsilon}, \, \rho E \right \}^{T}$, inside the flow domain, $\Omega$. The term, $\boldsymbol{\upsilon} = \boldsymbol{u}_{\Omega}$, denotes the no-slip condition on the surface, $S$. The expression, $ \partial_n T = 0 $, represents the adiabatic condition on the surface, $S$. The final condition, $ (W)_+ = W_{\infty} $, is the characteristic-based boundary condition at the far-field, $\Gamma_{\infty}$. The convective fluxes, viscous fluxes and source terms within Eq.~\ref{chap4:eq:rans} can then be respectively expressed as,

\begin{equation}
\boldsymbol{F}_{ale}^{c} = \begin{Bmatrix}
\rho (\boldsymbol{\upsilon} - \boldsymbol{u}_{\Omega}) \\
\rho \boldsymbol{\upsilon} \otimes (\boldsymbol{\upsilon} - \boldsymbol{u}_{\Omega}) + \bar{\bar{I}}_{p}\\
\rho E (\boldsymbol{\upsilon} - \boldsymbol{u}_{\Omega}) + p\boldsymbol{\upsilon}
\end{Bmatrix}
, \,
\boldsymbol{F}^{\upsilon1} = \begin{Bmatrix}
\cdot \\ 
\bar{\bar{\tau}}\\ 
\bar{\bar{\tau}} \cdot \boldsymbol{\upsilon}
\end{Bmatrix}
, \,
\boldsymbol{F}^{\upsilon2} = \begin{Bmatrix}
\cdot \\ 
\cdot \\ 
c_p \triangledown T
\end{Bmatrix}
, \,
\boldsymbol{Q} = \begin{Bmatrix}
q_{\rho}\\ 
\boldsymbol{q}_{\rho \upsilon}\\ 
q_{\rho E}
\end{Bmatrix}    
\end{equation}
where the fluid density, flow velocity vector, static pressure, temperature, and specific heat are respectively given by the orthodox notation $\rho, \, \boldsymbol{\upsilon}, \, p,\, T,\, c_p$. Of the remaining variables; $E$, is the total energy per unit mass; $\bar{\bar{\tau}}$ is the viscous stress tensor; $\boldsymbol{u}_{\Omega}$ is the grid velocity.

\subsection{Turbulence Models}
Here we are concerned with rotorcraft which during forward flight, climb, descent, and manoeuvre can be characterised by unsteady turbulent flows. With this in mind, the solution of the unsteady Reynolds average Navier-Stokes equations need to be solved which requires the inclusion of a turbulence model. Using the Boussinesq hypothesis \cite{wilcox1998turbulence}, the effect of turbulence is represented as an increase in the viscosity. The total viscosity is then separated into laminar and turbulent viscosity and can be respectively denoted as $\mu_{dyn}$ and $\mu_{tur}$. The laminar viscosity is determined based upon Sutherland's law \cite{sutherland1893lii}. The total viscosity as part of the momentum and energy equations in Eq.~\ref{chap4:eq:rans} is then substituted by

\begin{equation}
\mu_{tot}^{1} = \mu_{dyn} + \mu_{tur}, \quad \mu_{tot}^{2} = \frac{\mu_{dyn}}{Pr_{l}} + \frac{\mu_{tur}}{Pr_{t}}
\end{equation}
where the laminar and turbulent Prandtl numbers are given by $Pr_{l}$ and $Pr_{t}$. 

The turbulent viscosity is computed using a suitable turbulence model which itself is dependent upon the flow state and a new set of variables, $\hat{\nu}$, to represent the turbulence, such that, $\mu_{tur} = \mu_{tur}(\boldsymbol{U},\, \hat{\nu})$. One of the most widely used models for aeronautical attached flows is the one-equation Spallart-Allmaras (SA) turbulence model \cite{spalart1992one}. 

\subsection{Rotating Frame of Reference}
The unique ability of rotorcraft to be able to hover as well as to climb and descend vertically in axial flight helps to distinguish them from other aircraft. During these specific operational flight conditions, the flow around the main rotor can be considered as a steady rotation. Under this assumption, it is then possible to transform the unsteady problem into a steady problem to improve the efficiency of the simulation. This is possible by transforming the system of governing equations in Eq.~\ref{chap4:eq:rans} into a rotating frame of reference which rotates at the constant rotational velocity of the main rotor blades. This modification to the system was implemented by Economon et al. \cite{economon2013viscous} inside the SU2 code and can be written as,

\begin{equation}
\frac{\partial \boldsymbol{U}}{\partial t} = 0, \quad \boldsymbol{u}_{\Omega} = \boldsymbol{\omega} \times \boldsymbol{r}, \quad \boldsymbol{Q} =\begin{Bmatrix}
\cdot\\ 
-\rho (\boldsymbol{\omega} \times \boldsymbol{r})\\ 
\cdot
\end{Bmatrix} 
\end{equation}
where the rotational velocity vector of the rotating frame of reference is specified by $\boldsymbol{\omega} = \{ \omega_{x}, \, \omega_{y}, \, \omega_{z}\}^{T}$ and the radial distance from the center of rotation is specified by $\boldsymbol{r}$.

\section{Numerical Implementation}
\label{chap4:sec:numerical_implementation}
Alongside the physical modelling in the paper from Economon et al. \cite{economon2016su2} there is also a full description of the numerical implementation. In detail, it discusses the spatial and temporal integration which are the foundation for more application-specific methods such as the harmonic balance and non-conformal boundary interface treatment. The concurrent implementation for modelling rotor blade kinematics will finally be introduced. 

\subsection{Spatial Integration}
The flow equations are solved numerically via a finite volume method \cite{versteeg2007introduction} which is applied on unstructured grids with an edge-based structure. Integrating the governing equations over a control volume and using the divergence theorem to obtain the semi-discretized form gives,

\begin{equation}
\label{eq:FVM}
\int_{\Omega_{i}} \frac{\partial \boldsymbol U}{\partial t} d{\Omega} + \sum_{j \in \mathcal{N}(i)} (\tilde{\boldsymbol{F}_{ij}^{c}} + \tilde{\boldsymbol{F}_{ij}^{v}})\Delta S_{ij} - \boldsymbol{Q}\left | \Omega_{i} \right | = \int_{\Omega_{i}} \frac{\partial \boldsymbol{U}}{\partial t} d{\Omega} + \mathcal{R}_{i} (\boldsymbol{U}) = 0
\end{equation}
where $\mathcal{R}_{i} (\boldsymbol{U}) $ includes the convective and viscous fluxes integrated over the surface area of a control volume and any source terms. The numerical approximations of the convective and viscous fluxes are expressed respectively by, $\tilde{\boldsymbol{F}_{ij}^{c}} $ and $ \tilde{\boldsymbol{F}_{ij}^{v}}$. The area of the face belonging to the edge $ij$ is represented by $\Delta S_{ij}$ and the set of neighbouring vertices to vertex $i$ is expressed by $\mathcal{N}(i)$. The volume of control volume $i$ in the domain is denoted by $\left | \Omega_{i} \right |$. 

The convective and viscous fluxes are evaluated at the mid-point of an edge. Centered or upwind schemes such as the Jameson-Schmidt-Turkel scheme \cite{jameson1981numerical} or the approximate Riemann solver of Roe \cite{roe1981approximate} are then used to discretize the convective fluxes. Second-order of accuracy of upwind schemes is achieved by using Monotone UpStream-Centered schemes for conservation Laws (MUSCL) \cite{van1979towards} to reconstruct the variables on the cell interfaces. Green-Gauss or weighted least-squares methods are used approximate the spatial gradients of the flow at the cell faces to determine the viscous fluxes. 

\subsection{Time Integration}
In conjunction with the spatial discretization, there is also the requirement for temporal discretization of the governing equations and so Equation~\ref{eq:FVM} is further discretized over a control volume, $  \left | \Omega_{i} \right | $, in time such that it becomes,


\begin{equation}
\frac{d}{dt} (\left | \Omega_{i} \right | \boldsymbol{U})  + \mathcal{R}_{i} (\boldsymbol{U}) = 0
\end{equation}
and a backward Euler scheme can be used to evaluate the solution state at the updated time for steady problems.

The temporal discretization for unsteady simulations is achieved using a dual-time stepping approach which allows for second order accuracy in time \cite{jameson1991time}. Using this approach, an additional fictitious time is introduced so that the unsteady problem becomes a series of pseudo-steady problems. The fictitious time is introduced in front of the governing equations such that,


\begin{equation}
\label{chap4:eq:dual_time}
\left | \Omega_{i} \right |  \frac{\partial \boldsymbol{U}_{n}}{\partial \tau} + \mathcal{R}^{*}(\boldsymbol{U}_{n}) = 0
\end{equation}

with the subscript $ n $ denoting the physical time level and to achieve second-order backward difference in time,

\begin{equation}
\mathcal{R}^{*}(\boldsymbol{U}_{n}) = \frac{2}{3 \Delta t} \boldsymbol{U}_{n} + \frac{1}{\left | \Omega_{i} \right |_{n+1}} \left ( \mathcal{R} (\boldsymbol{U}_{n}) - \frac{2}{\Delta t} \left | \Omega_{i} \right |_{n} \boldsymbol{U}_{n} + \frac{1}{2 \Delta t}  \left | \Omega_{i} \right |_{n-1} \boldsymbol{U}_{n-1} \right )
\end{equation}
where the physical and fictitious time are expressed by $\Delta t$ and $\tau$ respectively. The convergence of every physical time in pseudo time results in the modified residual equates to $\mathcal{R}^{*}(\boldsymbol{U}_{n}) = 0 $ meaning it is equivalent to finding the state $\boldsymbol{U}_{n} = \boldsymbol{U}_{n+1}$.

\subsection{Harmonic Balance}
The high computational cost of time-accurate methods means that it is often desirable to use reduced-order models for simulating unsteady problems. Adopting the harmonic balance method is beneficial as it can be used for quasi-periodic flows dominated by a specific set of frequencies which need not be integral multiples of each other \cite{hall2002computation}. The subsequent harmonic balance method was implemented by Ref. \cite{rubino2018adjoint} in SU2 and will briefly be summarised.

In the harmonic balance method the time operator, $ \mathcal{D}_{t} $, is introduced and is approximated using spectral interpolation. Applying the spectral operator to the vector of conservative variables, $ \tilde{\boldsymbol{U}} $, which are now to be evaluated at $\mathcal{N}$ time instances, one obtains,

\begin{equation}
\mathcal{D}_{t} ( \boldsymbol{U} ) \approx  \mathcal{D}_{t} ( \tilde{\boldsymbol{U}} )
\end{equation}
and the harmonic balance operator $\mathcal{D}_{t}$ can eventually take the form,

\begin{equation}
\mathcal{D}_{t} ( \tilde{\boldsymbol{U}} ) = \boldsymbol{H} \tilde{\boldsymbol{U}}
\end{equation}
where $ \boldsymbol{H} $ is the spectral operator matrix. While now considering $ \tilde{\boldsymbol{U}} $ as the vector of conservative variables evaluated at $ \mathcal{N} $ time instances, Equation~\ref{chap4:eq:dual_time} can be rewritten for a single time instance as,

\begin{equation}
\label{chap4:eq:hb}
\left | \Omega_{i} \right |  \frac{\Delta \boldsymbol{U}_{n}^{q+1}}{\Delta \tau_{n}} + \boldsymbol{H} \tilde{\boldsymbol{U}} + \mathcal{R}(\boldsymbol{U}_{n}^{q+1}) = 0
\end{equation}			
where $q+1$ is the physical time step index and where $n$ is now the time instance. The expression can then be linearised and a semi-implicit approach can be used to solve for each time instance. An unsteady problem can then be characterized by $K$ frequencies. 
%
%
%

\subsection{Non-Conformal Boundary Interface Treatment}
\begin{figure}[hbht!]
	\centering
	\includegraphics[width=0.995\linewidth]{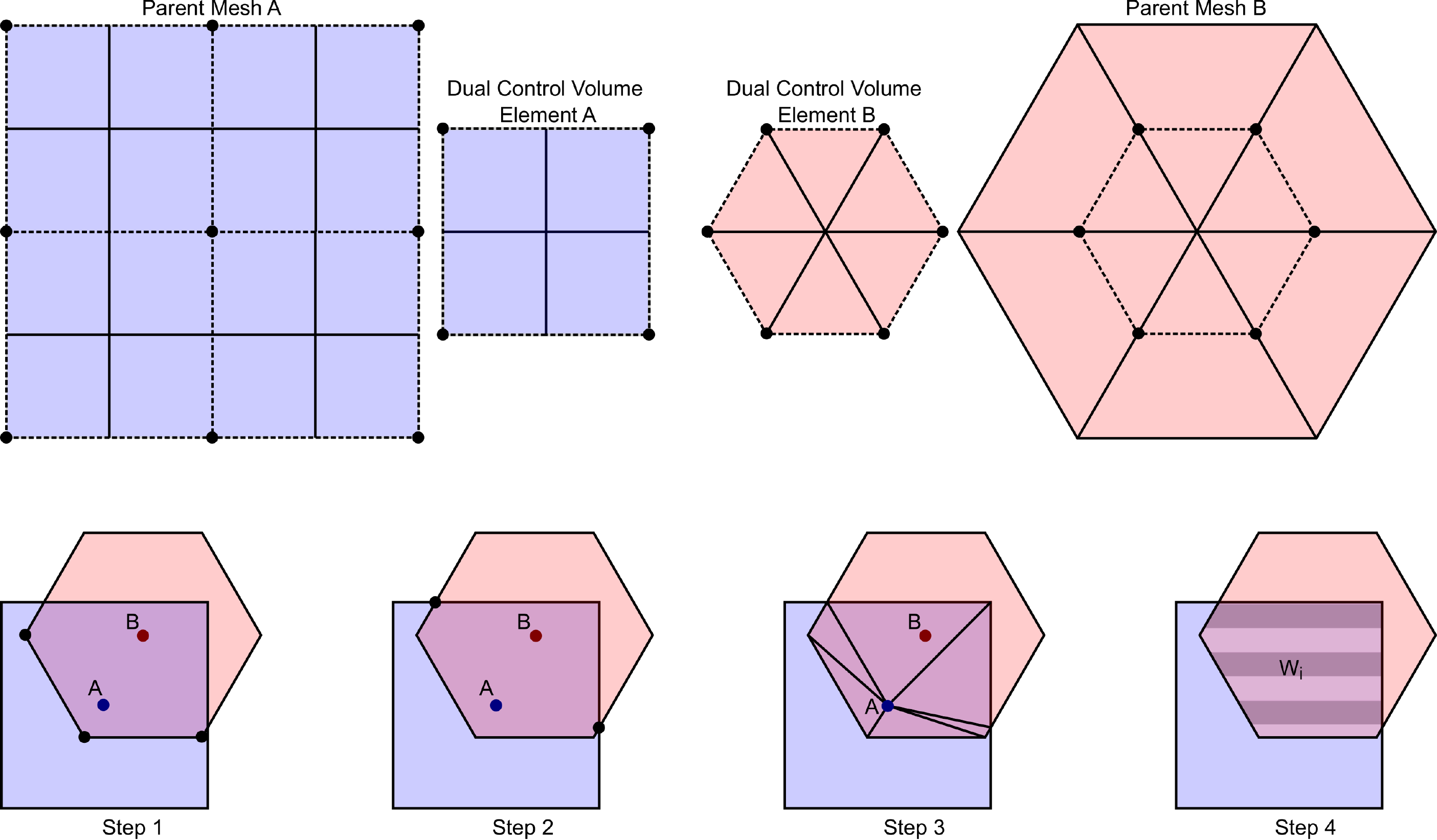}
	\caption[Supermesh schematic]{Schematic showing the steps for the construction of the supermesh.}
	\label{fig:supermesh}
\end{figure}

The numerical discretization of complex rotorcraft geometries is a non-trivial task. It requires the main and tail rotor blades to move in relative motion to the fuselage. To permit this, the computational domain should be split into separate sub-regions, thus introducing non-conformal boundary interface. Suitable treatment of each artificial boundary interface is achieved through the supermesh technique \cite{rinaldi2015flux} implemented by developers Gori et al. \cite{gori2017sliding} of SU2. The supermesh acts an additional auxiliary grid between two non-conformal interfaces and is of size $n - 1$ dimensional elements where $n$ is the dimension of the computational domain. The algorithm used to construct the supermesh of two generic parent mesh is shown in Fig.~\ref{fig:supermesh} and subsequently outlined:

\begin{enumerate}
	\item Identify the nodes of the control volume element $\mathcal{E}_B$ contained within $\mathcal{E}_A$.
	\item Determine the points of intersection between the edges of $\mathcal{E}_A$ and $\mathcal{E}_B$.
	\item Decompose the overlap region into triangles.
	\item Calculate the overlapping area of $\mathcal{E}_B$ inside of $\mathcal{E}_A$ and the weight $W_{i}$.
\end{enumerate}
this procedure is repeated for each of the overlapping neighbouring elements and its contribution to the numerical flux is assembled by $\sum^{N_{sf}}_{q=1}{W_{q}}$. Where the number of supermesh faces that are mapped to each control volume is denoted by $\mathcal{N}_{sf}$. The flux balance across the supermesh interface is obtained through incorporating the neighbouring cells either side of the boundary into Eq.~\ref{eq:FVM}, which can then be re-written in the semi-discrete form in the cell center of each control volume, $k$ as,

\begin{equation}
\frac{\partial \boldsymbol{U}_{k}}{\partial t} = 
-\frac{1}{\left | \Omega_{k} \right |} \left[
\sum_{p=1}^{\mathcal{N}_f} \boldsymbol{F} (\boldsymbol{U}_k, \, \boldsymbol{U}_p)\Delta S_{p}
+ 
\sum_{p=1}^{\mathcal{N}_{bf}} \Delta S_{p} 
\left(\sum_{q=1}^{\mathcal{N}_{sf}} W_{q} \boldsymbol{F} (\boldsymbol{U}_k, \, \boldsymbol{U}_q) \right) \right]
\end{equation}
where the number of internal faces and boundary faces of each control volume is represented by $\mathcal{N}_f$ and $\mathcal{N}_{bf}$. The area of the face is represented by $\Delta S$. The vector of fluxes are represented as $\boldsymbol{F}$. The conservative variables are once again represented by $\boldsymbol{U}$ and the subscripts $p$ and $q$ denote the neighbouring and boundary cells. 

\subsection{Main Rotor Blade Kinematics}
\label{sec:kinematics}
\begin{figure}[hbht!]
	\centering
	\includegraphics[width=0.99\linewidth]{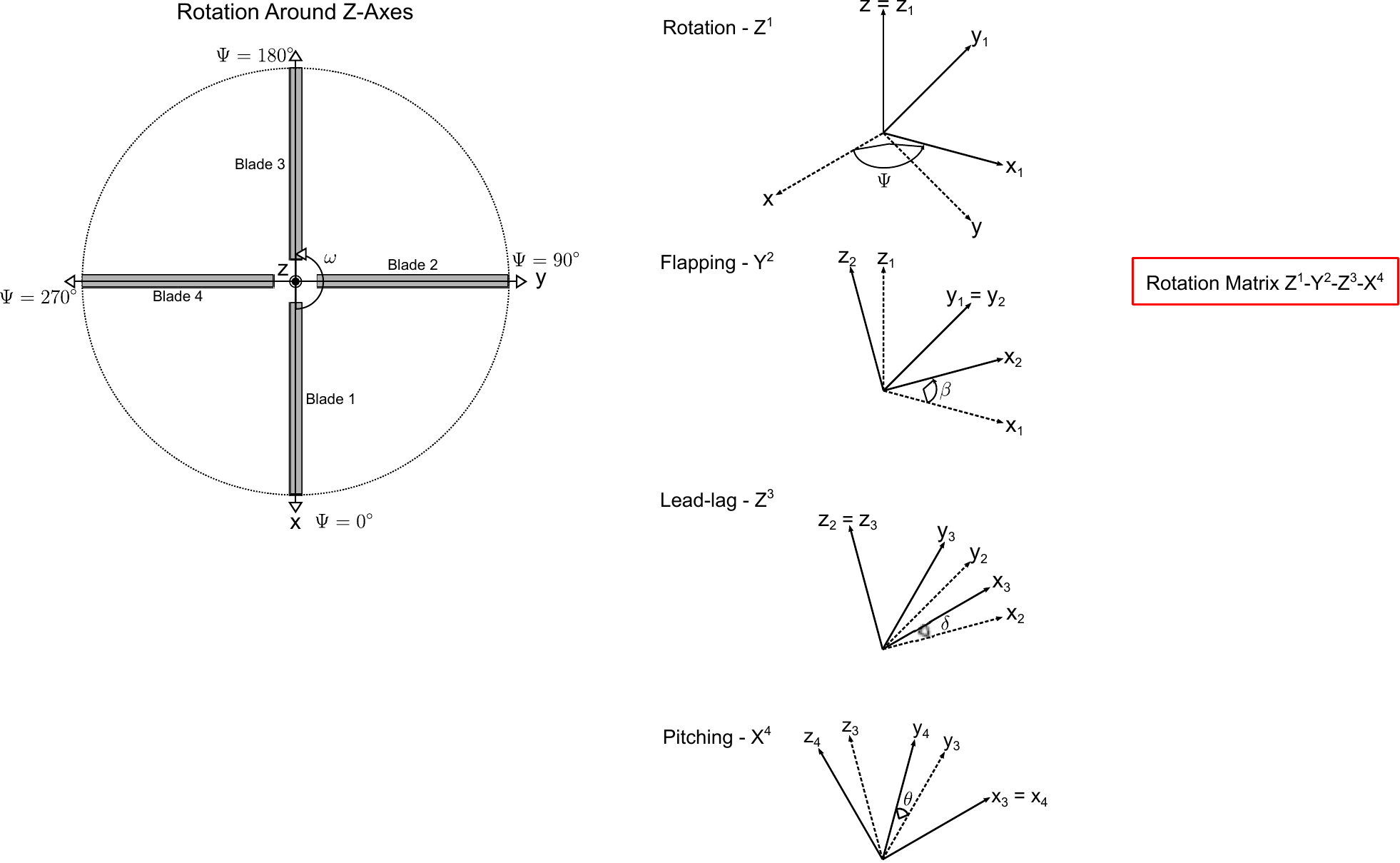}
	\caption[Main rotor schematic]{Schematic of the main rotor blade motion and respective rotations.}
	\label{fig:kinematics}
\end{figure}

The method in which conventional main-rotor/tail-rotor helicopters generate their vertical and propulsive forces, as well as moments to control the attitude and position of the helicopter in three-dimensional space, will now be introduced and its implementation within SU2 discussed. A schematic depicting the blade motion is shown in Fig.~\ref{fig:kinematics}. Unlike fixed-wing aircraft, the rotor blades alone must provide both the lifting forces and the control meaning the aerodynamics and dynamics of the blades are coupled so that it is hard to understand one without the other. During hover and axial flight, the velocity variation along the blade is azimuthally axisymmetric and depends solely upon the relationship of $\omega r$. In forward flight, however, this is not the case and a component of the free-stream velocity alters the blade velocity such that it now depends upon its azimuthal position, $\psi$ as given by,

\begin{equation}
M_{\textrm{n}}(\psi) = M_{\text{tip}}\frac{r}{R} + M_{\infty} \sin \psi = M_{\text{tip}} \left ( \frac{r}{R} + \mu \sin \psi \right )
\end{equation}
where $M_{\textrm{n}}$ represents the blade normal Mach number and $M_{\infty}$ is the freestream Mach number. The advance ratio, $\mu$, can be considered as the relationship of $M_{\infty}/M_{\text{tip}} $ . A repercussion of local variations in velocity is that each blade is now required to be individually controlled to eliminate rolling and pitching moments. With this, each blade is typically connected to a rotor head by a set of hinges which allows the blades to move independently. In reality, the rotor blades are also highly elastic in nature, however, for simplification this work considers the blades to be rigid with the ramifications of this being understood. Rigid blades are then able to move with respect to the hinge positions. Each of the hinges will now be introduced: a flap hinge allows the blade to move in the plane containing the blade and the shaft; a pitch hinge allows the blade to move around the quarter chord and spanwise axis; finally a lead-lag hinge allows the blade to move in-plane forward or backwards. The flapping $\beta$, lead-lag $\delta$, and pitching $\theta$ motion can then be described as a function of the azimuthal position $\psi$ of the blades as they rotate,

\begin{align}
\label{eq:blade_motion}
\begin{aligned}
\psi &= \omega t \\       
\beta (\psi) &= \beta_{0} -\beta_{1s}\sin(\psi) -\beta_{1c}\cos(\psi) -\beta_{2s}\sin(2\psi) -\beta_{2c}\cos(2\psi) - ... \\
\delta (\psi) &= \delta_{0} -\delta_{1s}\sin(\psi) -\delta_{1c}\cos(\psi) -\delta_{2s}\sin(2\psi) -\delta_{2c}\cos(2\psi) - ... \\
\theta (\psi) &= \theta_{0} -\theta_{1s}\sin(\psi) -\theta_{1c}\cos(\psi) -\theta_{2s}\sin(2\psi) -\theta_{2c}\cos(2\psi) - ...
\end{aligned}
\end{align}
where $\beta_{0}$ denotes the coning angle and $\theta_{0}$ denotes the collective pitch applied equally to all the blades. The lateral and longitudinal components of the motion are represented by subscripts $c$ and $s$ respectively. 

A rotation matrix then transforms the blade Cartesian coordinates, $\boldsymbol{x}$, from the fixed frame of reference into the hub frame of reference. The hub reference system is then transformed into each of the blade reference systems and the blade motion laws described in Equation~\ref{eq:blade_motion} are applied. If we now consider the definition of the reference system where; the rotation occurs around the $z$-axes, the flapping occurs around the $y$-axes, the lead-lag motion around the $z$-axes, and the pitching occurs around the $x$-axes as depicted in Fig.~\ref{fig:kinematics} the following transformation matrices can be introduced to prescribe the blade motion,

\begin{align}
\begin{aligned}
C_{\textrm{rot}} &= \begin{pmatrix}
\cos\psi & -\sin\psi & 0\\ 
\sin\psi & \cos\psi & 0\\ 
0 & 0 & 1
\end{pmatrix}, \\       C_{\textrm{lead-lag}} &= \begin{pmatrix}
\cos\delta & -\sin\delta  & 0\\ 
\sin\delta & \cos\delta & 0\\ 
0 & 0 & 1
\end{pmatrix},
\end{aligned}
&&
\begin{aligned}
C_{\textrm{flap}} &= \begin{pmatrix}
\cos\beta & 0 & -\sin\beta\\ 
0 & 1 & 0\\ 
\sin\beta & 0 & \cos\beta 
\end{pmatrix} \\       C_{\textrm{pitch}} &= \begin{pmatrix}
1 & 0 & 0\\ 
0 & \cos\theta & -\sin\theta \\ 
0 & \sin\theta & \cos\theta
\end{pmatrix}
\end{aligned}
\end{align}

The transformation from the fixed $\rightarrow$ hub reference system using $C_{\textrm{rot}}$ is implemented by a rigid rotation of the entire grid as the azimuthal rotation of the blades far exceeds the blade deflections. The rate of rotation is currently constrained and is maintained constant such that,

\begin{equation}
C_{\textrm{rot}} = \begin{pmatrix}
\cos\psi & -\sin\psi & 0\\ 
\sin\psi & \cos\psi  & 0\\ 
0 & 0 & 1
\end{pmatrix}
\end{equation}

The transformation from the hub $\rightarrow$ blade reference system can then be combined and applied to the blade surface. This is achieved through deforming the mesh at each physical time step and is required as each blade has its own independent motion. The subsequent transformation matrix, $C_{\textrm{flap-leadlag-pitch}}$ can then be introduced,

\begin{equation}
C_{\textrm{flap-leadlag-pitch}} = \begin{pmatrix}
\cos\beta \cos\delta & \sin\beta \sin\theta - \cos\beta \cos\theta \sin\delta & \cos\theta \sin\beta + \cos\beta \sin\delta \sin\theta\\ 
\sin\delta & \cos\delta \cos\theta & -\cos\delta \sin\theta\\ 
-\cos\delta \sin\beta & \cos\beta \sin\theta + \cos\theta \sin\beta \sin\delta& \cos\beta \cos\theta - \sin\beta \sin\delta \sin\theta
\end{pmatrix}
\end{equation}

The motion can then be applied to the coordinates, $\boldsymbol{x} = \{{x, y, z}\}^{T}$ at an arbitrary point, $i$ on the blade as, 

\begin{equation}
\boldsymbol{x}_{i} = C_{\textrm{rot}}C_{\textrm{flap-leadlag-pitch}}(\boldsymbol{x}_i - \boldsymbol{x}_{\textrm{hinge}})
\end{equation}	
where $\boldsymbol{x}_{\textrm{hinge}}$ is the position of the hinge about which the flapping, lead-lag and pitching occurs. Furthermore the grid velocity is updated using a second order finite difference scheme,
\begin{equation}
\boldsymbol{u}_{\Omega_i} = \frac{3 \boldsymbol{x}_{i}^{n + 1} - 4 \boldsymbol{x}_{i}^{n} + \boldsymbol{x}_{i}^{n -1} }{2 \Delta t}
\end{equation}

\subsection{Radial Basis Function Mesh Deformation}
To allow the rotor blades to move independently requires the mesh to deform at each physical time step. Rendall and Allen, to the authors' knowledge, were the first and only group to demonstrated how Radial Basis Functions (RBF) could be used to account for the blade motion \cite{rendall2010parallel}. The outlook of using RBF mesh deformation for rotor blade motion was promising however one major drawback was the high computational cost for large meshes. Despite this, RBF mesh deformation methods are robust and preserve high-quality mesh even during large deformations. They also have the unique property that they do not require the grid connectivity meaning that even for three-dimensional problems they are relatively simple to implement. 

The term radial basis function refers to a series of functions whose values depends on their distance to a supporting position. In the most general of forms, radial basis functions can be written as, $ \phi (\textbf{r}, \textbf{r}_{i}) = \phi\left( \| \textbf{r} - \textbf{r}_{i} \|\right) $, where the distance corresponds to the radial basis centre, $ \textbf{r}_{i} $. This distance is frequently taken as the Euclidean distance, meaning it becomes the spatial distance between two nodes. 

An interpolation function, $f(\textbf{r})$ can be introduced as a method for describing the displacement of a set of nodes in space and can be approximated by a weighted sum of basis functions. However, the interpolation relies on the weight coefficients of the basis points, $\alpha$. The interpolation takes the form

\begin{equation}
\label{eq:2}
f(\textbf{r}) = \sum_{i=1}^{N} \alpha_{i} \phi \left ( \left \| \textbf{r} - \textbf{r}_{i} \right \| \right )
\end{equation}
The primary expense of RBF mesh deformation is associated to the solving the linear system to obtain the weight coefficients. To address the concerns of the computation cost of RBF for rotor blades, multi-level greedy surface point selection algorithms \cite{wang2015improved} and volume point reduction methods \cite{xie2017efficient} for large scale problems are introduced which greatly improve the efficiency. These numerically enhanced methods were implemented by Ref.~\cite{morelli2020radial} in SU2 and are now used within this work.

\section{Results}
\label{chap4:sec:validation}

\subsection{Hover}
The benchmark Caradonna and Tung experimental tests of a model rotor in hover \cite{caradonna1981experimental} are used for hover performance validation. The experiment was carried out in the U.S. Army Aeromechanics Laboratory's hover test facility which included special ducting designed to eliminate room recirculation. The rotor employed two manually adjustable cantilever-mounted blades attached to a drive shaft at the top of a column test stand. The blades were untwisted and untapered with an aspect ratio of 6. The profile of the blades used a symmetric NACA0012 airfoil. The rotor radius was $1.143$ m and the blade chord was $0.191$ m. The data recorded during the tests consisted of pressure measurements along the chord of the blade at various radial positions. The tests were conducted at a variety of collective pitch angles and tip Mach numbers. The conditions used for the validation of the blade loading predictions are detailed in Table~\ref{chap4:tab:hover}.  

\begin{table}[htb!]
	\caption[Hover test conditions taken from the Caradonna-Tung experiment]{Hover test conditions taken from the Caradonna-Tung experiment \cite{caradonna1981experimental}.}
	\centering
	\label{chap4:tab:hover}
	\begin{tabular}{lll} 		
		\hline \hline
		\begin{tabular}[c]{@{}l@{}}Collective Pitch\\ ($^{\circ}$)\end{tabular} &
		\begin{tabular}[c]{@{}l@{}}Rotational Speed\\ ($\textrm{RPM}$)\end{tabular} &
		\begin{tabular}[c]{@{}l@{}}Tip Mach No.\\ ( )\end{tabular} \\ \hline 
		8 &
		1250 &
		0.439 \\ \hline \hline
	\end{tabular}
\end{table}

\begin{figure}[hbt!]
	\centering
	\includegraphics[width=0.7\linewidth]{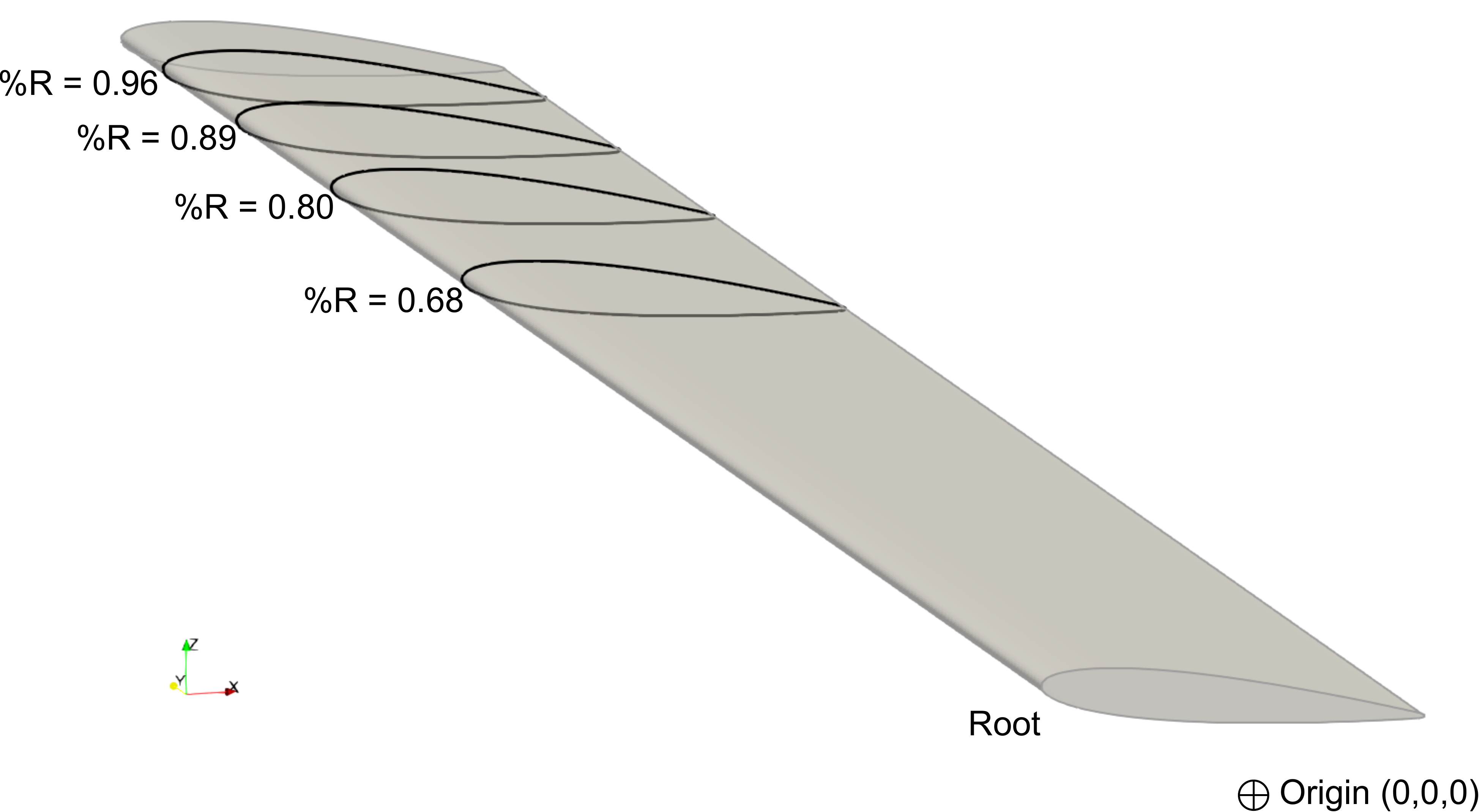}
	\caption[Blade static pressure measurements]{Blade stations where the static pressure measurements were taken during Ref.~\cite{caradonna1981experimental}. Select radial positions  $r/R = 0.68\%$, $r/R = 0.80\%$, $r/R = 0.89\%$ and $r/R = 0.96\%$ are used for the assessment of numerical predictions.}
	\label{fig.x:CT_measurementLocations}
\end{figure}

The problem is transformed into a rotating frame of reference for computational efficiency and the RANS equations are solved to first assess the influence of the radial position on the surface pressure distribution. Computed and experimental surface pressure distributions are compared at select radial positions $r/R = 0.68\%$, $r/R = 0.80\%$, $r/R = 0.89\%$ and $r/R = 0.96\%$ as represented by Fig.~\ref{fig.x:CT_measurementLocations}. 

The standard one-equation SA turbulence model is used with the RANS equations. The SA turbulence variable is convected using a first-order scalar upwind method. The convective fluxes are computed using the Roe scheme and second-order accuracy is achieved using the MUSCL scheme. Spatial gradients for the viscous fluxes and second order reconstruction are approximated using the Green-Guass method. The flow solution was considered converged when there was a reduction of 6 orders of magnitude on the density residual. 

\begin{figure}[!htb]
	\subfloat[$r/R = 0.68\%$.
	\label{fig.3:subfig-1:CT_surfacePressure}]{%
		\includegraphics[width=0.495\textwidth]{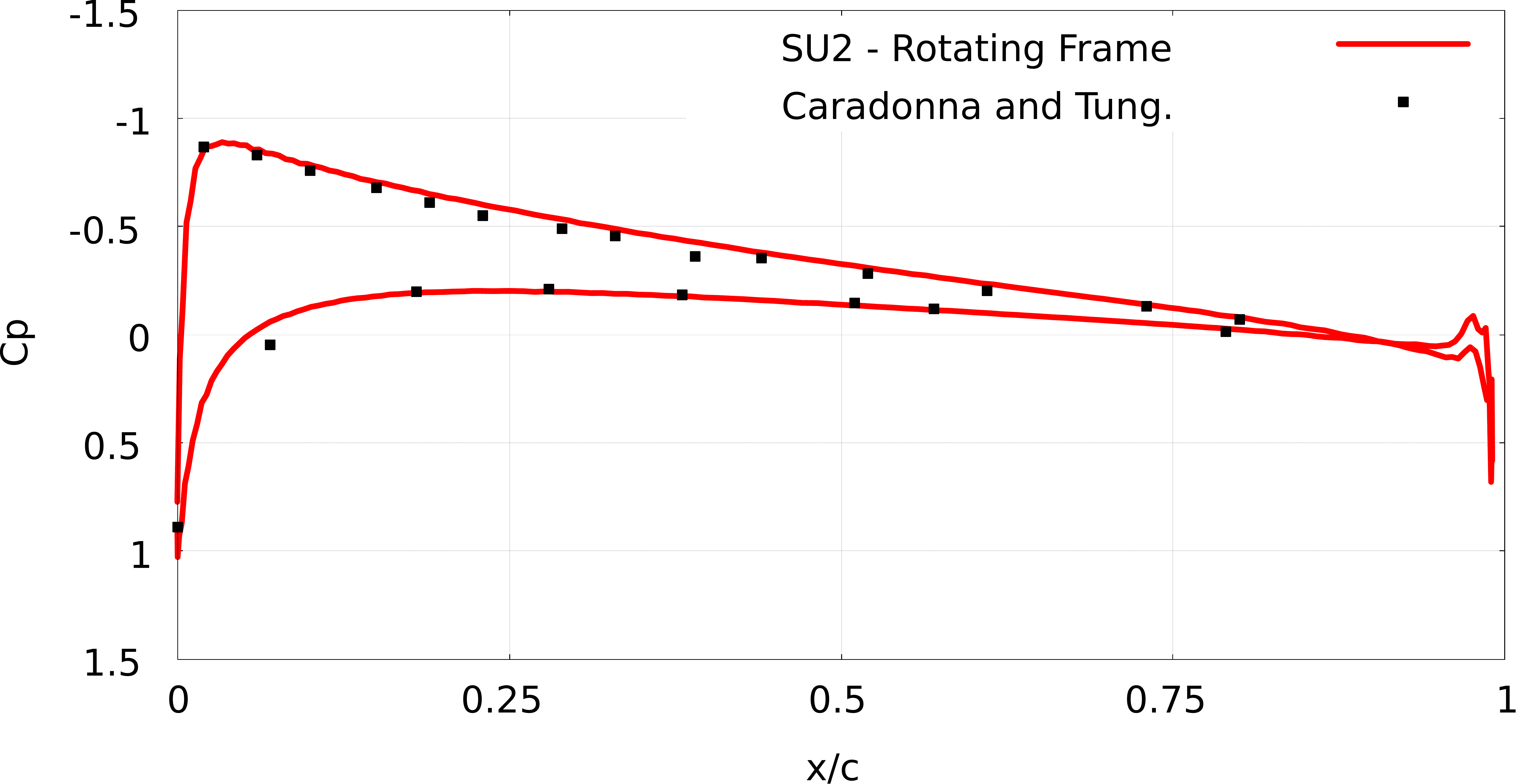}
	}
	\hfill
	\subfloat[$r/R = 0.80\%$.
	\label{fig.3:subfig-2:CT_surfacePressure}]{%
		\includegraphics[width=0.495\textwidth]{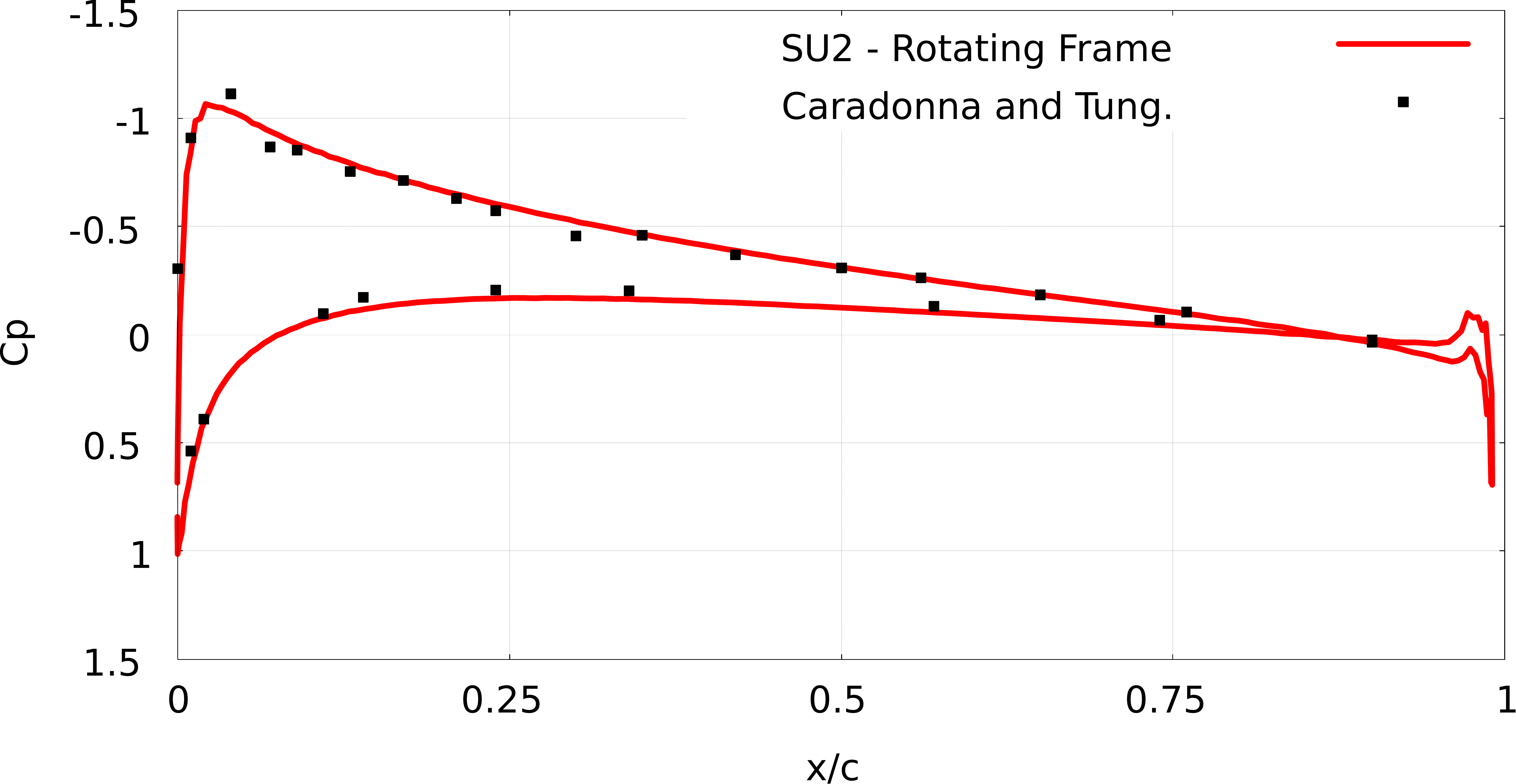}
	}\\
	\subfloat[$r/R = 0.89\%$.
	\label{fig.3:subfig-3:CT_surfacePressure}]{%
		\includegraphics[width=0.495\textwidth]{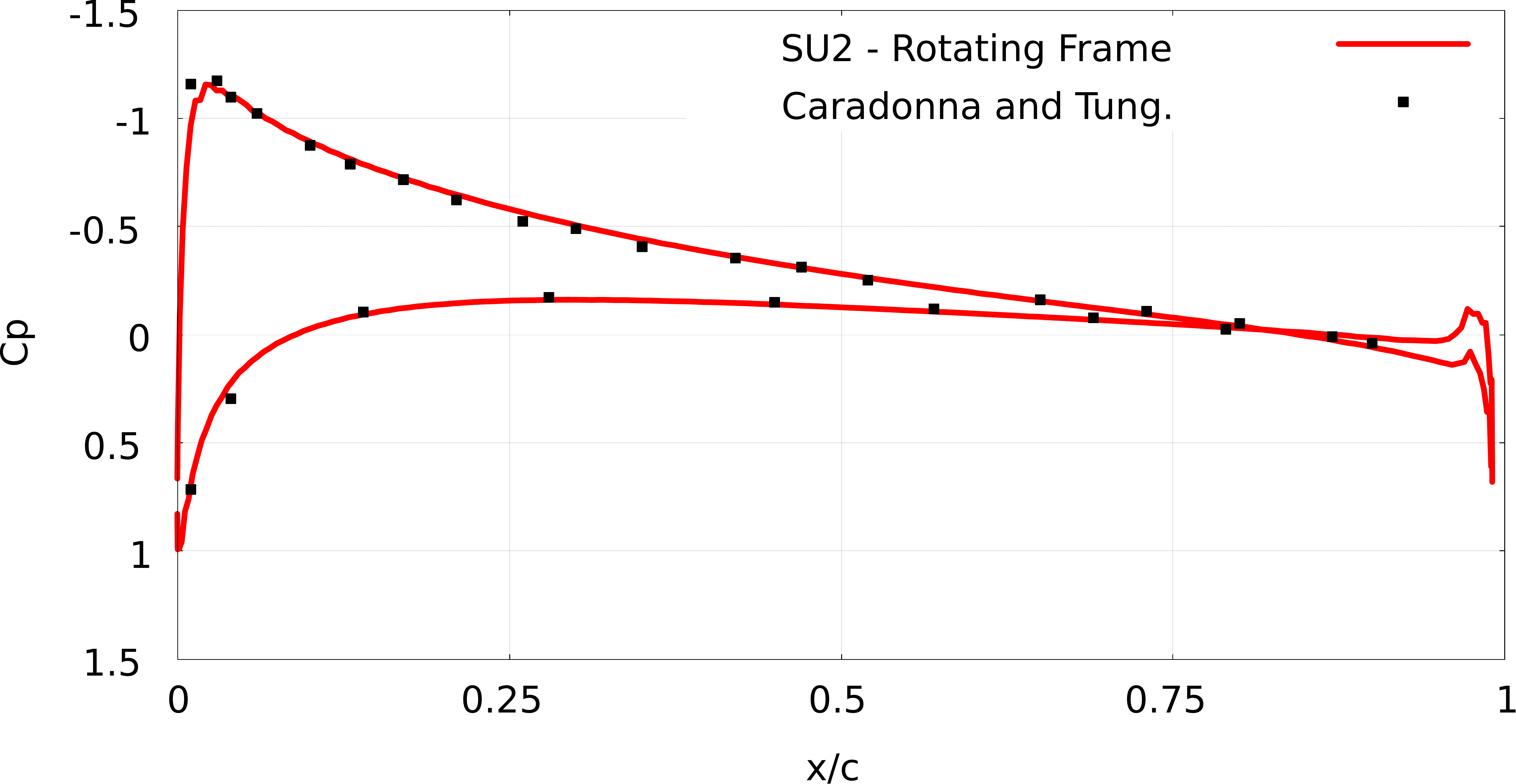}
	}
	\hfill
	\subfloat[$r/R = 0.96\%$.
	\label{fig.3:subfig-4:CT_surfacePressure}]{%
		\includegraphics[width=0.495\textwidth]{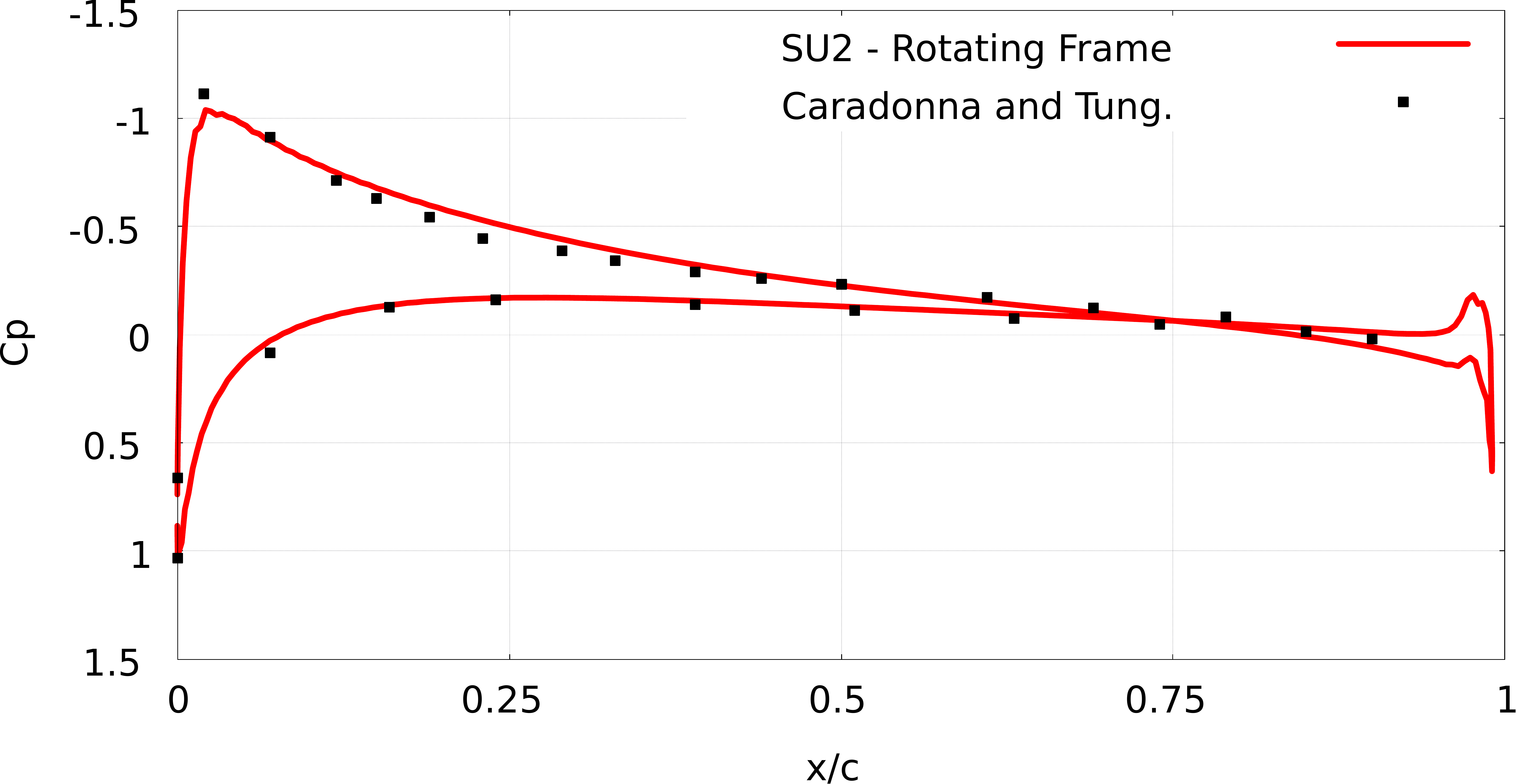}
	}
	\caption[$C_p$ distributions at various radial positions during hover]{Pressure distributions at different radial positions computed using the rotating frame method during hover. Predictions compared against measured data taken from Caradonna-Tung \cite{caradonna1981experimental}.}
	\label{fig.3:CT_surfacePressure}
\end{figure}

The predicted results in Fig.~\ref{fig.3:CT_surfacePressure} are in close agreement with the measured data at all radial positions. Both the upper and lower surface pressure profiles follow the measured data. The suction peak on the upper surface is captured well. There are however discrepancies at the trailing edge and this is expected to be caused by modelling the trailing edge as flat which was done to aid with the mesh generation process and convergence of the solver. The results reaffirm the calculations from Ref.~\cite{palacios2014stanford}. 

\begin{figure}[!htb]
	\subfloat[Grid sensitivity.
	\label{fig.4:subfig-1:CT_meshRes}]{%
		\includegraphics[width=0.495\textwidth]{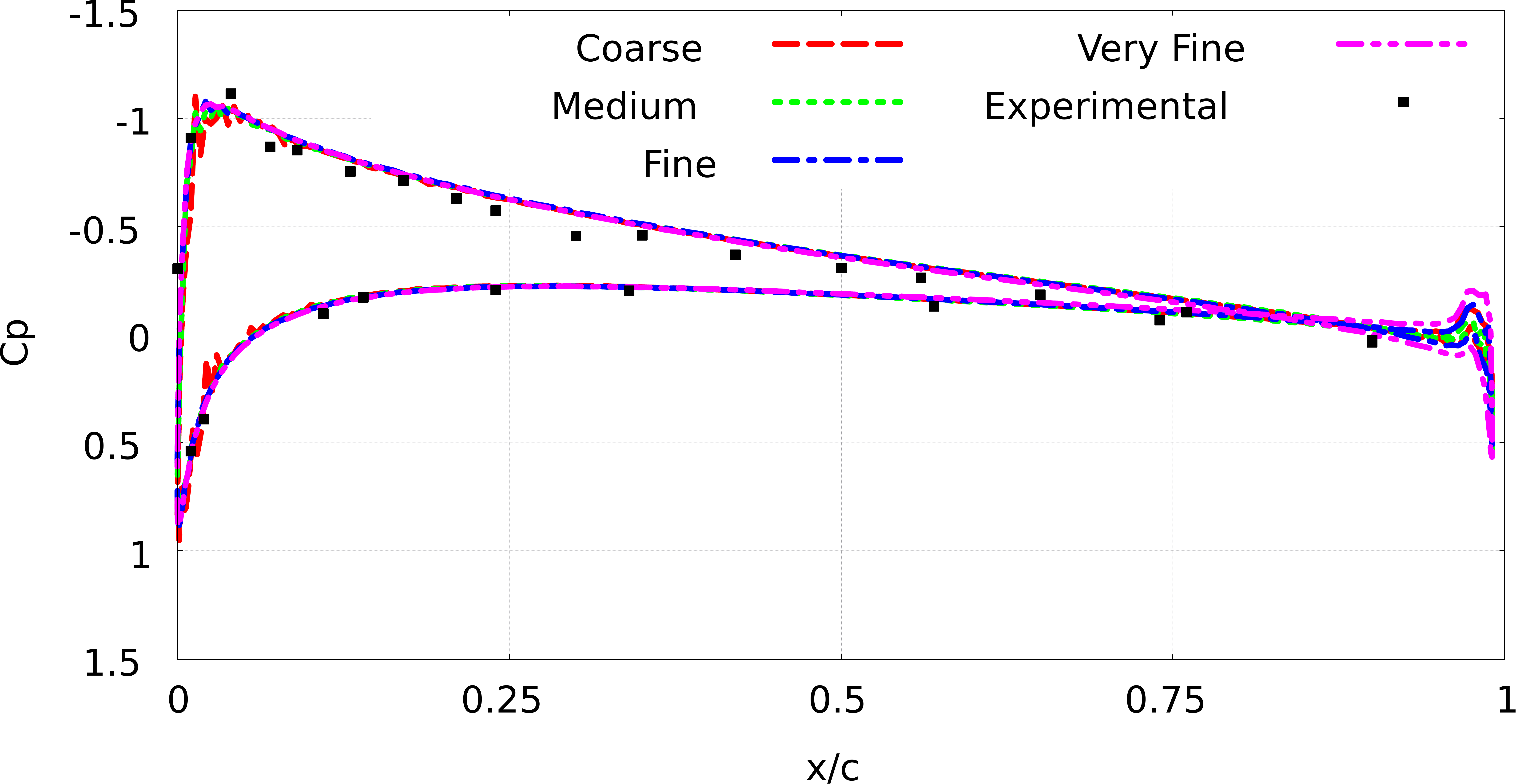}
	}
	\hfill
	\subfloat[Region around the leading-edge.
	\label{fig.4:subfig-1:CT_meshRes2}]{%
		\includegraphics[width=0.495\textwidth]{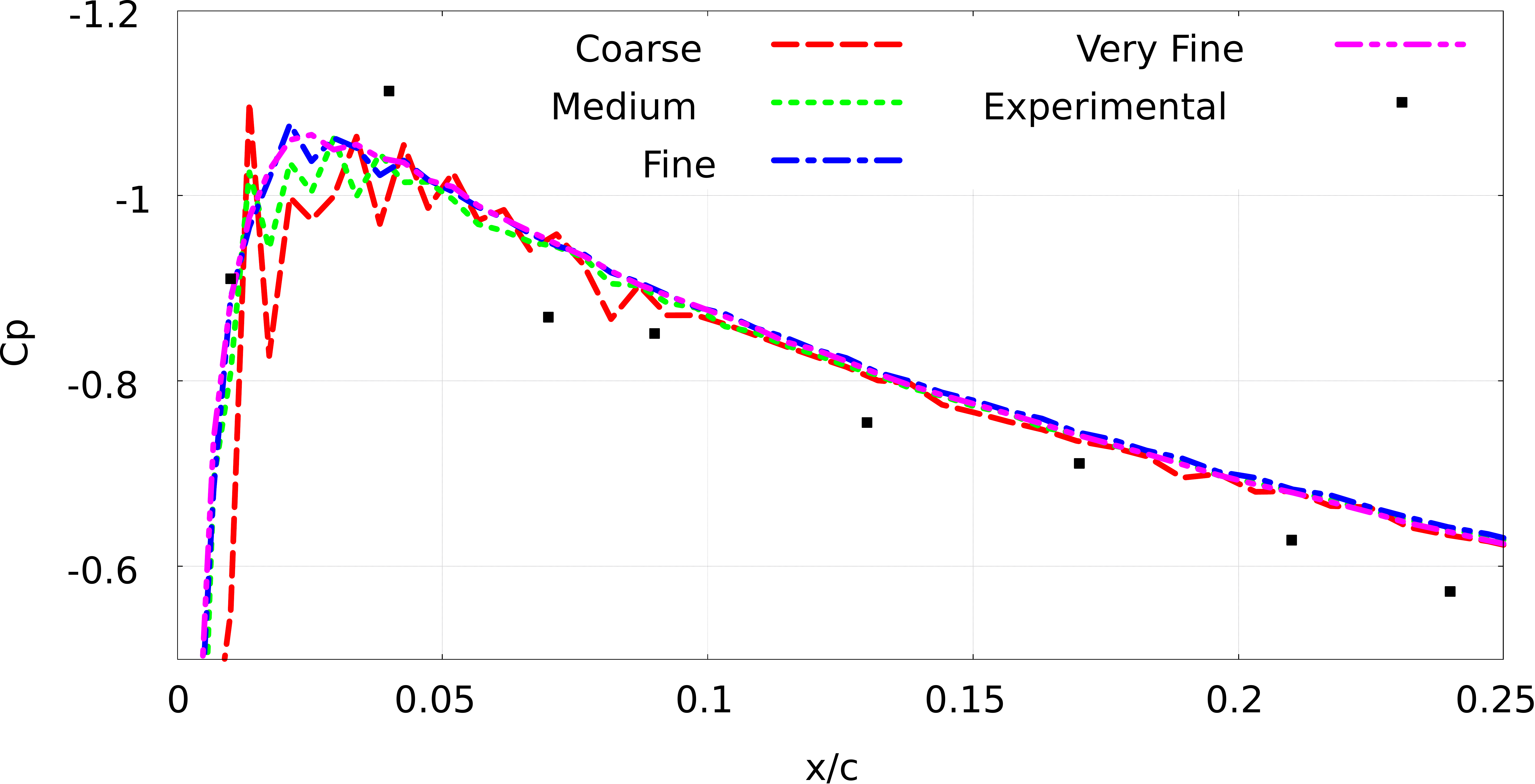}
	}
	\\
	\subfloat[Rotor revolutions.
	\label{fig.4:subfig-2:CT_nRevs}]{%
		\includegraphics[width=0.495\textwidth]{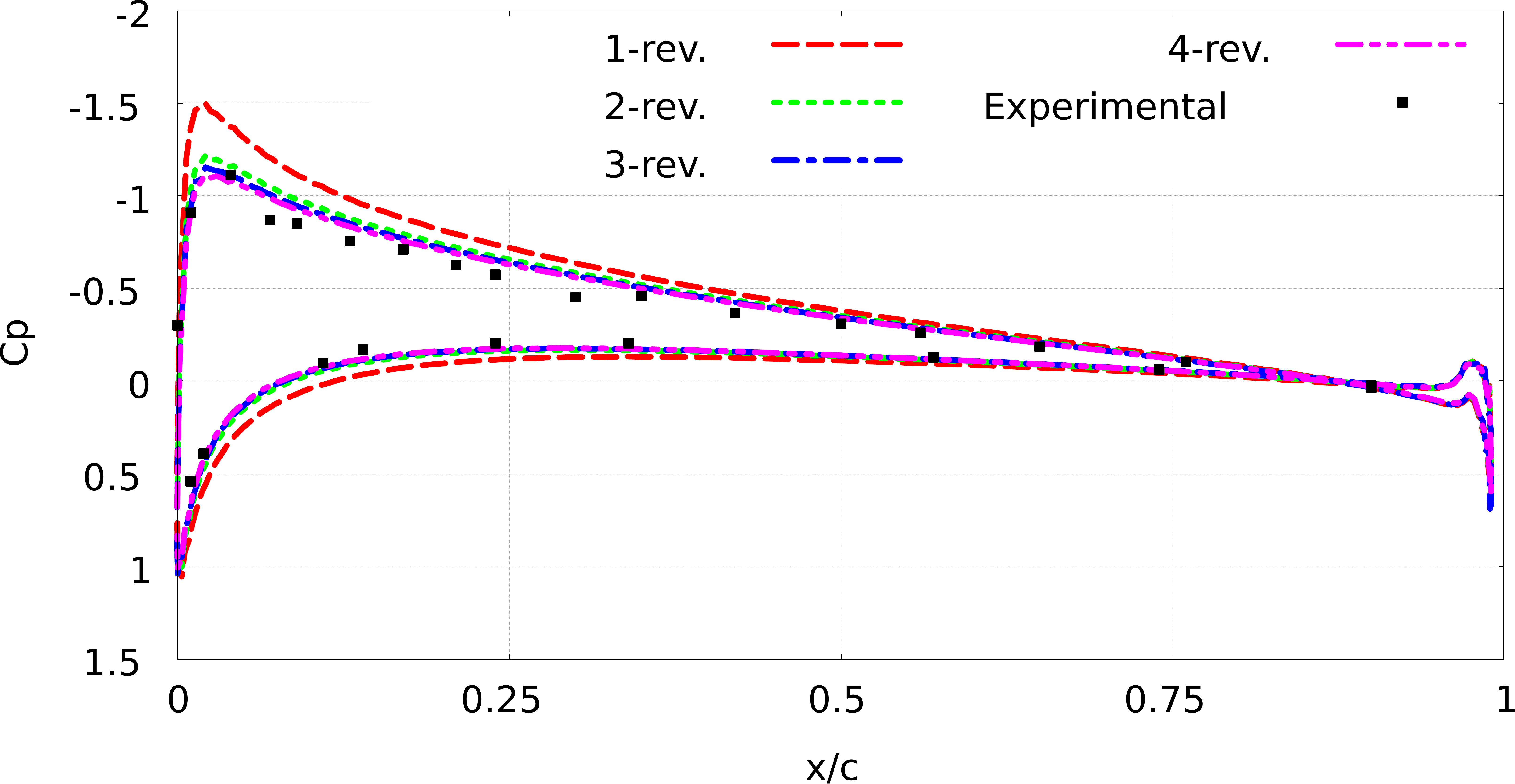}
	}
	\hfill
	\subfloat[Region around the leading-edge.
	\label{fig.4:subfig-2:CT_nRevs2}]{%
		\includegraphics[width=0.495\textwidth]{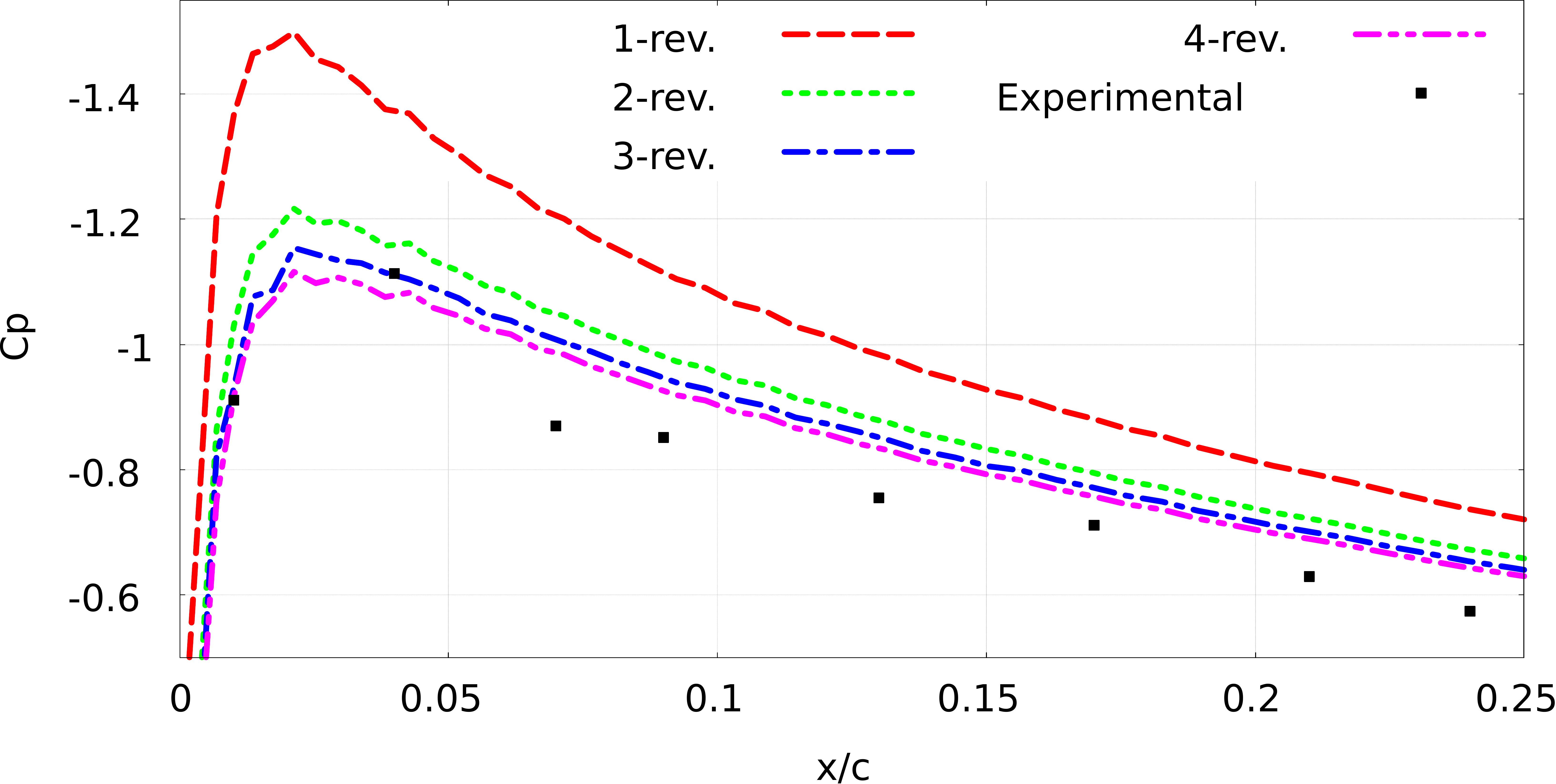}
	}
	\\
	\subfloat[Steady-state vs. time accurate approaches.
	\label{fig.4:subfig-3:CT_moveType}]{%
		\includegraphics[width=0.495\textwidth]{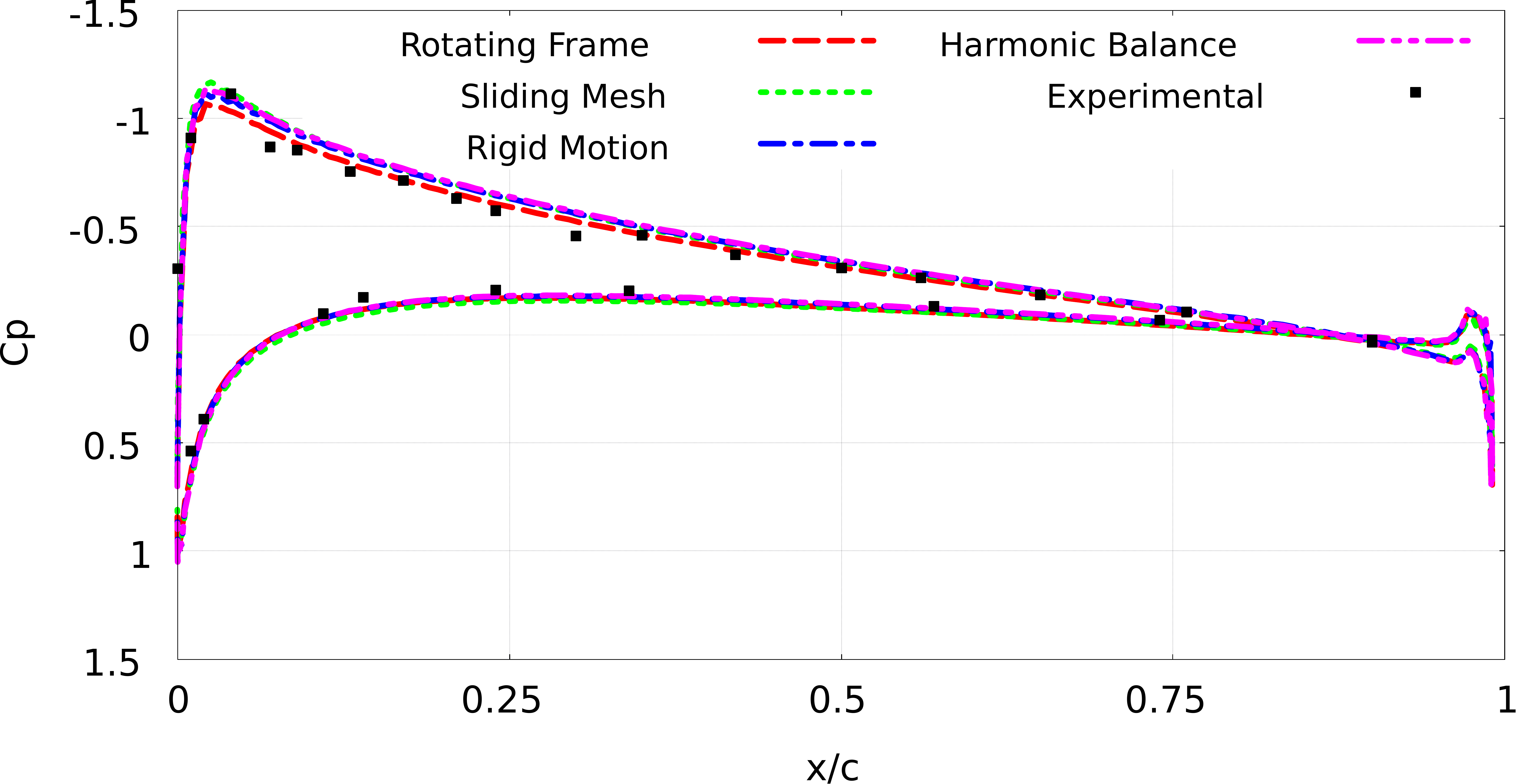}
	}
	\hfill
	\subfloat[Region around the leading-edge.
	\label{fig.4:subfig-3:CT_moveType2}]{%
		\includegraphics[width=0.495\textwidth]{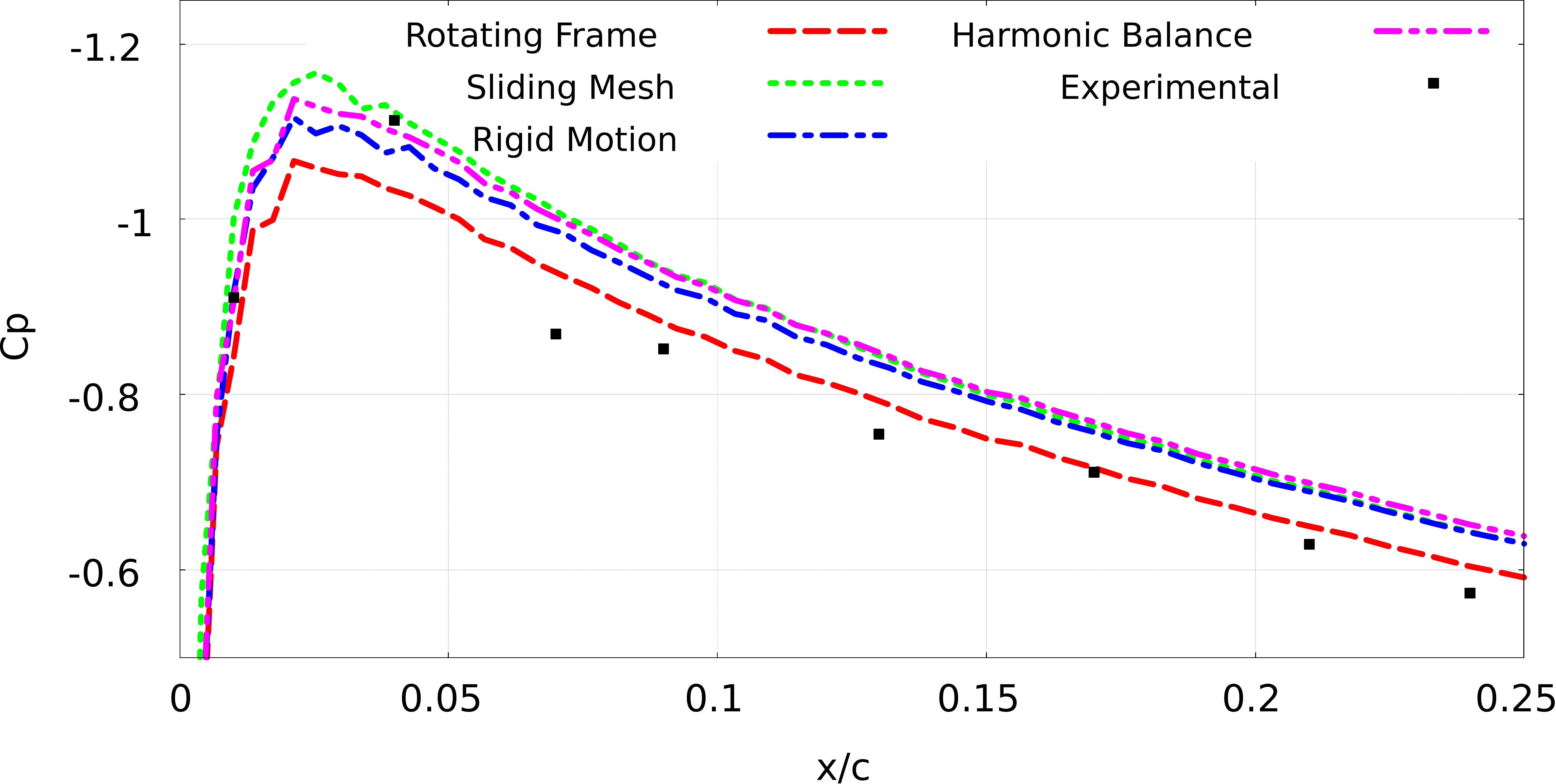}
	}
	\\
	\subfloat[Inviscid, laminar and turbulent-SA.
	\label{fig.4:subfig-4:CT_turbModel}]{%
		\includegraphics[width=0.495\textwidth]{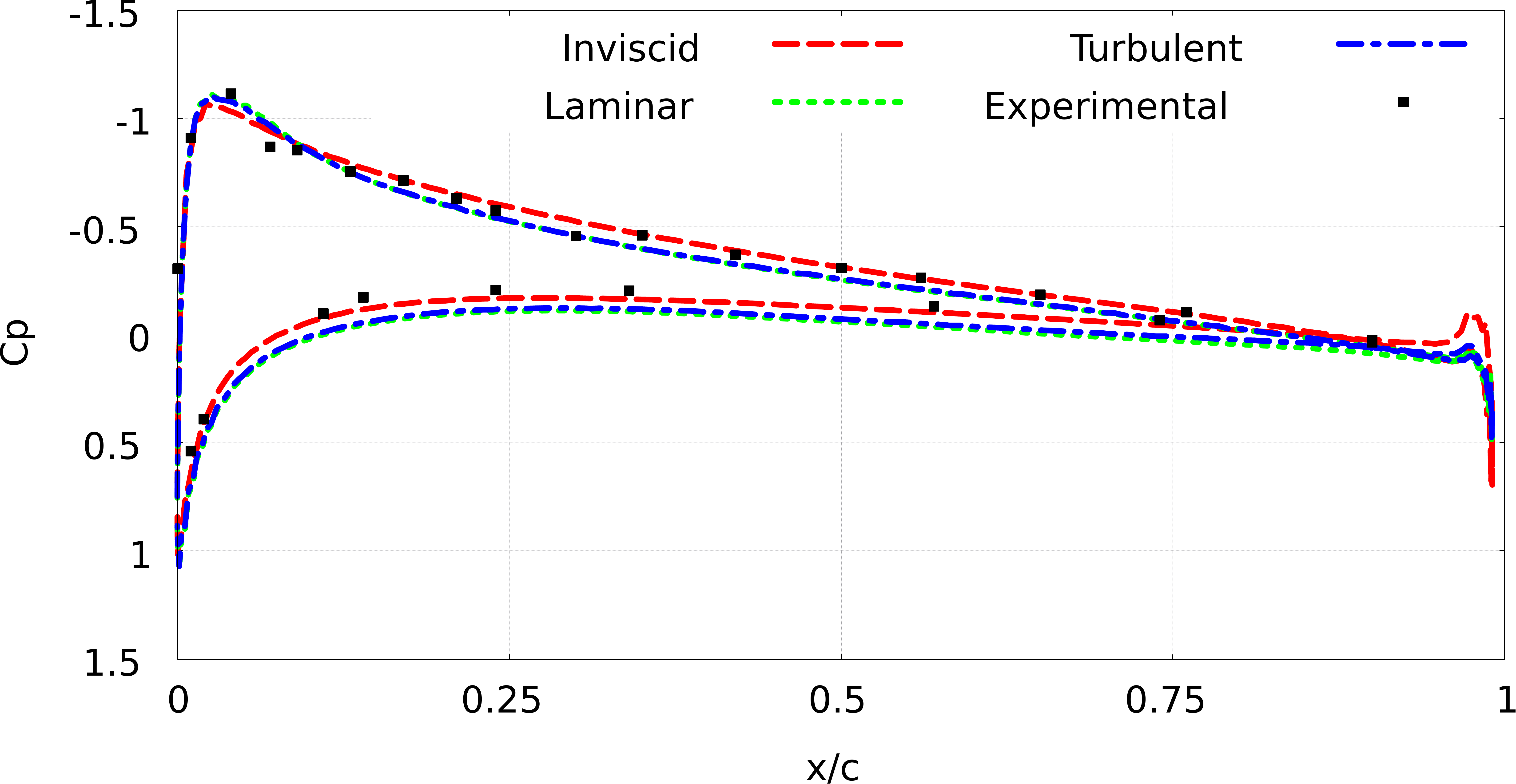}
	}
	\hfill
	\subfloat[Region around the leading-edge.
	\label{fig.4:subfig-4:CT_turbModel2}]{%
		\includegraphics[width=0.495\textwidth]{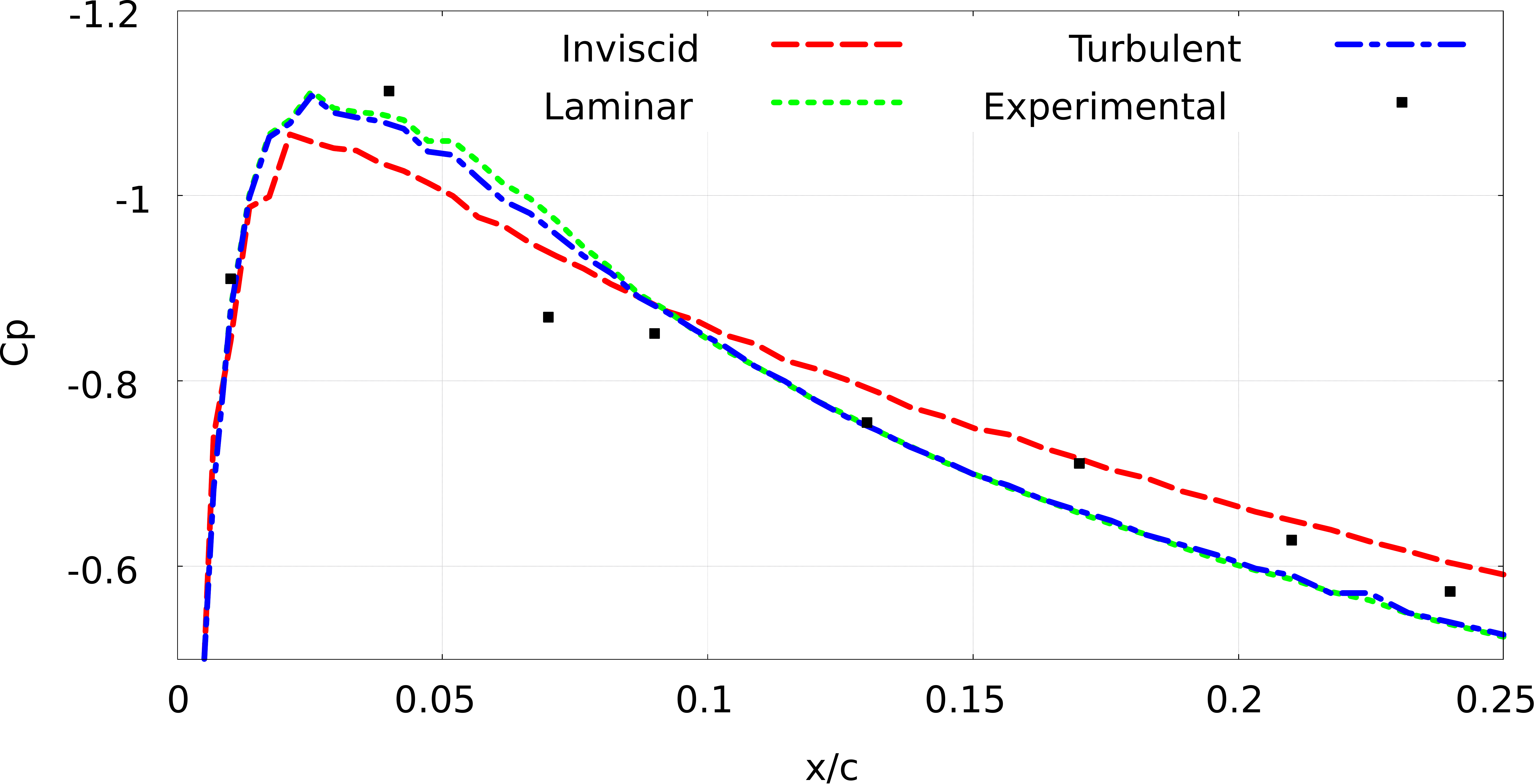}
	}
	\caption[Influential parameters on the performance prediction during hover]{Comparison of different influential parameters on the performance prediction during hover. Pressure distributions shown at $r/R = 0.80\%$ and compared against measured data taken from Caradonna-Tung \cite{caradonna1981experimental}.}
	\label{fig.4:CT_surfacePressure}
\end{figure}

An assessment of the most influential parameters during hover simulations are displayed in Fig.~\ref{fig.4:CT_surfacePressure}. Each of the parameters are compared against the measured data at $r/R = 0.80\%$. The spatial mesh resolution is assessed in Fig.~\ref{fig.4:subfig-1:CT_meshRes}. Four levels of mesh refinement are used ranging from coarse $\rightarrow$ very fine. The `coarse' mesh has 1.56 million volume elements and 22 thousand surface elements. The `medium' mesh has 2.62 million volume elements and 41 thousand surface elements. The `fine' mesh has 5.19 million volume elements and 80 thousand surface elements. The `very fine' mesh has 11.09 million volume elements and 157 thousand surface elements. The coarse mesh exhibits instabilities in the pressure profile at the leading-edge. The medium mesh reduces this and shows a smoother profile however the suction peak remains slightly below the measured data. The fine mesh improves the suction peak compared to the medium mesh. There appear to be no discernible differences between the fine and very fine mesh suggesting the mesh has reached convergence. 

It is frequently being desirable to model the rotor as a time-accurate problem where the blades are physically moving. In this scenario, the number of revolutions the blades have completed when the flow is initialized from the freestream is important and is shown in Fig.~\ref{fig.4:subfig-2:CT_nRevs}. It shows that the time for convergence of the solution is at least 3 full rotor revolutions.

The influence of the choice of steady or time-accurate modelling on the pressure profiles is shown in Fig.~\ref{fig.4:subfig-3:CT_moveType}. All methods show good agreement with the measured data. The rotating frame method shows the closest agreement. There is a very slight difference on the upper surface for the time-accurate sliding mesh, rigid motion and harmonic balance methods. Furthermore, all of the time-accurate methods appear overlapping.  

Finally the flow regime is evaluated in Fig.~\ref{fig.4:subfig-4:CT_turbModel}. It depicts the influence of inviscid, laminar and turbulent flow on the pressure distribution. As excepted there appears to be significant difference between inviscid flow and laminar and turbulent flow regimes.

\begin{figure}[!htb]
	\subfloat[Unsteady rigid motion at $10^{\circ}$ azimuth increments.
	\label{fig.5:subfig-1:CT_surfPressure}]{%
		\includegraphics[width=0.999\textwidth]{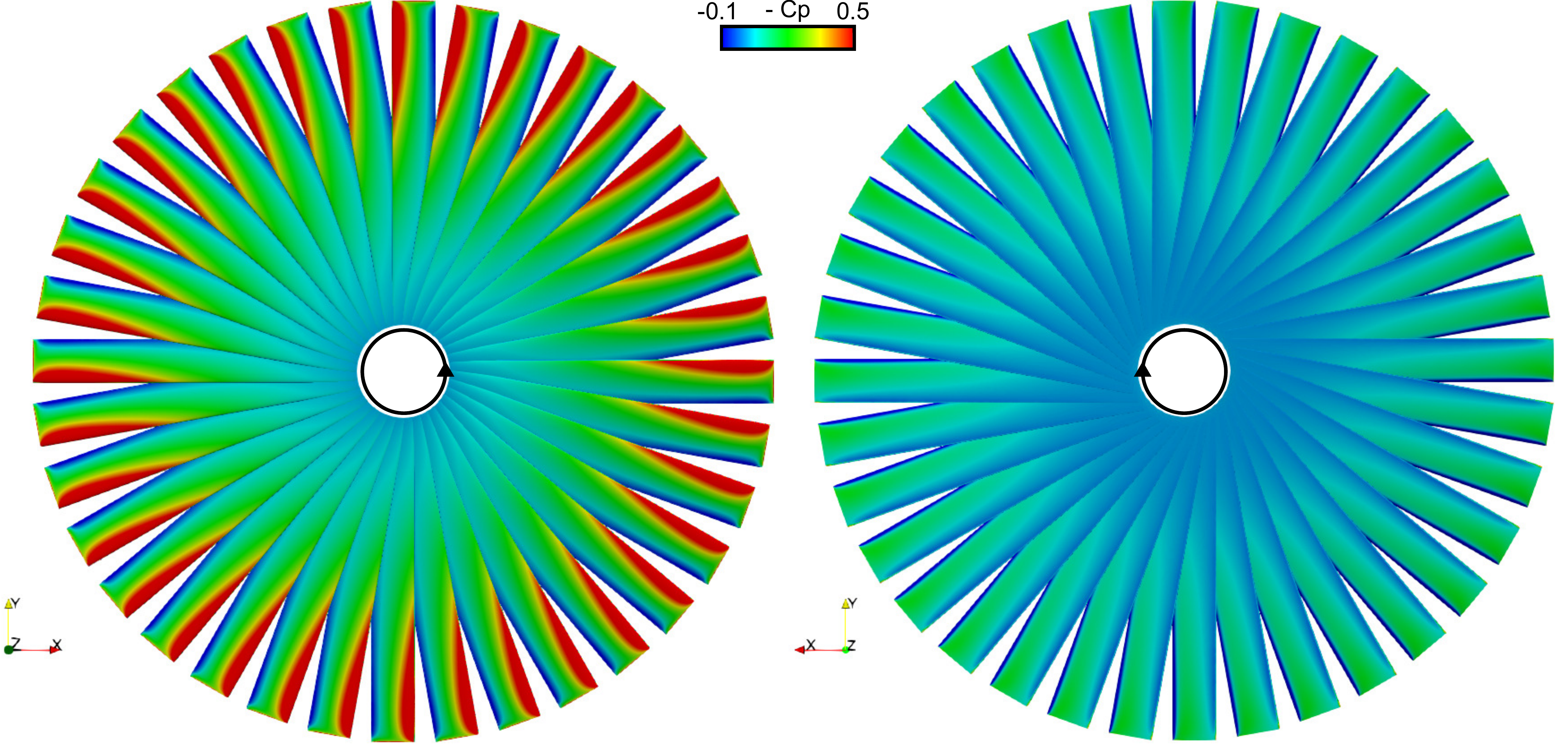}
	}
	\\
	\subfloat[Harmonic balance at 3 time instances with the input frequencies of $ \boldsymbol{\omega} = ( 0, \, \pm \omega_{1} )  $.
	\label{fig.5:subfig-2:CT_surfPressure}]{%
		\includegraphics[width=0.999\textwidth]{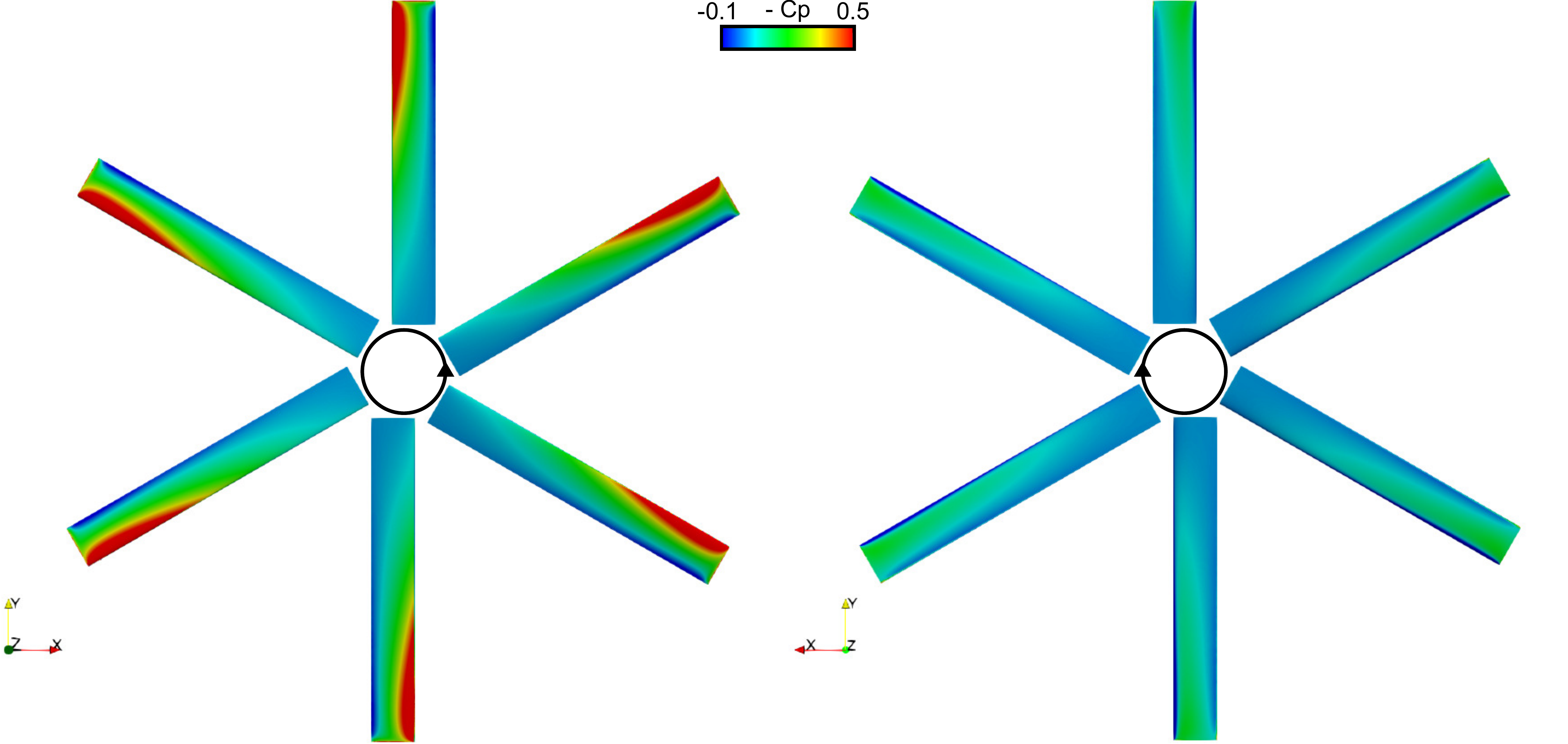}
	}
	\caption[$C_p$ distribution on the upper and lower blade surface during hover]{Contour maps of the pressure coefficient on the upper and lower blade surface during hover. Comparing the reduced order harmonic balance method against a fully time-accurate method.}
	\label{fig.5:CT_surfPressure}
\end{figure}

A contour map of the pressure coefficient on both the suction and pressure sides of the blades for time-accurate problems is shown in Fig.~\ref{fig.5:CT_surfPressure}. The approach using unsteady rigid motion is displayed at $10^{\circ}$ azimuth increments in Fig.~\ref{fig.5:subfig-1:CT_surfPressure}. The approach using the harmonic balance is displayed at 3 time instances with the input frequencies of $ \boldsymbol{\omega} = ( 0, \, \pm \omega_{1} ) $ in Fig.~\ref{fig.5:subfig-2:CT_surfPressure}. Both sets of results show that the pressure coefficient during hover is axisymmetric about the out-of-plane axis.

\begin{figure}[hbt!]
	\centering
	\includegraphics[width=0.7\linewidth]{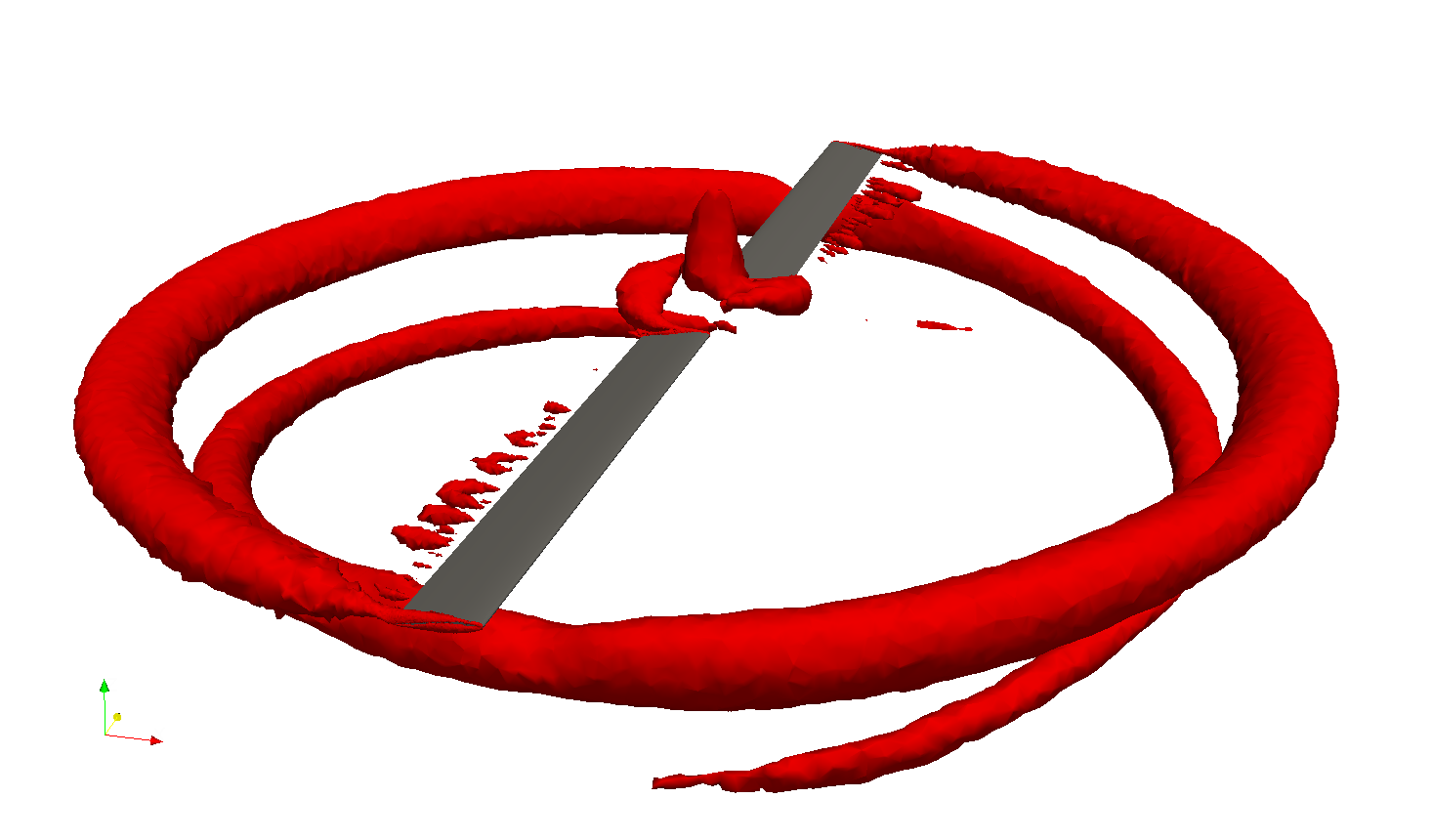}
	\caption[Q-criterion iso-surface during hover]{Iso-surface of the Q-criterion visualizing the near-field wake and blade tip vortices as they are convected downstream in the wake of the rotor during hover.}
	\label{fig6:CT_qCrit}
\end{figure}

The iso-surface of the Q criterion visualizing the near-field wake and blade tip vortices during hover are displayed in Fig.~\ref{fig6:CT_qCrit}. The wake exhibits no interaction effects due to the preceding blade tip vortices passing above previous blade tip vortices and due to the isolated blades being modelled without a hub or test stand.

\subsection{Forward Flight}

The extensive report and data survey from Cross and Watts on tip aerodynamics and acoustics \cite{cross1988tip} is used for forward flight performance validation. The report describes the Tip Aerodynamic and Acoustics Tests (TAAT) carried out at the NASA Ames Research Center. The TAAT used a highly instrumented AH-1G cobra helicopter and measured the rotor airloads at multiple radial locations. The AH-1G is a two-person single-engine helicopter. It has a two-bladed teetering rotor. The blades are untapered and have a linear twist of $-10^{\circ}$ from root to tip. The rotor radius is 6.71 m and the blade chord is 0.686 m resulting in an aspect ratio of 9.8. The profile of the blades use a symmetric highly modified 540 airfoil section stemming from the NACA0012 family. In the aerodynamics phase of testing several specific thrust coefficients, tip Mach numbers, and advance ratios were flown to study the correlation between these parameters and the pressure distributions. 

The flight tests chosen for the present validation work are detailed in Table~\ref{chap4:tab:ah1g}. Two subsets of the flight tests are used. The first being a low-speed test which correlates to counter number 2157 from the database and the second being a high-speed test which correlates to counter number 2152 from the database. During the 150 km/h low-speed condition the rotor is operated at a tip Mach number of 0.65 and an advance ratio of 0.19. The collective pitch of $11.7^{\circ}$ is set to be congenial with the measured rotor thrust of $\textrm{C}_{\textrm{T}} = 0.00464$ and the cyclic pitch is trimmed to eliminate pitching and rolling moments. The sine and cosine cyclic pitch angles are respectively set to $1.7^{\circ}$ and $-5.5^{\circ}$. The sine and cosine cyclic flap angles are respectively set to $-0.15^{\circ}$ and $2.13^{\circ}$. During the 290 km/h high-speed condition the rotor is operated at a tip Mach number of 0.64 and an advance ratio of 0.24. The collective pitch of $18.0^{\circ}$ is set to be congenial with the measured rotor thrust of $\textrm{C}_{\textrm{T}} = 0.00464$ and the cyclic pitch is trimmed to eliminate pitching and rolling moments. The sine and cosine cyclic pitch angles are respectively set to $3.6^{\circ}$ and $-11.8^{\circ}$. The sine and cosine cyclic flap angles are respectively set to $1.11^{\circ}$ and $2.13^{\circ}$.

\begin{table}[htb!]
	\caption[Forward flight conditions taken from the TAAT data survey]{Forward flight conditions taken from the TAAT data survey \cite{cross1988tip}.}
	\centering
	\label{chap4:tab:ah1g}
	\begin{tabular}{ccc} 		
		\hline \hline
		Variable & \begin{tabular}[c]{@{}l@{}}Low Speed\\ Counter 2157\end{tabular} & \begin{tabular}[c]{@{}l@{}}High Speed\\ Counter 2152\end{tabular} 												\\ \hline
		$\textrm{V}_{\infty}$     		& 150 km/h   			& 290 km/h   			\\
		$\textrm{M}_{\infty}$     		& 0.12       			& 0.24       			\\
		$\textrm{M}_{\textrm{T}}$     	& 0.65       			& 0.64       			\\
		$\textrm{Re}_{\textrm{T}}$      & $9.7 \times 10^6$ 	& $10.2 \times 10^6$ 	\\
		$\mu$    						& 0.19       			& 0.38       			\\
		$\theta_{0}$ 					& $11.7^{\circ}$       	& $18.0^{\circ}$    	\\
		$\theta_{0.75}$ 				& $5.7^{\circ}$        	& $12.0^{\circ}$    	\\
		$\theta_{s}$ 					& $1.7^{\circ}$       	& $3.6^{\circ}$     	\\
		$\theta_{c}$ 					& $-5.5^{\circ}$       	& $-11.8^{\circ}$   	\\
		$\beta_{0}$  					& $2.75^{\circ}$       	& $2.75^{\circ}$     	\\
		$\beta_{c}$  					& $2.13^{\circ}$       	& $2.13^{\circ}$     	\\
		$\beta_{s}$  					& $-0.15^{\circ}$      	& $1.11^{\circ}$     	\\
		$\textrm{C}_{\textrm{T}}$    	& 0.00464    			& 0.00474    			\\ \hline \hline
	\end{tabular}
\end{table}

The spatial discretization is achieved using a single-zone mesh of the two main rotor blades. The outer far-field is placed 5 radii away from the blades. A mixed-element grid composed of 30.25 million elements and 760 thousand vertices is used. Each blade has 105 thousand surface elements. A density region surrounding the rotor blades is used to sufficiently resolve the near-field wake. Elements within the density region have a maximum element size of $0.0075 \, x/c$.

The standard one-equation SA turbulence model is used with the RANS equations. The SA turbulence variable is convected using a first-order scalar upwind method. The convective fluxes are computed using the Roe scheme and second-order accuracy is achieved using the MUSCL scheme. Spatial gradients are approximated using the Green-Guass numerical method. The flow solution was considered converged when there was a reduction of 6 orders of magnitude on the density residual. Per rotor revolution, there are 360 physical time-steps equating to $1^{\circ}$ azimuth increments. At each physical time-step, there are a maximum of 30 internal pseudo time-steps. Each simulation was run for a total of 5 rotor revolutions corresponding to 1800 physical time steps.  

\begin{figure}[hbt!]
	\centering
	\includegraphics[width=0.7\linewidth]{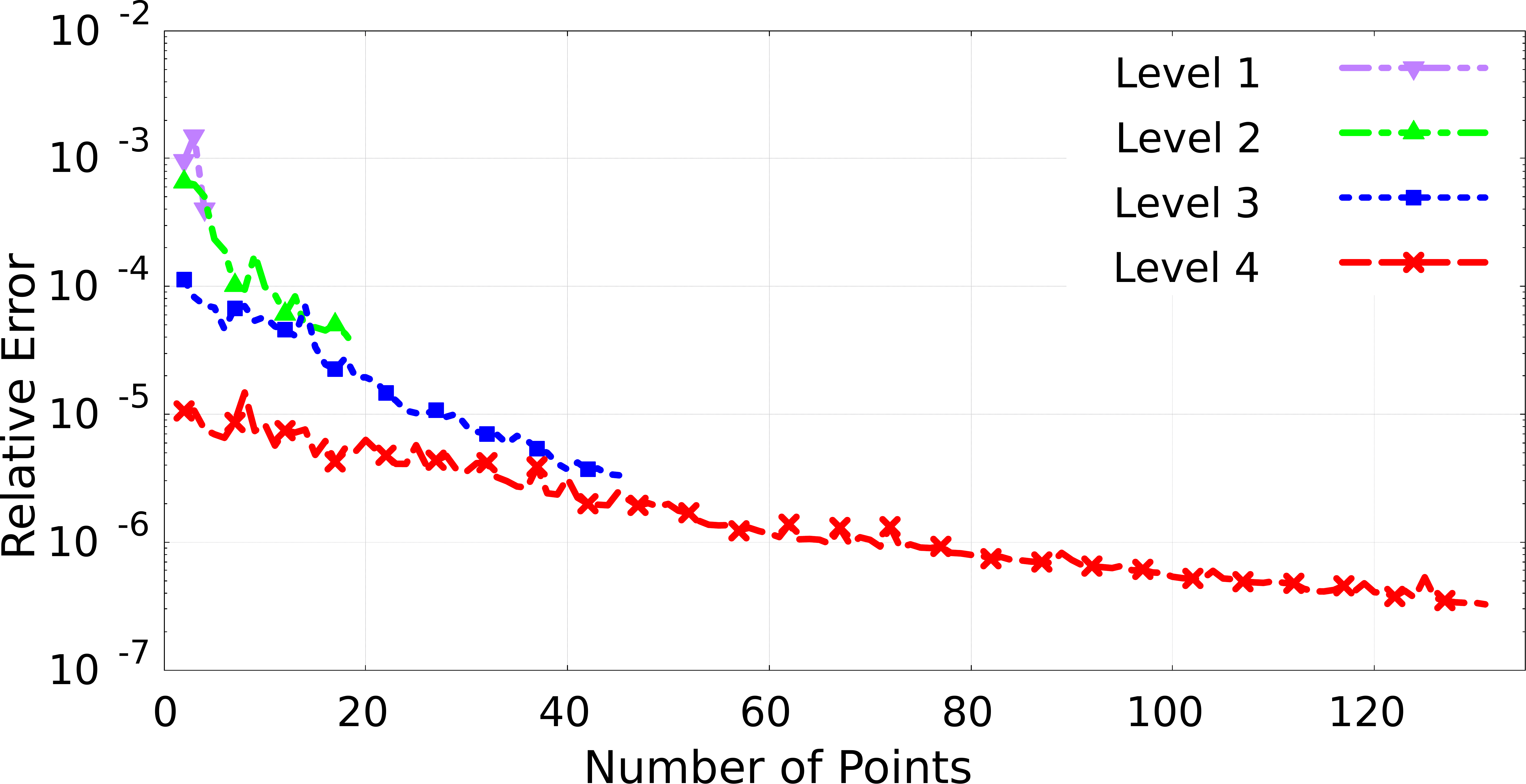}
	\caption[Multi-level convergence history on the AH-1G rotor]{Multi-level greedy surface point selection error reduction rates for a single deformation on the AH-1G rotor during forward flight.}
	\label{fig.14:AH1G_RBF_Efficiency}
\end{figure}

\begin{figure}[!htb]
	\subfloat[Level 1 - avg. 6 control points.
	\label{fig.13:subfig-1:AH1G_RBF}]{%
		\includegraphics[width=0.495\textwidth]{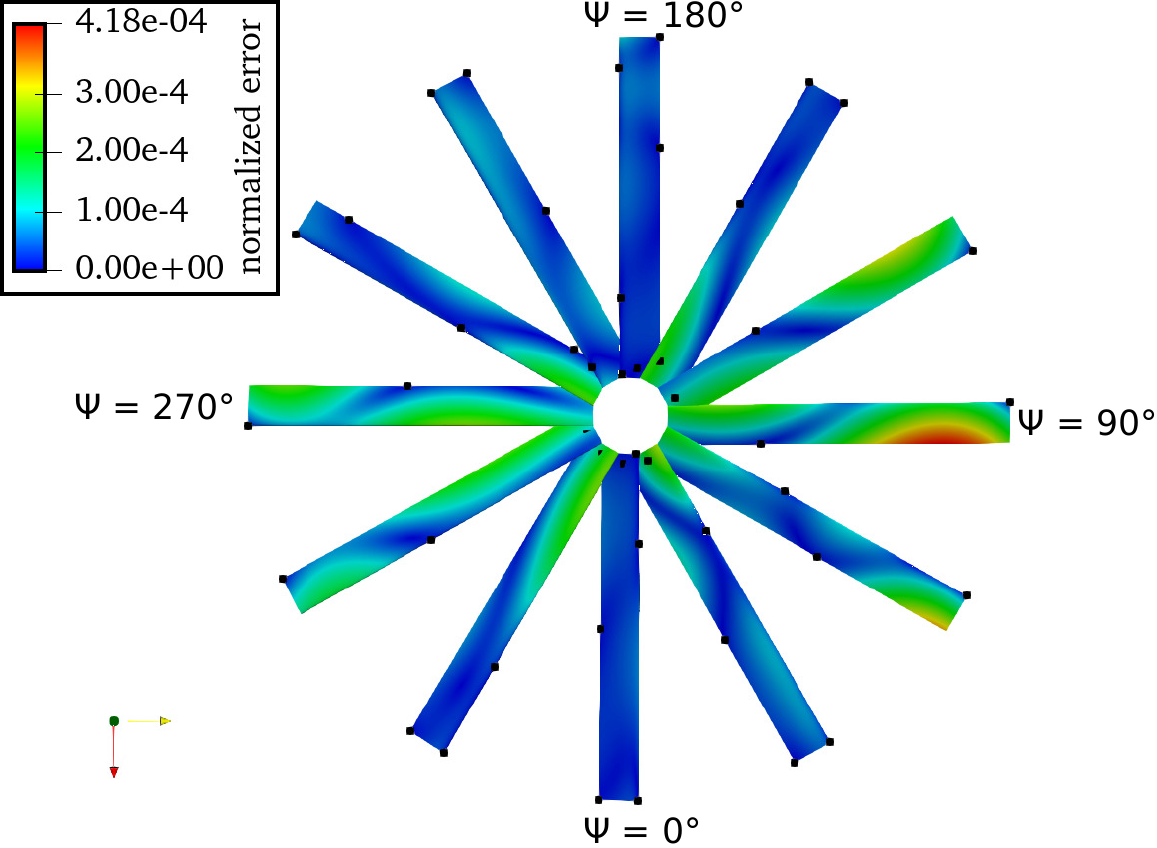}
	}
	\hfill
	\subfloat[Level 2 - avg. 21 control points.
	\label{fig.13:subfig-2:AH1G_RBF}]{%
		\includegraphics[width=0.495\textwidth]{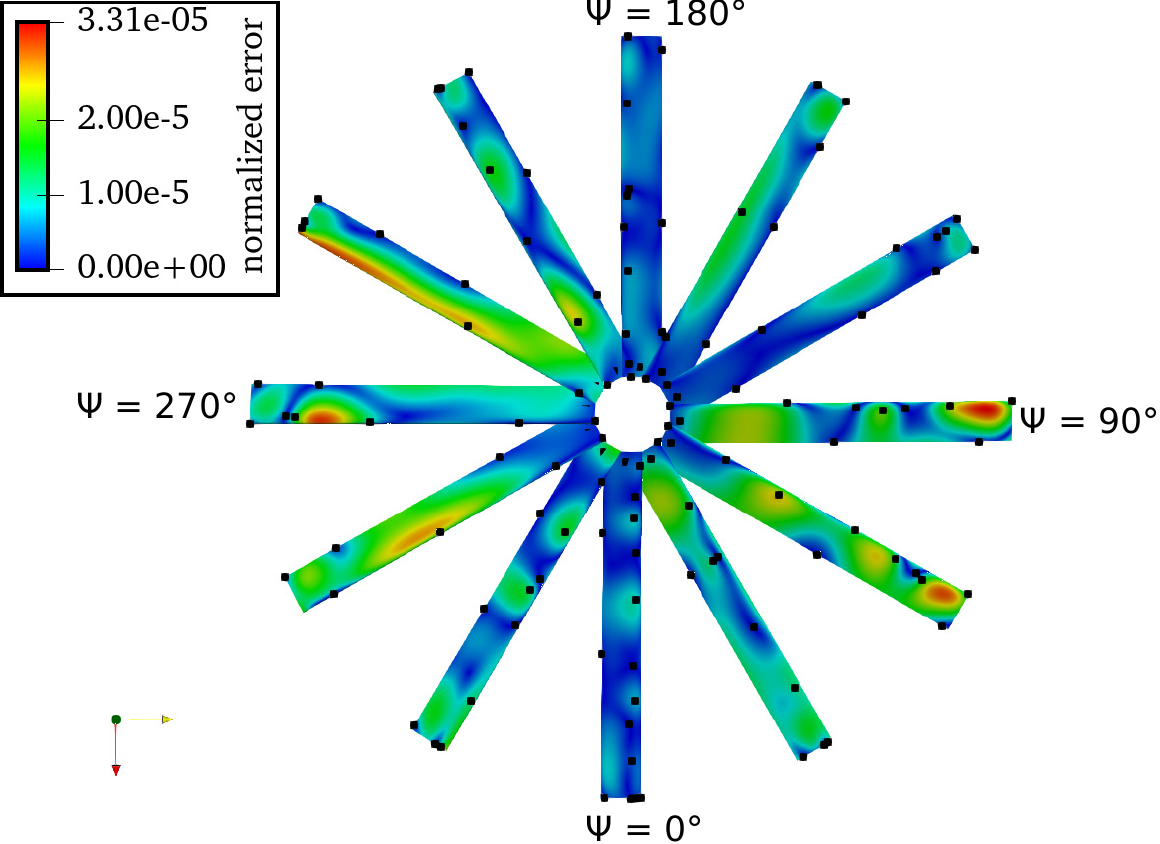}
	}\\
	\subfloat[Level 3 - avg. 42 control points.
	\label{fig.13:subfig-3:AH1G_RBF}]{%
		\includegraphics[width=0.495\textwidth]{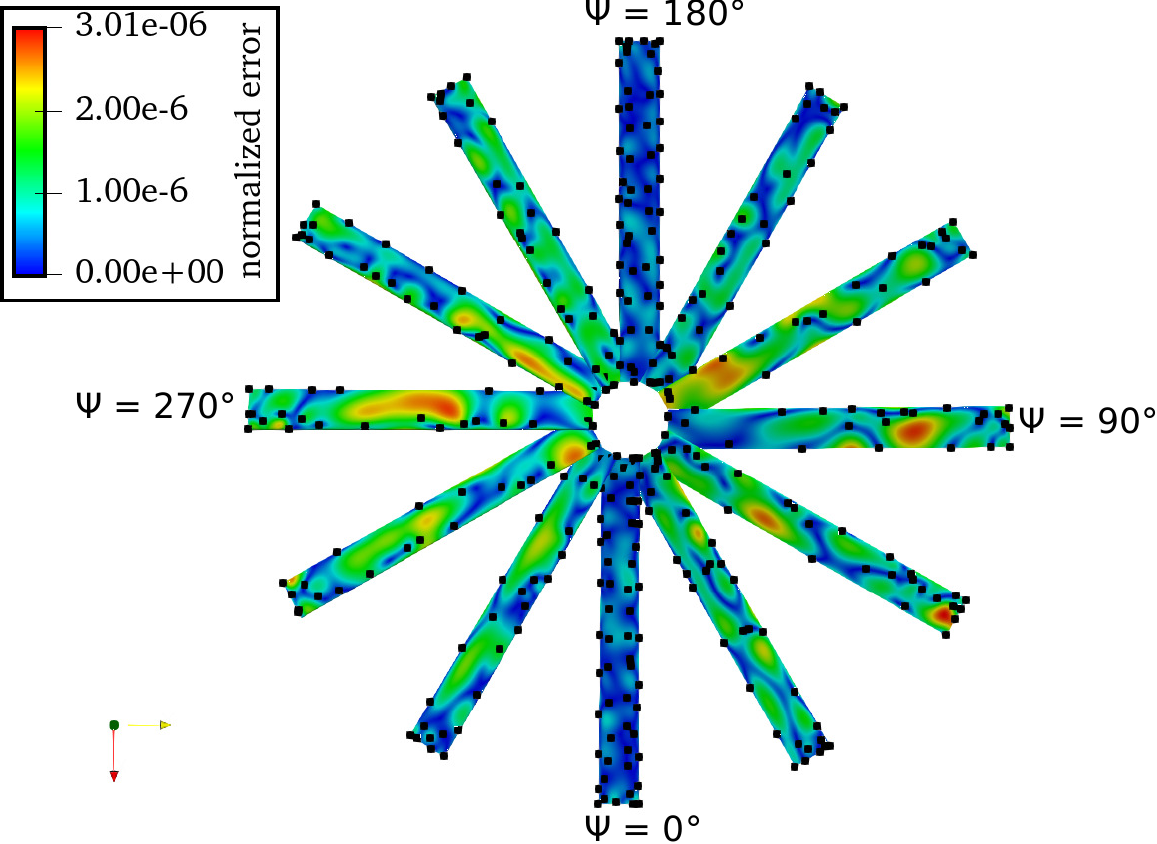}
	}
	\hfill
	\subfloat[Level 4 - avg. 136 control points.
	\label{fig.13:subfig-4:AH1G_RBF}]{%
		\includegraphics[width=0.495\textwidth]{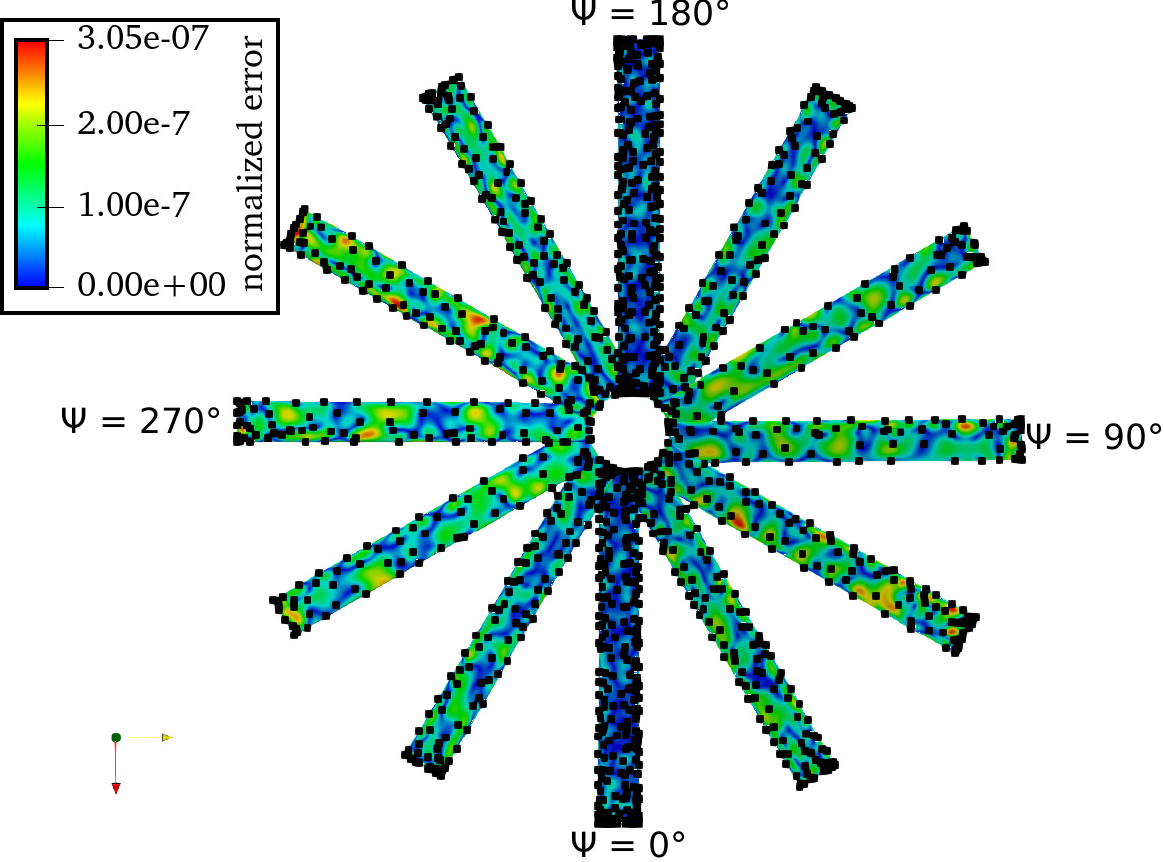}
	}
	\caption[Selected control points on the AH-1G rotor]{Evolution of the multi-level greedy surface point selection algorithm used for the RBF mesh deformation. Displaying the control points as black dots and the contours of the normalized surface error displacement. Where the normalized error is computed as the difference between the computed displacement and the actual displacement.}
	\label{fig.13:AH1G_RBF}
\end{figure}

The blade motion is prescribed using the implementation introduced in Section~\ref{sec:kinematics} with the values outlined in Table~\ref{chap4:tab:ah1g}. The mesh is updated at each physical time-step using the RBF mesh deformation strategy. Multi-level surface point reduction and volume point reduction algorithms are used to improve the efficiency. Four levels are employed with the surface point selection error reduction rates and efficiency being shown in Fig.~\ref{fig.14:AH1G_RBF_Efficiency}. It shows that the multi-level method for updating the blade position is highly effective in minimizing the computational cost while simultaneously reducing the surface error. The control points and associated normalized surface error of the computed displacement at each level are shown in Fig.~\ref{fig.13:AH1G_RBF}. The minimum orthogonality angle of the grid prior to deformation was $40.979^{\circ}$. The minimum orthogonality angle of the grid post 1800 deformations is $40.997^{\circ}$. This indicates there is no detrimental effect on the grid quality from the continuous deformation.  

\begin{figure}[!htb]
	\subfloat[$\psi = 30^{\circ}$.
	\label{fig.7:subfig-1:AH1G_surfacePressure}]{%
		\includegraphics[width=0.32\textwidth]{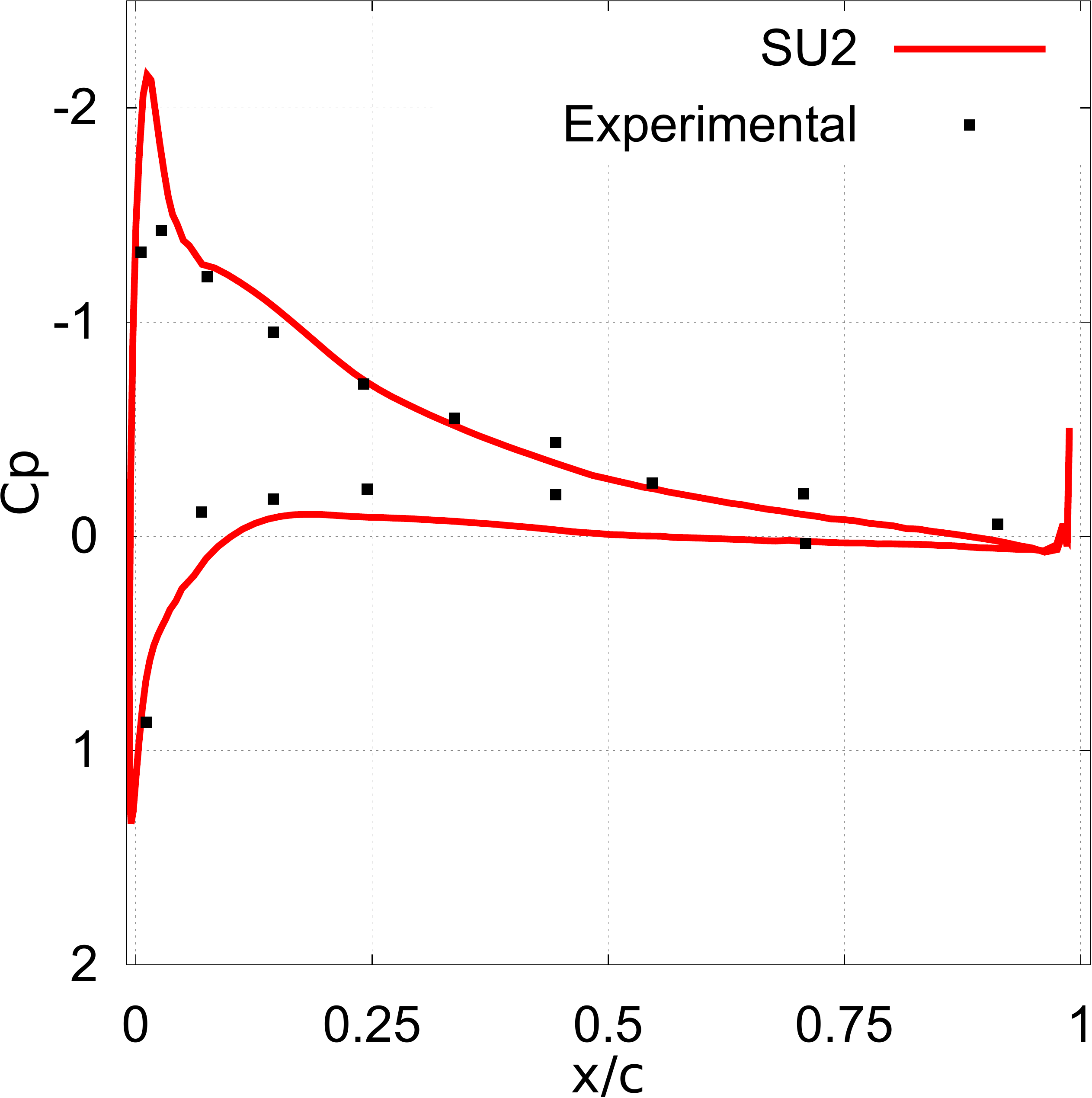}
	}
	\hfill
	\subfloat[$\psi = 90^{\circ}$.
	\label{fig.7:subfig-2:AH1G_surfacePressure}]{%
		\includegraphics[width=0.32\textwidth]{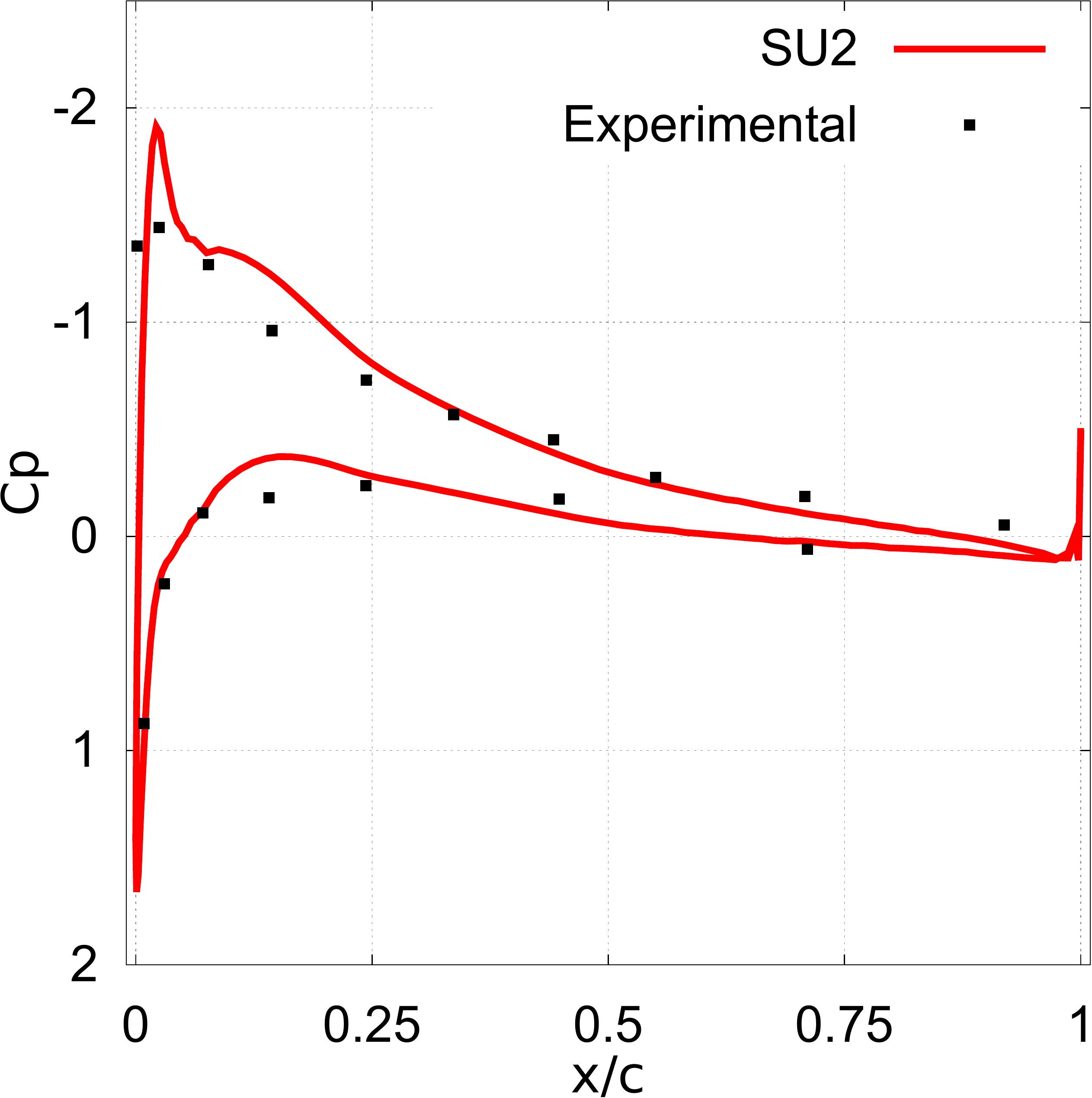}
	}
	\hfill
	\subfloat[$\psi = 180^{\circ}$.
	\label{fig.7:subfig-3:AH1G_surfacePressure}]{%
		\includegraphics[width=0.32\textwidth]{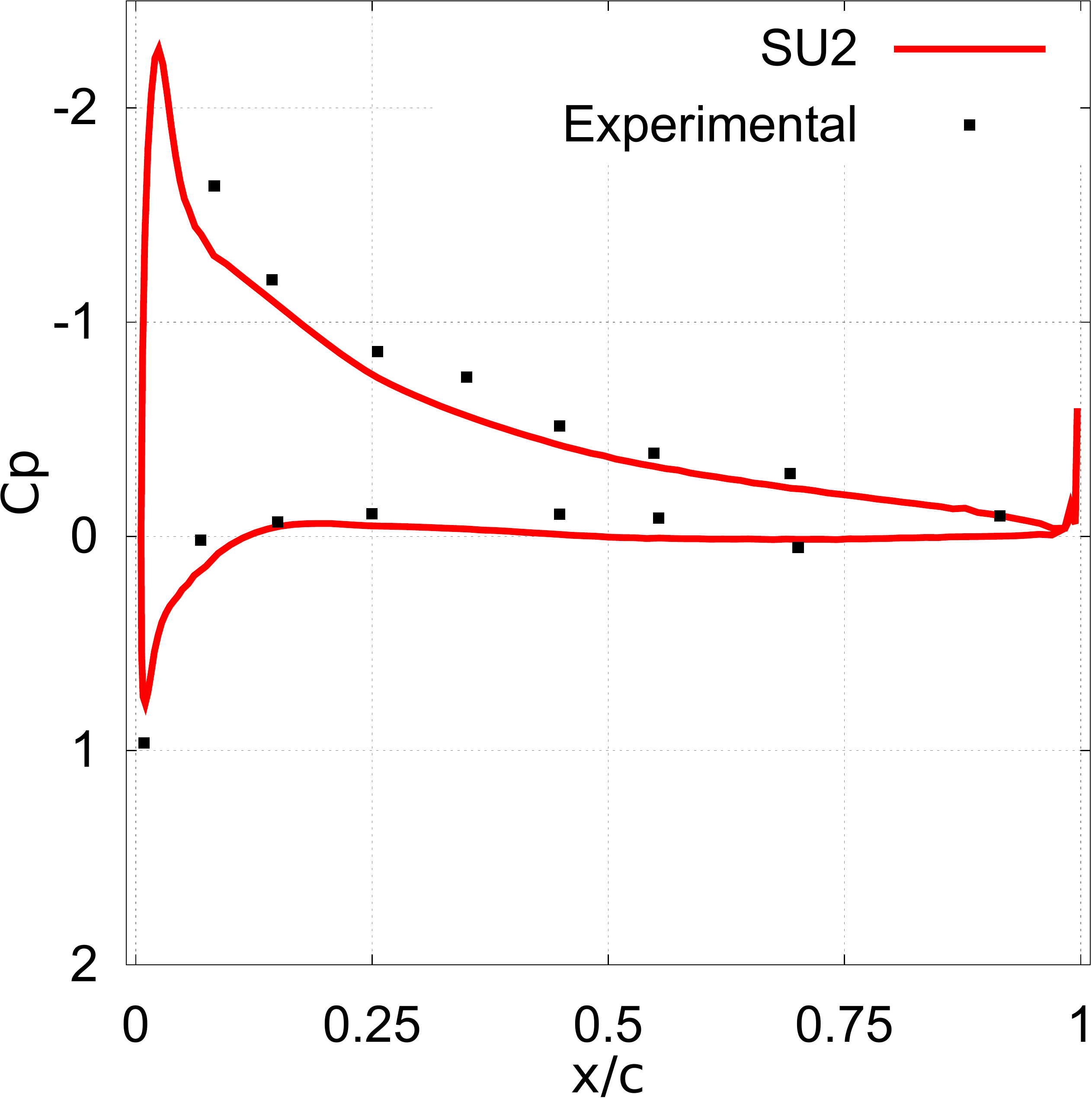}
	}
	\caption[Low-speed performance prediction on the advancing side]{Pressure distributions at $r/R = 0.60$ on the advancing side of the rotor during low-speed forward flight. Predictions compared against the measured data taken from the TAAT \cite{cross1988tip}.}
	\label{fig.7:AH1G_surfacePressure}
\end{figure}

\begin{figure}[!htb]
	\subfloat[$\psi = 270^{\circ}$.
	\label{fig.X:subfig-4:AH1G_surfacePressure}]{%
		\includegraphics[width=0.32\textwidth]{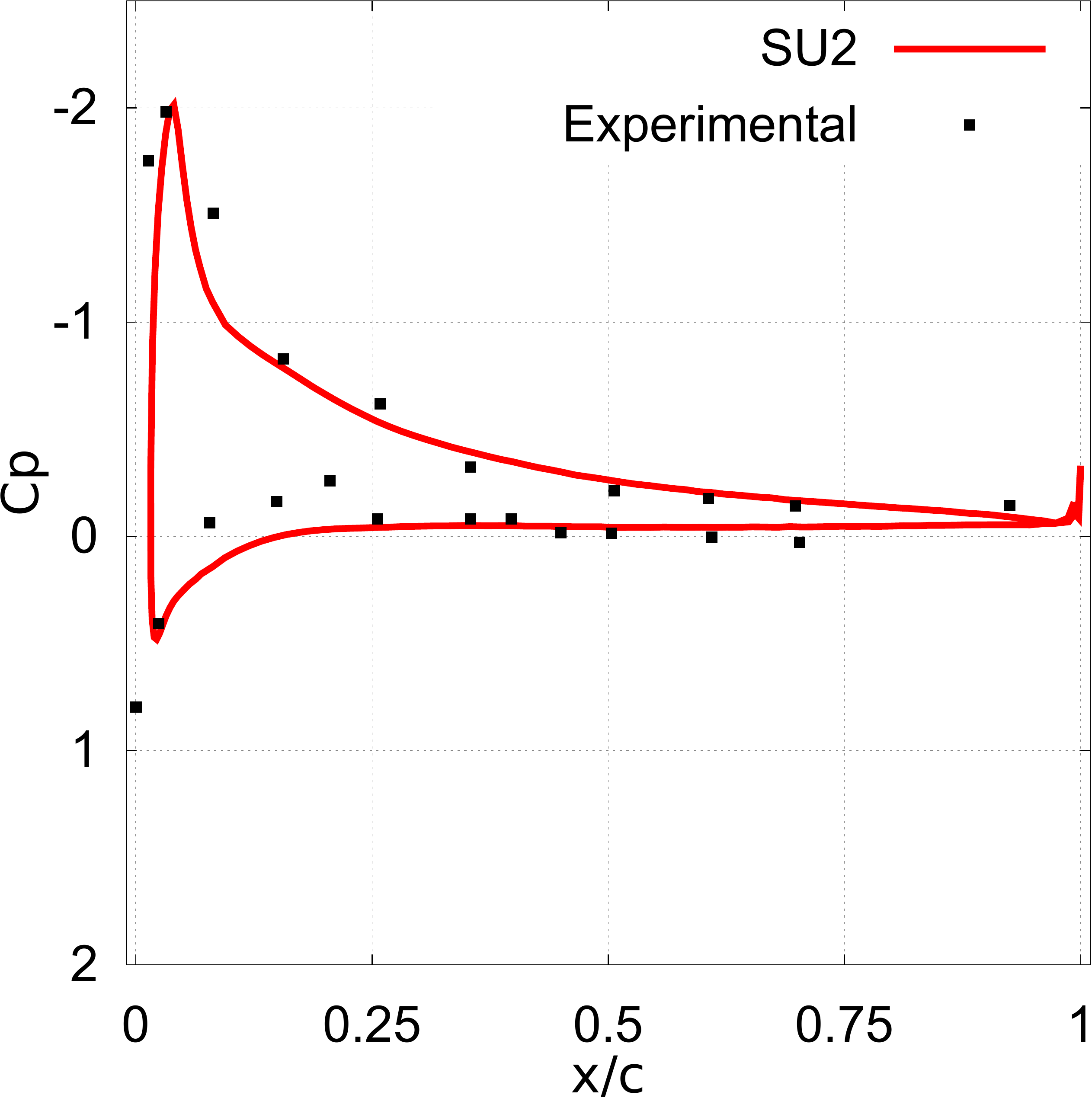}
	}
	\hfill
	\subfloat[$\psi = 285^{\circ}$.
	\label{fig.X:subfig-5:AH1G_surfacePressure}]{%
		\includegraphics[width=0.32\textwidth]{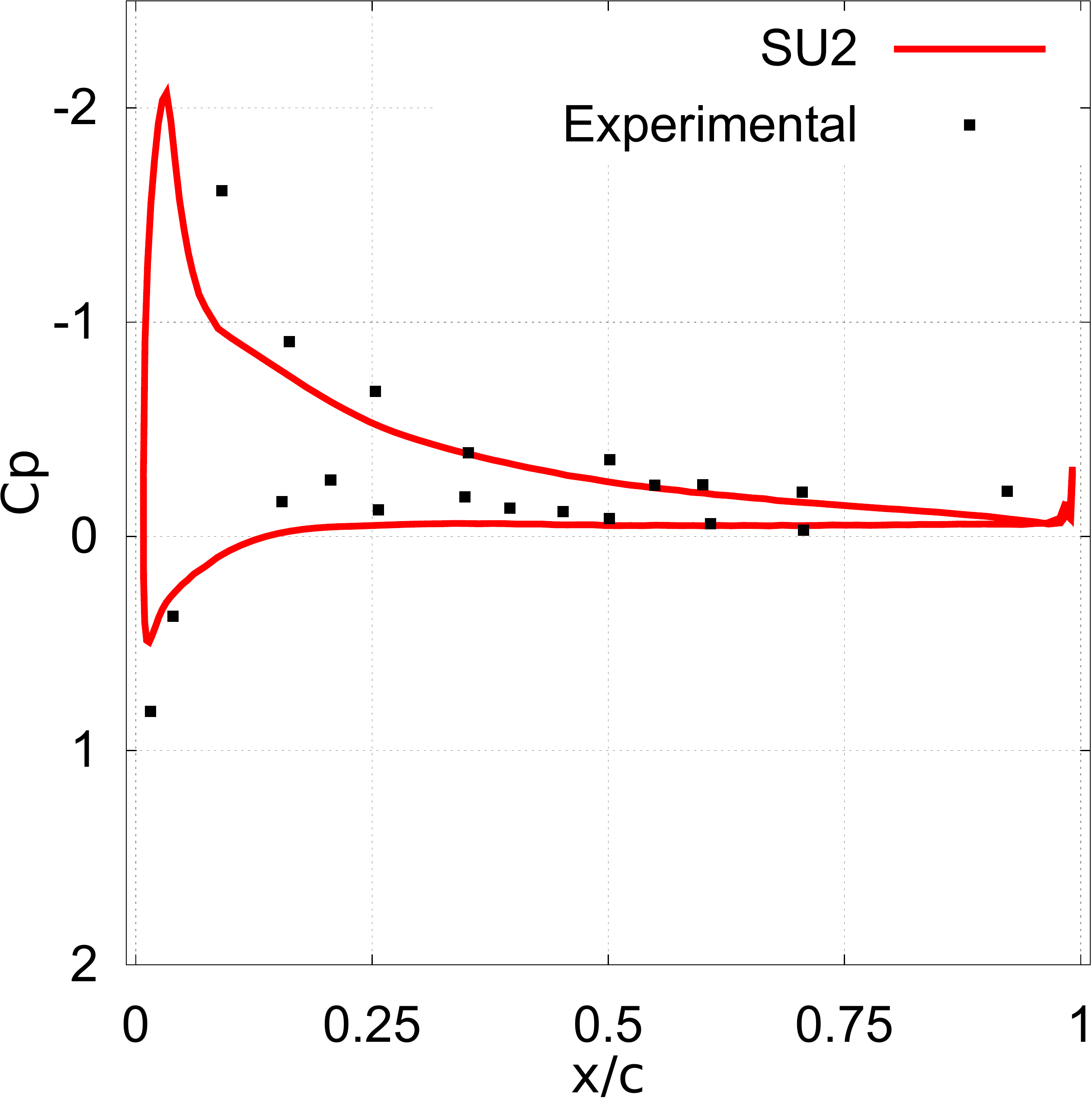}
	}
	\hfill
	\subfloat[$\psi = 300^{\circ}$.
	\label{fig.X:subfig-6:AH1G_surfacePressure}]{%
		\includegraphics[width=0.32\textwidth]{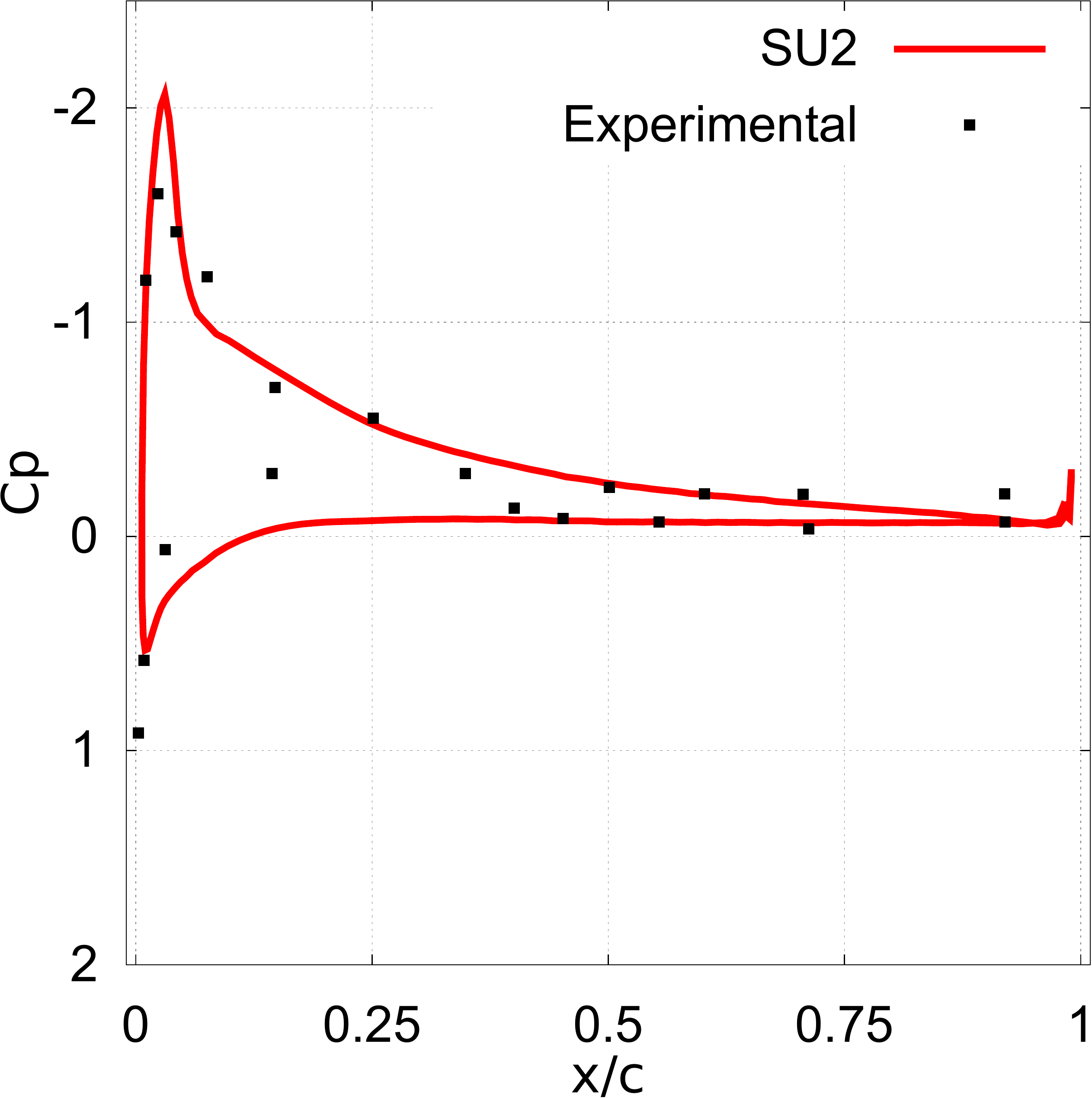}
	}
	\caption[Low-speed performance prediction on the retreating side]{Pressure distributions at $r/R = 0.91$ on the retreating side of the rotor during low-speed forward flight. Predictions compared against the measured data taken from the TAAT \cite{cross1988tip}.}
	\label{fig.8:AH1G_surfacePressure}
\end{figure}

The pressure coefficient distributions at selected advancing and retreating azimuth positions as well as radial positions are obtained from the data survey and are used for validation. The low-speed condition uses data at $r/R = 0.60$ for the advancing side of the rotor at azimuth angles $\psi = 30^{\circ}$, $\psi = 90^{\circ}$, and $\psi = 180^{\circ}$ with the results being shown in Fig.~\ref{fig.7:AH1G_surfacePressure}. Data at $r/R = 0.91$ on the retreating side of the rotor at azimuth angles $\psi = 270^{\circ}$, $\psi = 285^{\circ}$, and $\psi = 300^{\circ}$ are also used and the results are shown in Fig.~\ref{fig.8:AH1G_surfacePressure}. In general, the trends are relatively well captured when compared against the measured data. On the advancing side of the rotor the on the upper surface the suction peaks are marginally over-exaggerated. On the retreating side of the rotor this is improved slightly.

\begin{figure}[!htb]
	\subfloat[$\psi = 70^{\circ}$.
	\label{fig.9:subfig-1:AH1G_surfacePressure}]{%
		\includegraphics[width=0.32\textwidth]{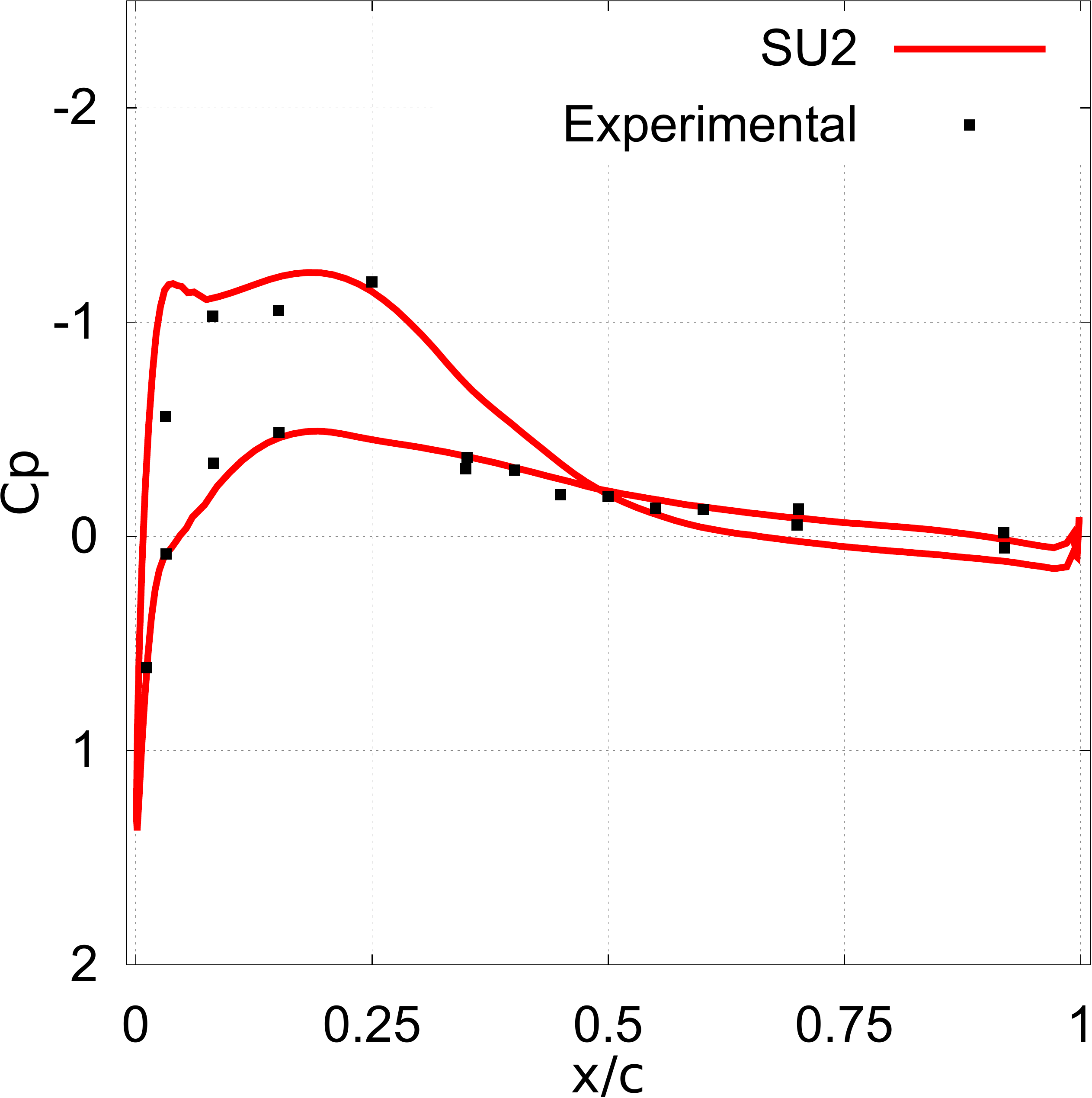}
	}
	\hfill
	\subfloat[$\psi = 90^{\circ}$.
	\label{fig.9:subfig-2:AH1G_surfacePressure}]{%
		\includegraphics[width=0.32\textwidth]{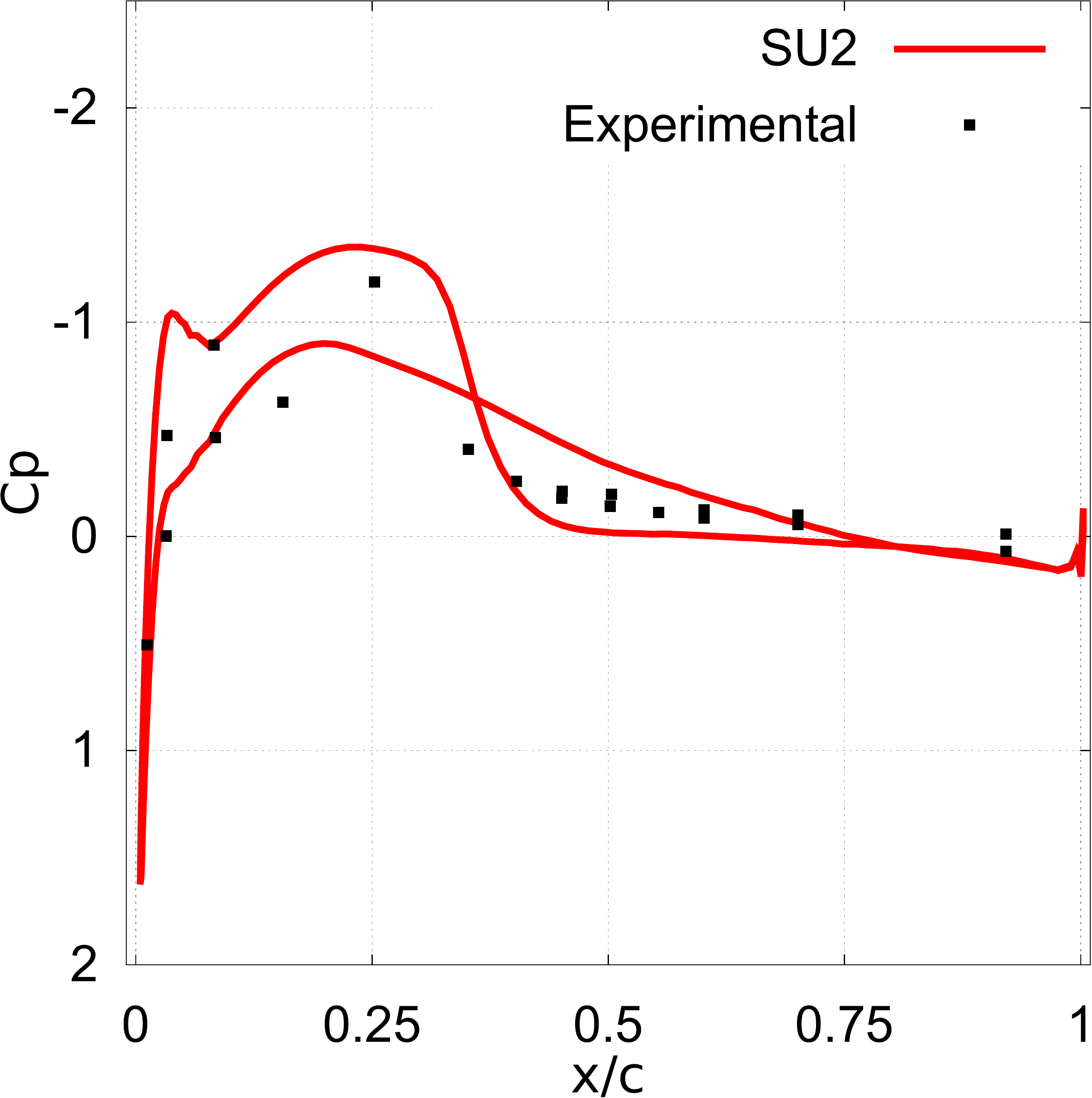}
	}
	\hfill
	\subfloat[$\psi = 110^{\circ}$.
	\label{fig.9:subfig-3:AH1G_surfacePressure}]{%
		\includegraphics[width=0.32\textwidth]{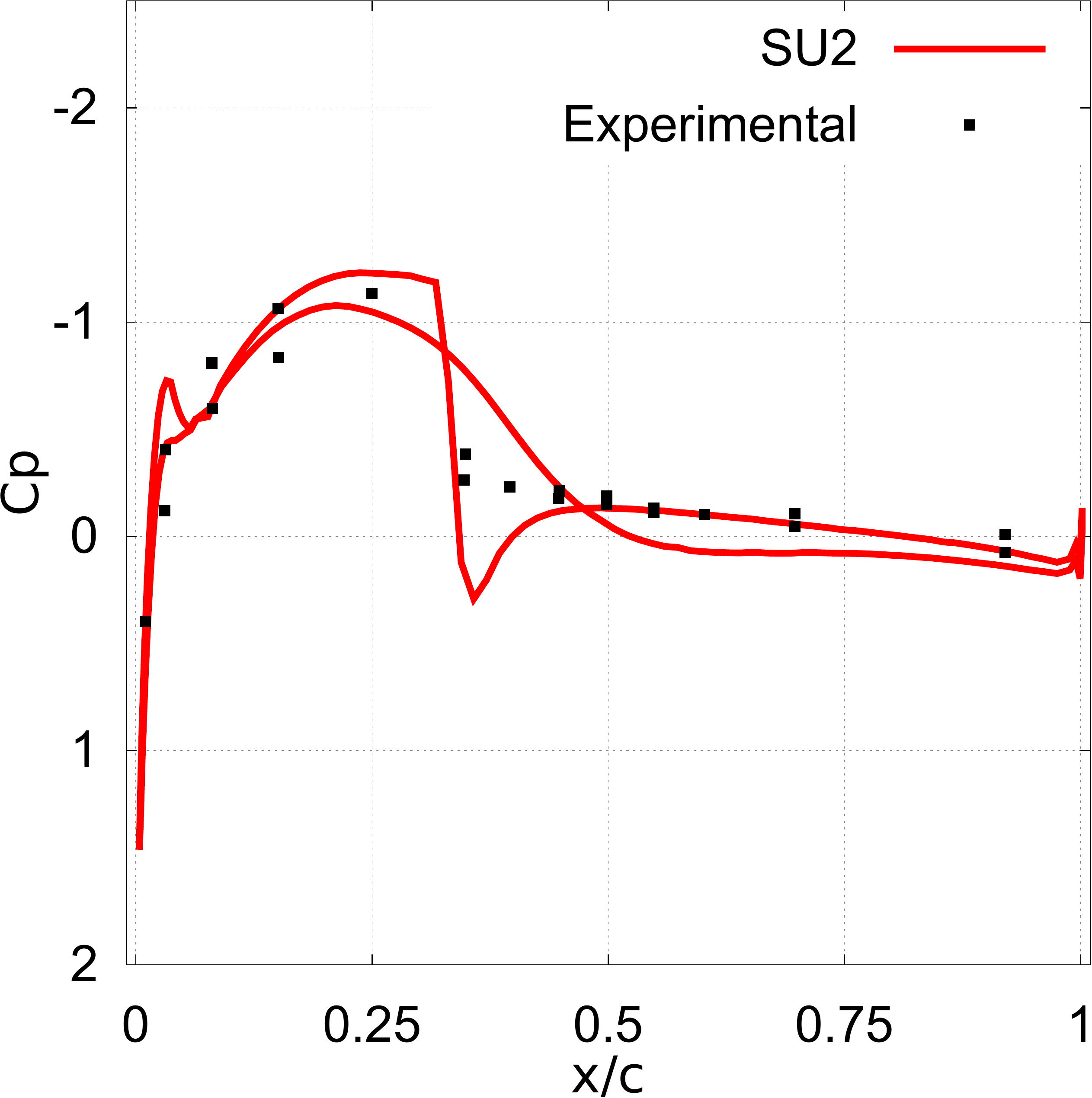}
	}
	\caption[High-speed performance prediction on the advancing side]{Pressure distributions at $r/R = 0.86$ on the advancing side of the rotor during high-speed forward flight. Predictions compared against the measured data taken from the TAAT \cite{cross1988tip}.}
	\label{fig.9:AH1G_surfacePressure}
\end{figure}

The high-speed conditions uses data at $r/R = 0.86 $ for the advancing side of the rotor at azimuth angles $\psi = 70^{\circ}$, $\psi = 90^{\circ}$, and $\psi = 110^{\circ}$ with the results being shown in Fig.~\ref{fig.9:AH1G_surfacePressure}. Data at $r/R = 0.96$ on the retreating side of the rotor at azimuth angles $\psi = 250^{\circ}$, $\psi = 270^{\circ}$, and $\psi = 290^{\circ}$ are also used and the results are shown in Fig~\ref{fig.10:AH1G_surfacePressure}. The higher flight speed initiates the presence of shocks on the advancing side of the blade and dynamic stall on the retreating side of the blade which makes replicating the measured data more challenging than the low-speed condition. Consequently there are significantly different pressure coefficient distributions. The suction peaks on the upper surface of the blade on the advancing side of the rotor are represented well. The shock positions however show discrepancies when compared to the measured data. The pressure coefficient distributions on the retreating side of the blade close to the tip appear to be better represented. 

A contour map of the pressure coefficient on both the suction and pressure sides of the blades for both the low- and high-speed conditions are shown in Fig.~\ref{fig.11:AH1G_surfPressure}. Both sets of results show that the pressure coefficient during forward flight is no longer axisymmetric about the out-of-plane axis. To achieve the higher flight speed an increase in the collective pitch is required. This however results in a low speed stalled region on the retreating side of the rotor near the blade root where there is a high angle-of-attack. This effect can be seen in Fig.~\ref{fig.11:subfig-2:AH1G_surfPressureHS} which is not present in Fig.~\ref{fig.11:subfig-1:AH1G_surfPressureLS}.   

The iso-surface of the Q-criterion visualizing the near-field wake and blade tip vortices is shown in Fig.~\ref{fig12:AH1G_qCrit}. It displays the different flow field behaviour in low-speed and high-speed forward flight. In low-speed forward flight, the vortex upstream of the preceding blade passes below the advancing blade close to the tip. Due to the low flight speed, the vortex downstream remains close to the rotor blades. In high-speed forward flight the vortex upstream of the preceding blade passes below the advancing blade closer to the root of the blade. With the flight speed being much higher the vortex downstream is no longer in the vicinity of the rotor blades. Due to both cases being straight and level flight there are no blade-vortex interactions present. It is of note that even with a density region and such a large grid that it is challenging to reduce artificial numerical dissipation of the blade tip vortices.   

\section{Conclusion}
\label{chap4:sec:conclusion}
This work provides a unique approach for the modelling of rotorcraft aerodynamics using the open-source toolkit SU2. The performance prediction capabilities are demonstrated on two test cases for preliminary validation. The first test case of the Caradonna-Tung model rotor in hover allowed for the assessment of the pressure coefficient distribution along four different radial positions. Models with varying fidelity ranging from steady-state to fully time-accurate were assessed alongside significant influential parameters. The second test case of the AH-1G rotor in forward flight required the introduction of the blade motion. Low and high speed forward flight tests were modelled. The low-speed test condition was found to be easier to simulate and the performance predictions are in good agreement. The high-speed test condition was more challenging due to the presence of a shock near the blade tip on the advancing side and dynamic stall on the retreating side of the blade near the root. Despite the results of this work being complementary to the measured data, it is not yet sufficient to suggest this is satisfactory for complete validation. It is import that more tests are carried out with different configurations and test conditions. Moving forward, likely areas of research will include; rotor-fuselage interaction effects, blade design optimization, rotor acoustics, mesh adaptation, and the extension of rigid blades to elastic blades. 

\newpage
\section*{Figures}

\begin{figure}[!htb]
	\subfloat[$\psi = 250^{\circ}$.
	\label{fig.10:subfig-1:AH1G_surfacePressure}]{%
		\includegraphics[width=0.32\textwidth]{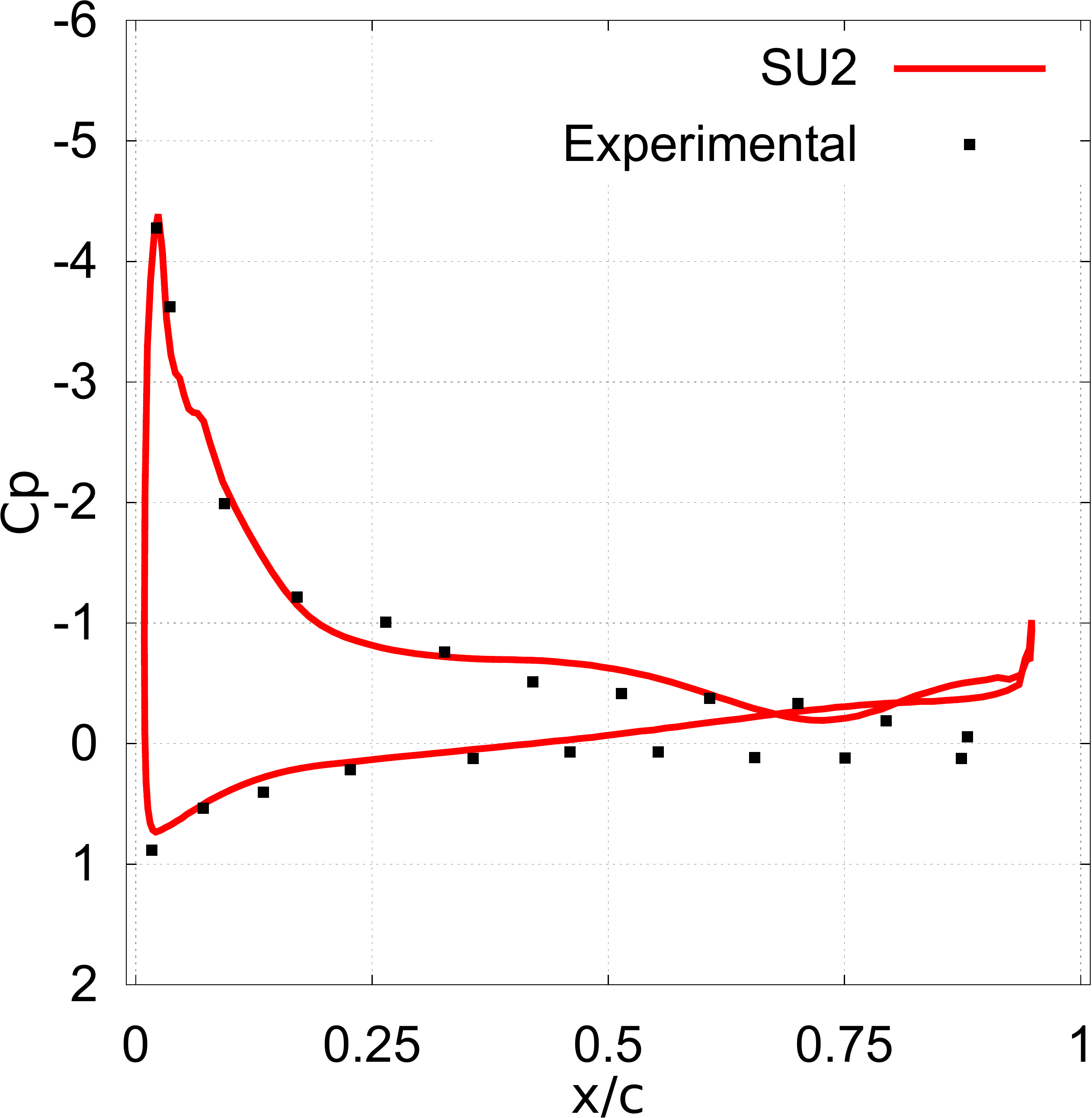}
	}
	\hfill
	\subfloat[$\psi = 270^{\circ}$.
	\label{fig.10:subfig-2:AH1G_surfacePressure}]{%
		\includegraphics[width=0.32\textwidth]{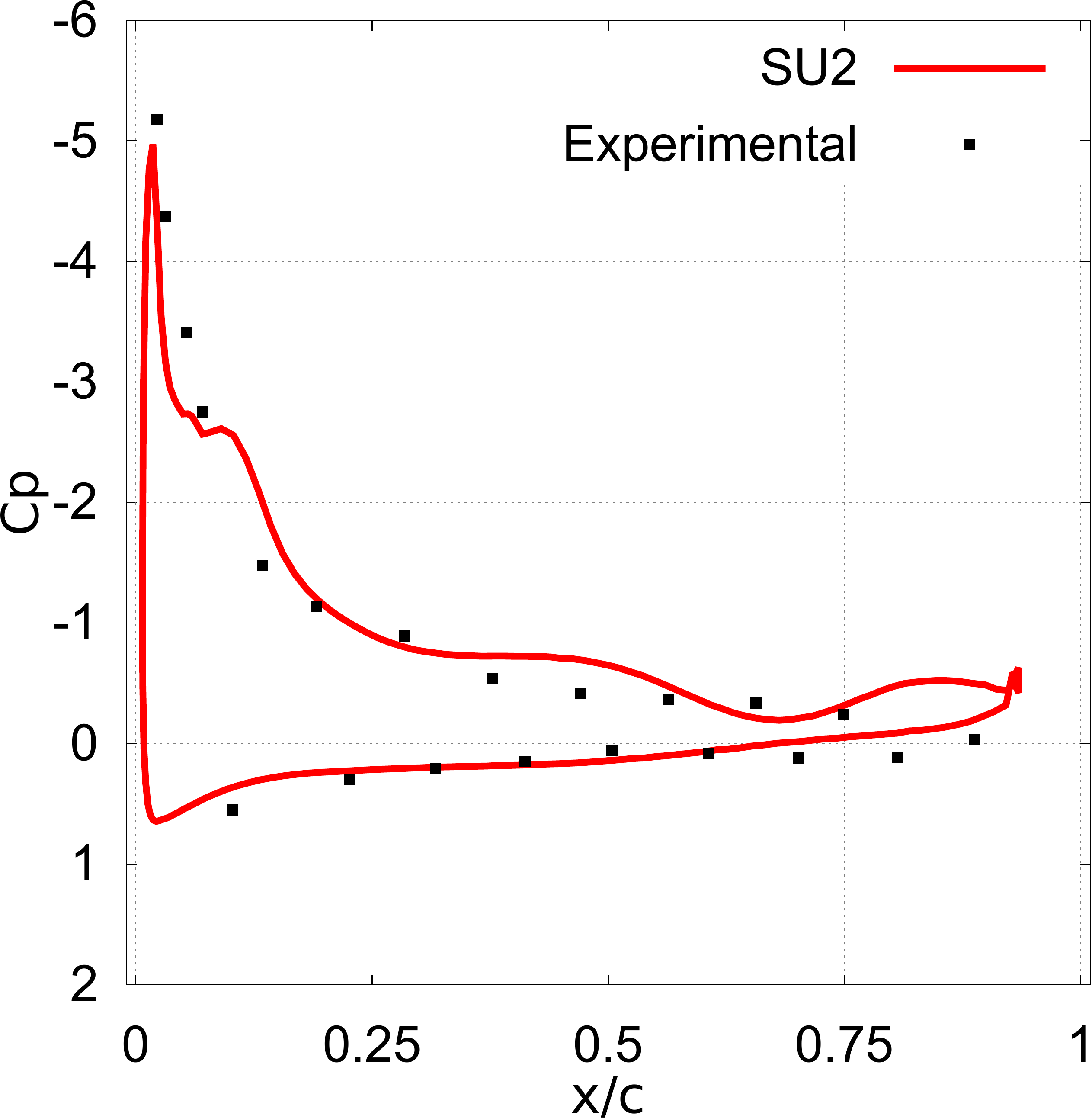}
	}
	\hfill
	\subfloat[$\psi = 290^{\circ}$.
	\label{fig.10:subfig-3:AH1G_surfacePressure}]{%
		\includegraphics[width=0.32\textwidth]{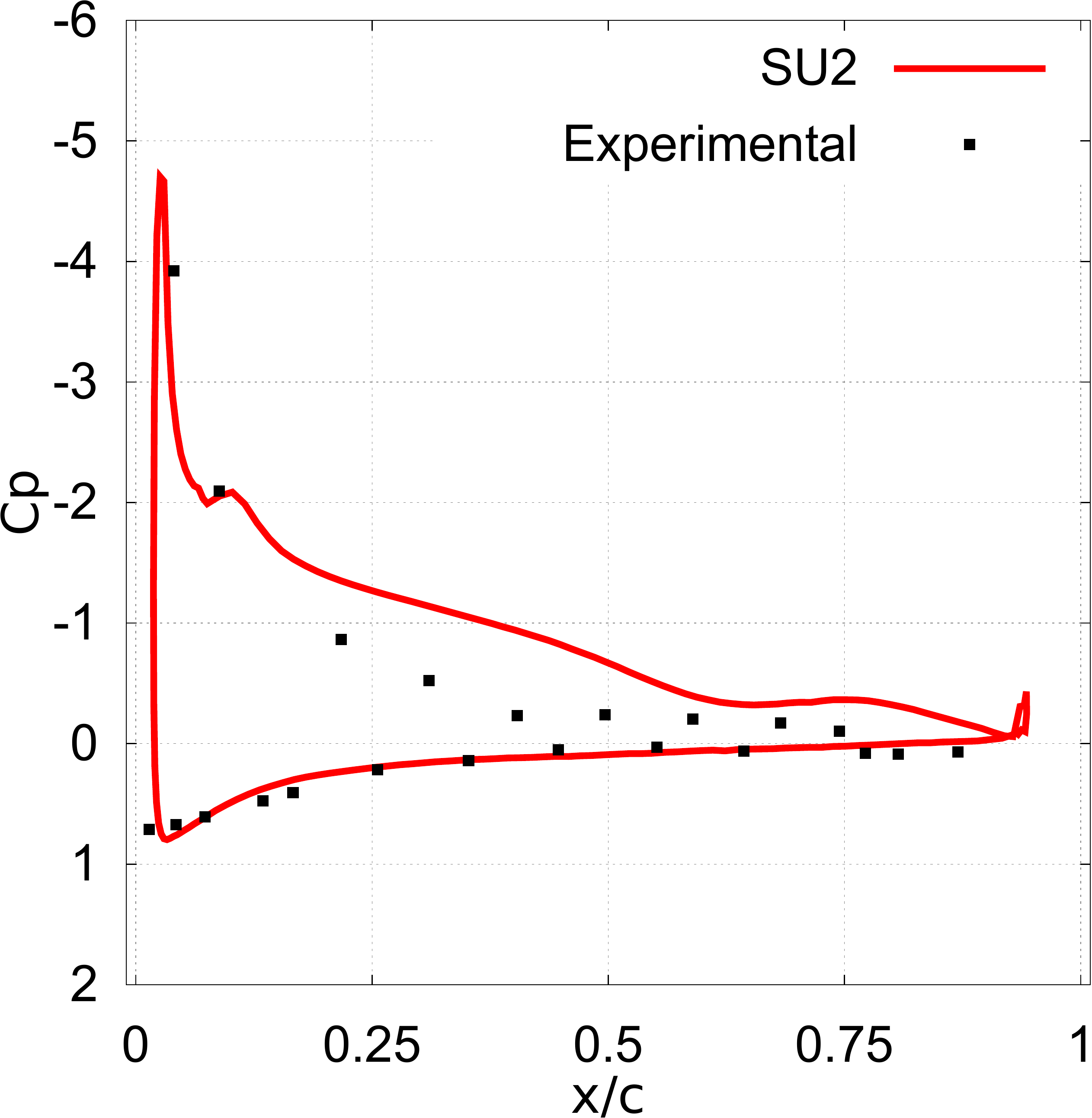}
	}
	\caption[High-speed performance prediction on the retreating side]{Pressure distributions at $r/R = 0.96$ on the advancing side of the rotor during high-speed forward flight. Predictions compared against the measured data taken from the TAAT \cite{cross1988tip}.}
	\label{fig.10:AH1G_surfacePressure}
\end{figure}

\begin{figure}[!htb]
	\subfloat[$r/R = 0.75$.
	\label{fig.10:subfig-1:AH1G_surfacePressure}]{%
		\includegraphics[trim={0.35cm 0.35cm 0.4cm 0.2cm},clip, width=0.495\textwidth]{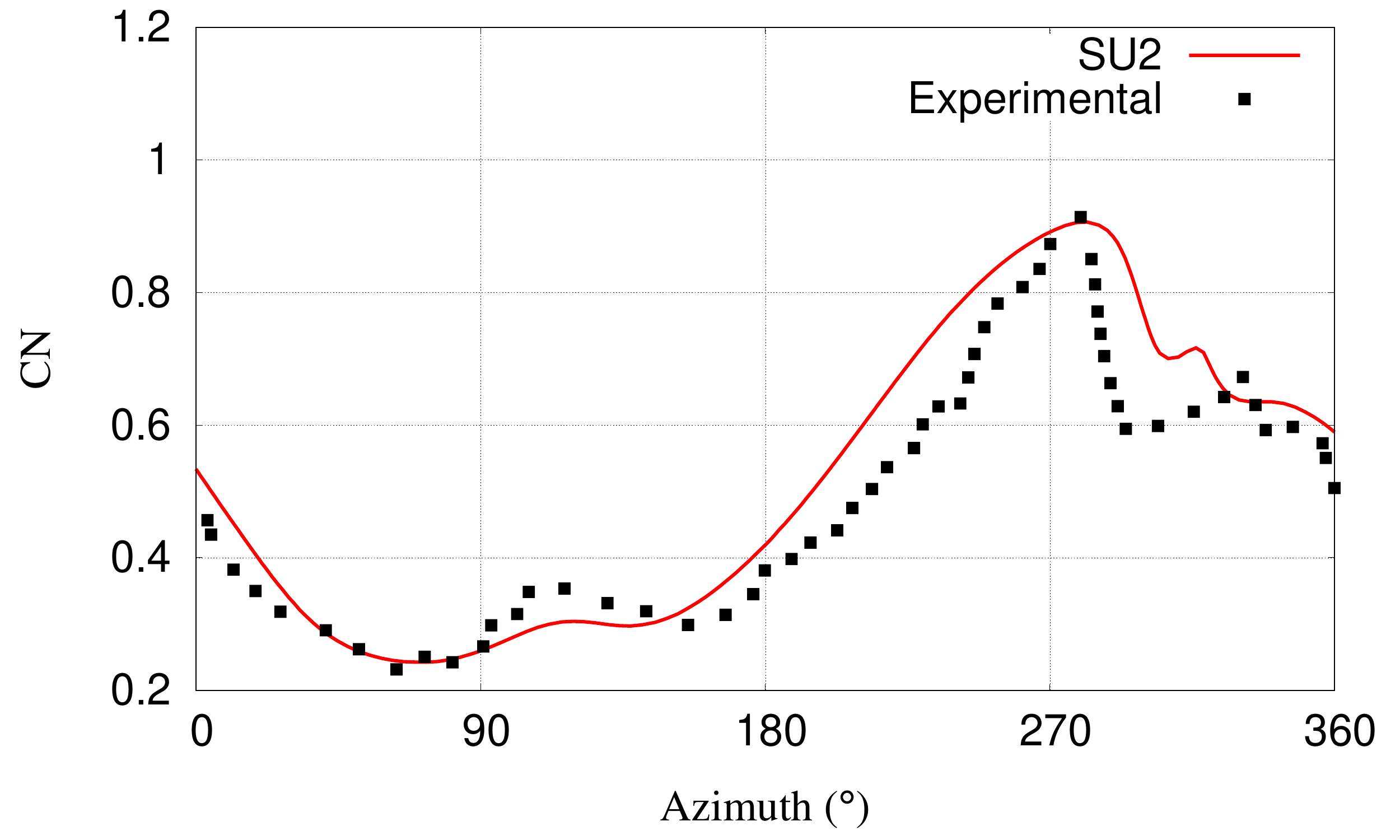}
	}
	\hfill
	\subfloat[$r/R = 0.86$.
	\label{fig.10:subfig-2:AH1G_surfacePressure}]{%
		\includegraphics[trim={0.35cm 0.35cm 0.4cm 0.2cm},clip, width=0.495\textwidth]{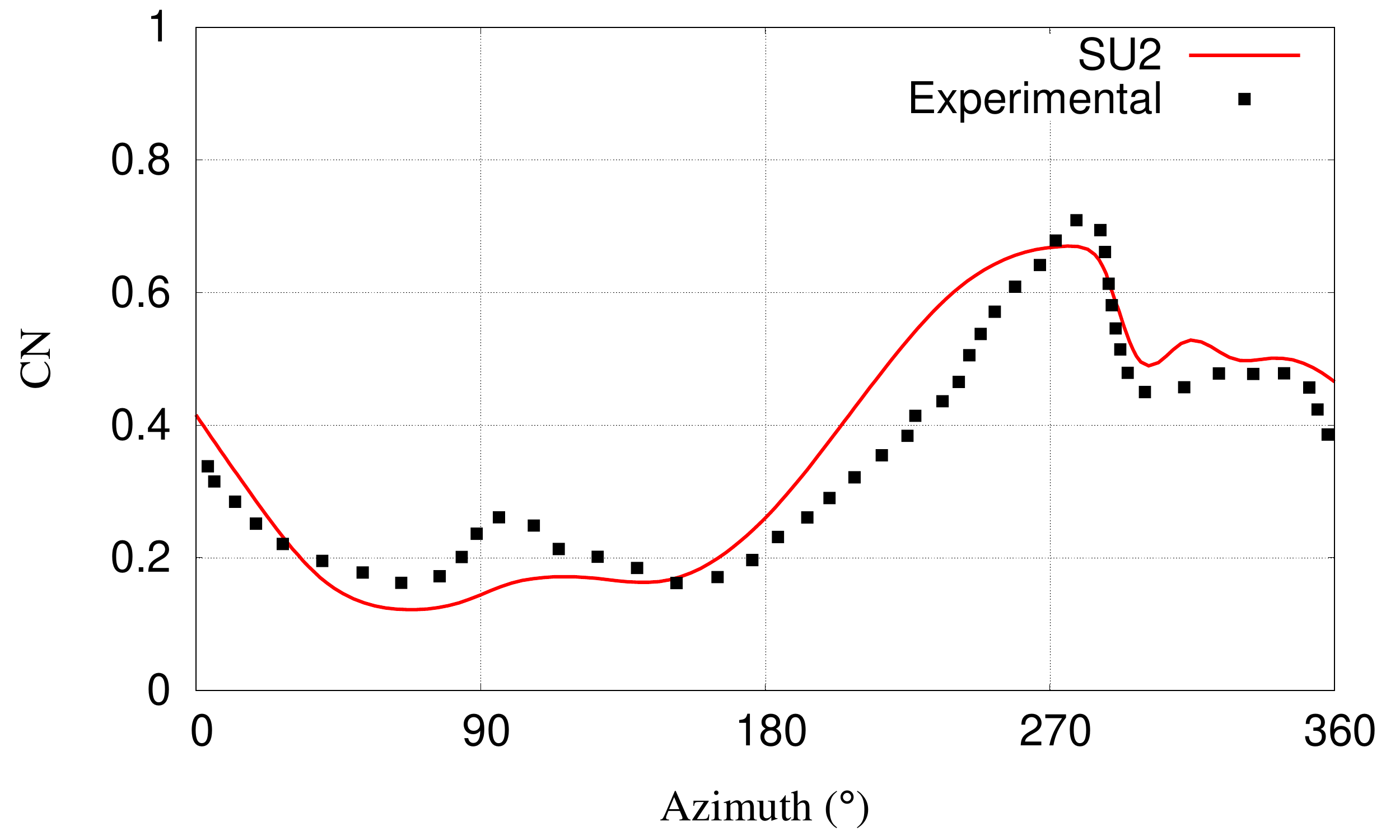}
	}
	\caption{Predictions compared against the measured data taken from the TAAT \cite{cross1988tip}.}
	\label{fig.10:AH1G_surfacePressure}
\end{figure}

\begin{figure}[!htb]
	\subfloat[$r/R = 0.75$.
	\label{fig.10:subfig-1:AH1G_surfacePressure}]{%
		\includegraphics[trim={0.35cm 0.35cm 0.4cm 0.2cm},clip, width=0.495\textwidth]{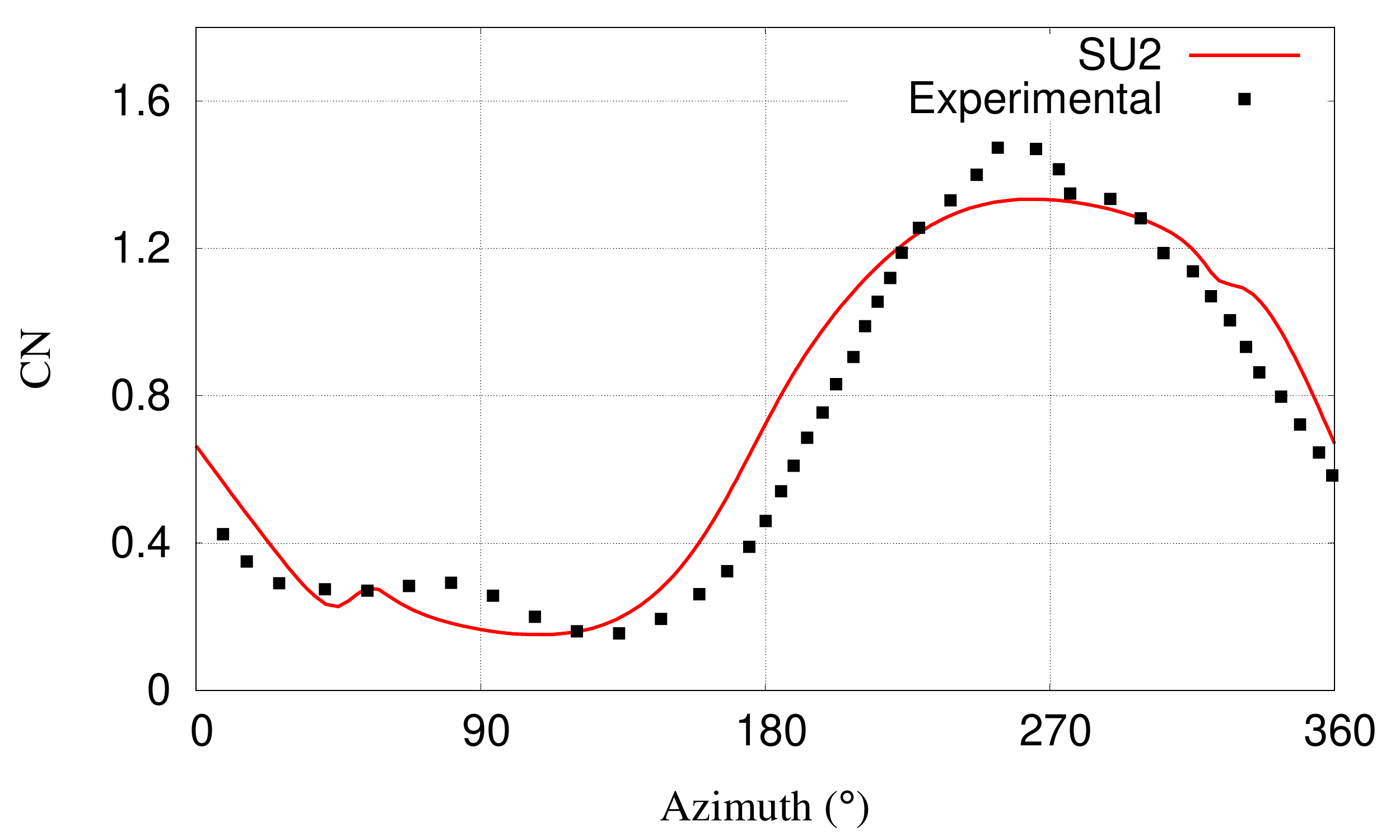}
	}
	\hfill
	\subfloat[$r/R = 0.86$.
	\label{fig.10:subfig-2:AH1G_surfacePressure}]{%
		\includegraphics[trim={0.35cm 0.35cm 0.4cm 0.2cm},clip, width=0.495\textwidth]{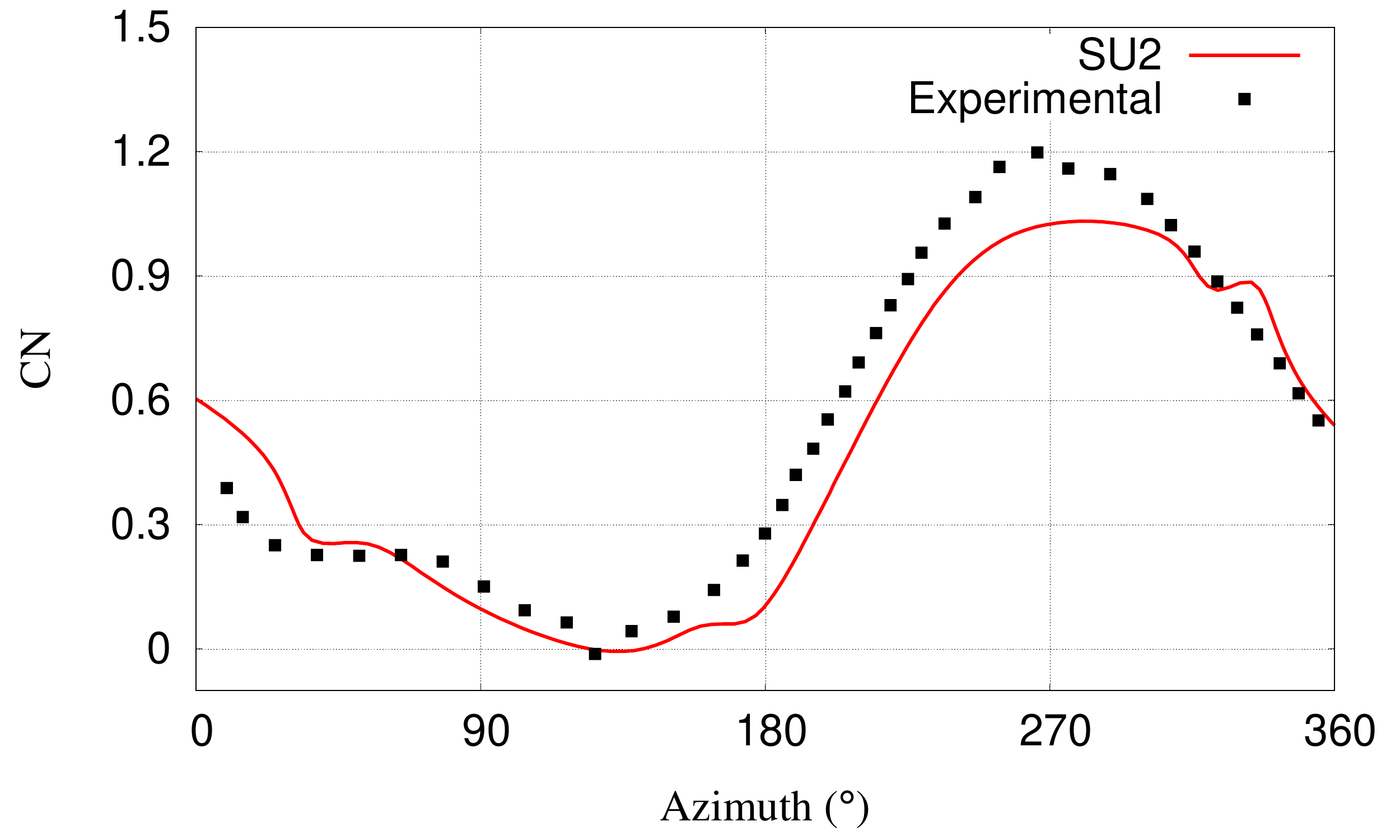}
	}
	\caption{Predictions compared against the measured data taken from the TAAT \cite{cross1988tip}.}
	\label{fig.10:AH1G_surfacePressure}
\end{figure}

\begin{figure}[!htb]
	\subfloat[Low speed forward flight.
	\label{fig.11:subfig-1:AH1G_surfPressureLS}]{%
		\includegraphics[width=1\textwidth]{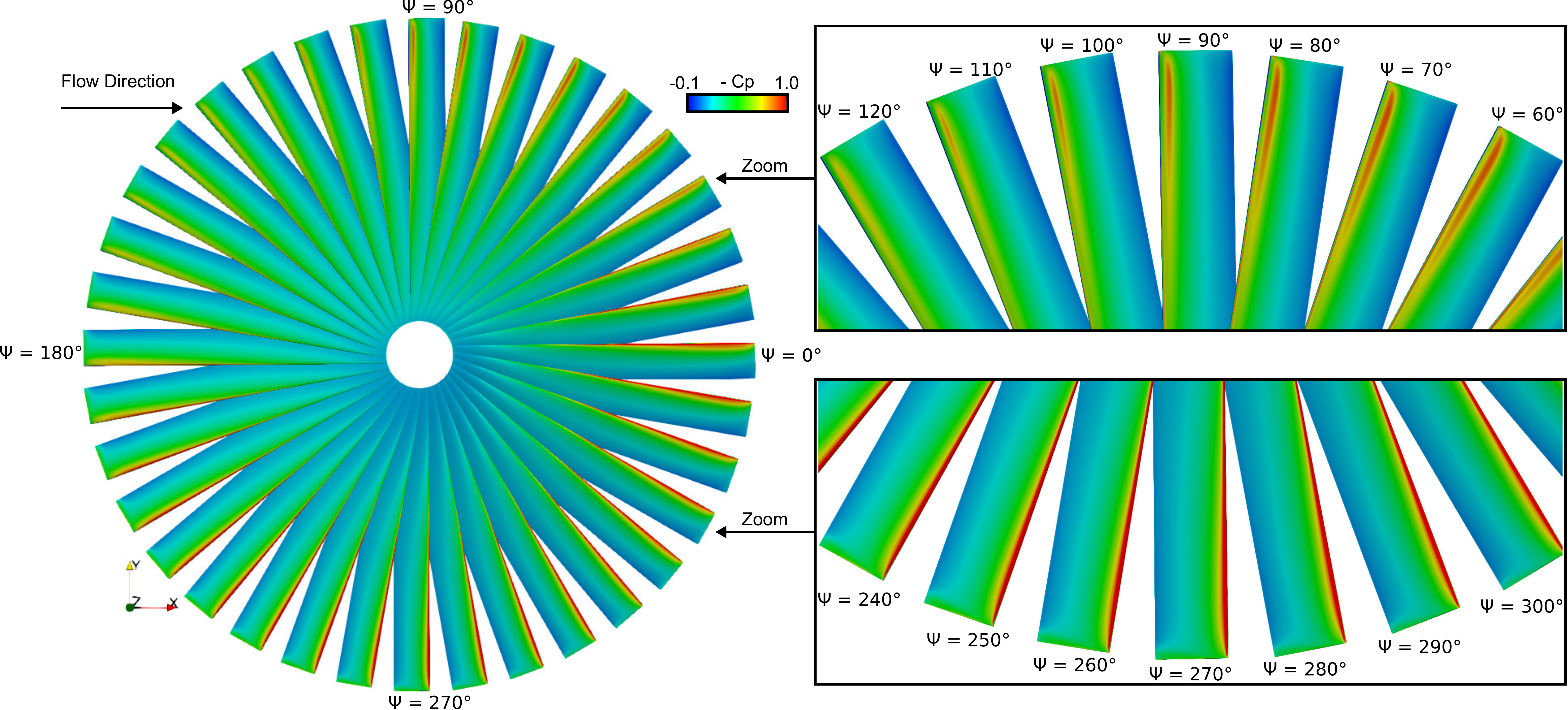}
	}
	\hfill
	\subfloat[High speed forward flight.
	\label{fig.11:subfig-2:AH1G_surfPressureHS}]{%
		\includegraphics[width=1\textwidth]{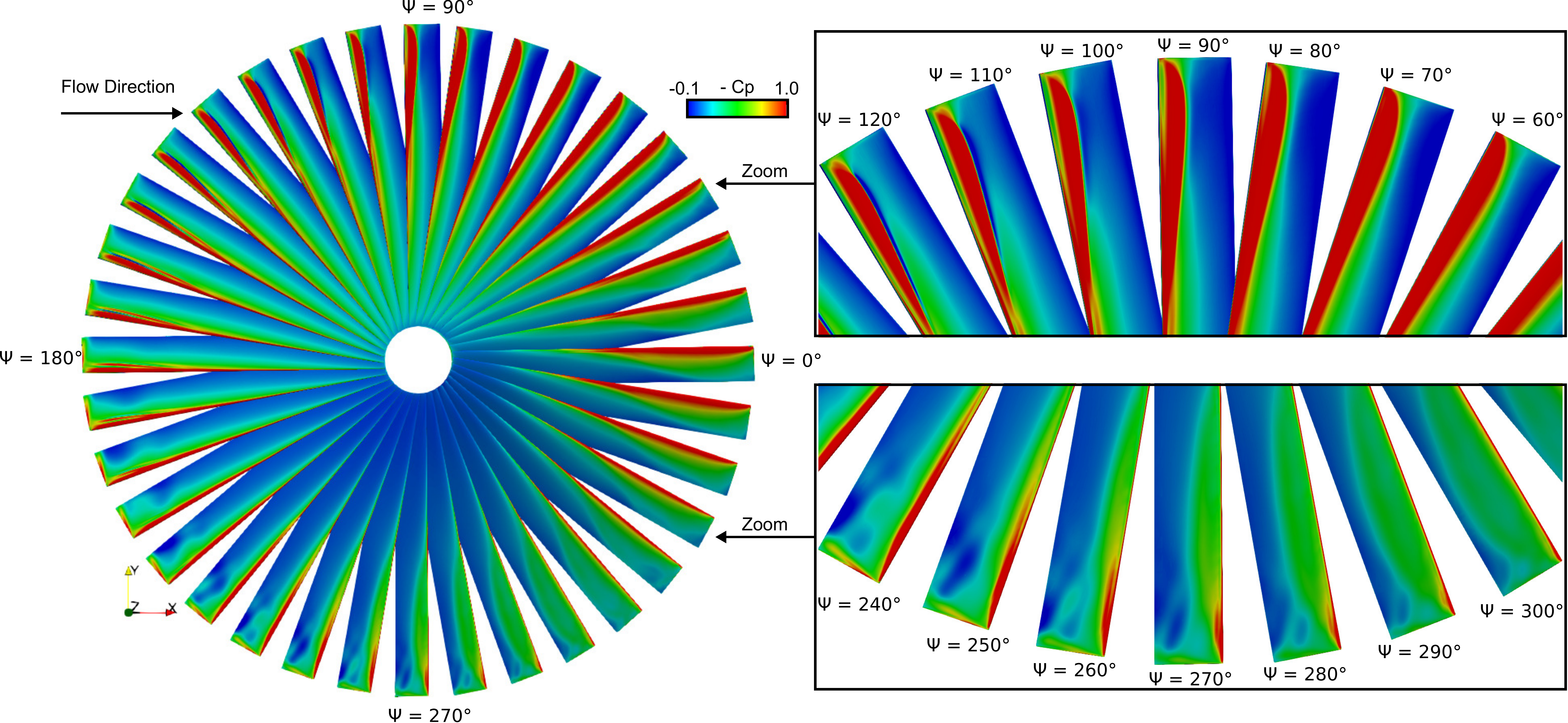}
	}
	\caption[$C_p$ distribution on the blade surface in forward flight]{Contour maps of the pressure coefficient on the upper blade surface during forward flight. Displaying the both low-speed and high-speed test conditions.}
	\label{fig.11:AH1G_surfPressure}
\end{figure}

\begin{figure}[hbt!]
	\centering
	\includegraphics[width=0.99\linewidth]{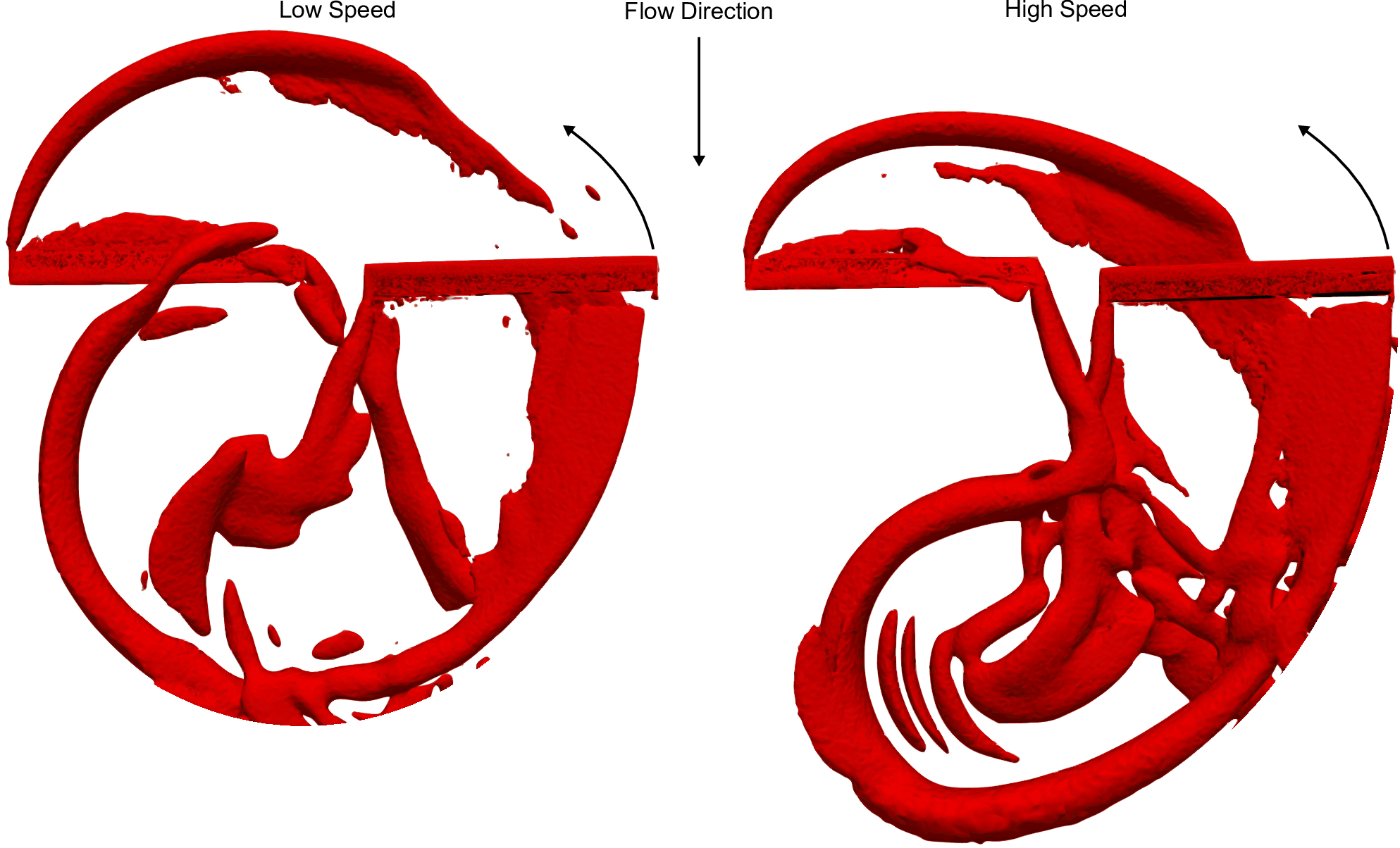}
	\caption[Q-criterion iso-surface during forward flight]{Iso-surface of the Q-criterion visualizing the near-field wake and blade tip vortices. Displaying the different flow field behaviour in low speed and high speed forward flight.}
	\label{fig12:AH1G_qCrit}
\end{figure}

\clearpage 
\listoffigures

\newpage
\listoftables

\newpage
\section*{References}
\bibliography{mybibfile}

\begin{thebibliography}{10}
\expandafter\ifx\csname url\endcsname\relax
  \def\url#1{\texttt{#1}}\fi
\expandafter\ifx\csname urlprefix\endcsname\relax\def\urlprefix{URL }\fi
\expandafter\ifx\csname href\endcsname\relax
  \def\href#1#2{#2} \def\path#1{#1}\fi

\bibitem{NTSB}
{Aviation Accident Reports - National Transportation Safety Board},
  \url{https://www.ntsb.gov/investigations/AccidentReports/Pages/aviation.aspx}
  (Accessed: May, 2020).

\bibitem{strawn200630}
R.~C. Strawn, F.~X. Caradonna, E.~P. Duque, 30 years of rotorcraft
  computational fluid dynamics research and development, Journal of the
  American Helicopter Society 51~(1) (2006) 5--21,
  https://doi.org/10.4050/1.3092875.

\bibitem{jespersen1997recent}
D.~Jespersen, T.~Pulliam, P.~Buning, Recent enhancements to overflow, in: 35th
  Aerospace Sciences Meeting and Exhibit, Reno, Navada, U.S.A, 06-09 January,
  1997, p. 644, https://doi.org/10.2514/6.1997-644.

\bibitem{anderson1994implicit}
W.~K. Anderson, D.~L. Bonhaus, An implicit upwind algorithm for computing
  turbulent flows on unstructured grids, Computers \& Fluids 23~(1) (1994)
  1--21, https://doi.org/10.1016/0045-7930(94)90023-X.

\bibitem{steijl2006framework}
R.~Steijl, G.~Barakos, K.~Badcock, A framework for cfd analysis of helicopter
  rotors in hover and forward flight, International journal for numerical
  methods in fluids 51~(8) (2006) 819--847, https://doi.org/10.1002/fld.1086.

\bibitem{biava2012simulation}
M.~Biava, L.~Vigevano, Simulation of a complete helicopter: A cfd approach to
  the study of interference effects, Aerospace Science and Technology 19~(1)
  (2012) 37--49, https://doi.org/10.1016/j.ast.2011.08.006.

\bibitem{antoniadis2012assessment}
A.~Antoniadis, D.~Drikakis, B.~Zhong, G.~Barakos, R.~Steijl, M.~Biava,
  L.~Vigevano, A.~Brocklehurst, O.~Boelens, M.~Dietz, et~al., Assessment of cfd
  methods against experimental flow measurements for helicopter flows,
  Aerospace Science and Technology 19~(1) (2012) 86--100,
  https://doi.org/10.1016/j.ast.2011.09.003.

\bibitem{srinivasan1992flowfield}
G.~R. Srinivasan, J.~Baeder, S.~Obayashi, W.~McCroskey, Flowfield of a lifting
  rotor in hover-a navier-stokes simulation, AIAA journal 30~(10) (1992)
  2371--2378, https://doi.org/10.2514/3.11236.

\bibitem{srinivasan1993turns}
G.~Srinivasan, J.~Baeder, Turns: A free-wake eule/navier-stokes numerical
  method for helicopter rotors, AIAA journal 31~(5) (1993) 959--962,
  https://doi.org/10.2514/3.49036.

\bibitem{gazaix2002elsa}
M.~Gazaix, A.~Jolles, M.~Lazareff, The elsa object-oriented computational tool
  for industrial applications, in: 23rd Congress of {ICAS}, Toronto, Canada,
  08-13 September, 2002, p. 220, {Corpus ID}: 56582723.

\bibitem{raddatz2005block}
J.~Raddatz, J.~K. Fassbender, {Block structured navier-stokes solver FLOWer},
  {MEGAFLOW - Numerical Flow Simulation for Aircraft Design. Notes on Numerical
  Fluid Mechanics and Multidisciplinary Design (NNFM)} Edition, Vol.~89,
  Springer, Berlin, Heidelberg, 2005, https://doi.org/10.1007/3-540-32382-1.

\bibitem{oberkampf1998issues}
W.~L. Oberkampf, F.~G. Blottner, Issues in computational fluid dynamics code
  verification and validation, AIAA journal 36~(5) (1998) 687--695,
  https://doi.org/10.2514/2.456.

\bibitem{economon2016su2}
T.~D. Economon, F.~Palacios, S.~R. Copeland, T.~W. Lukaczyk, J.~J. Alonso, Su2:
  An open-source suite for multiphysics simulation and design, Aiaa Journal
  54~(3) (2016) 828--846, https://doi.org/10.2514/1.J053813.

\bibitem{gori2017sliding}
G.~Gori, E.~van~der Weide, A.~Guardone, On conservation in compressible flow
  simulations using sliding mesh coupling, in: {VII International Conference on
  Computational Methods for Coupled Problems in Science and Engineering
  conference}, Rhodes Island, Greece, June 12-14, 2020.

\bibitem{wilcox1998turbulence}
D.~C. Wilcox, et~al., Turbulence modeling for CFD, 2nd Edition, DCW industries
  La Canada, CA, 1998, {ISBN 13: 9780963605153}.

\bibitem{sutherland1893lii}
W.~Sutherland, The viscosity of gases and molecular force, The London,
  Edinburgh, and Dublin Philosophical Magazine and Journal of Science 36~(223)
  (1893) 507--531, https://doi.org/10.1080/14786449308620508.

\bibitem{spalart1992one}
P.~Spalart, S.~Allmaras, A one-equation turbulence model for aerodynamic flows,
  in: 30th aerospace sciences meeting and exhibit, Reno, Nevada, U.S.A, 06-09
  January, 1992, p. 439, https://doi.org/10.2514/6.1992-439.

\bibitem{economon2013viscous}
T.~D. Economon, F.~Palacios, J.~J. Alonso, A viscous continuous adjoint
  approach for the design of rotating engineering applications, in: 21st {AIAA}
  computational fluid dynamics conference, San Diego, California, U.S.A, 24-27
  June, 2013, p. 2580, https://doi.org/10.2514/6.2013-2580.

\bibitem{versteeg2007introduction}
H.~K. Versteeg, W.~Malalasekera, An introduction to computational fluid
  dynamics: the finite volume method, 2nd Edition, Pearson education, 2007,
  {ISBN: 978-0-13-127498-3}.

\bibitem{jameson1981numerical}
A.~Jameson, W.~Schmidt, E.~Turkel, Numerical solution of the euler equations by
  finite volume methods using runge kutta time stepping schemes, in: 14th fluid
  and plasma dynamics conference, Palo Alto, California, U.S.A, 23-25 June,
  1981, p. 1259, https://doi.org/10.2514/6.1981-1259.

\bibitem{roe1981approximate}
P.~L. Roe, Approximate riemann solvers, parameter vectors, and difference
  schemes, Journal of computational physics 43~(2) (1981) 357--372,
  https://doi.org/10.1016/0021-9991(81)90128-5.

\bibitem{van1979towards}
B.~Van~Leer, Towards the ultimate conservative difference scheme. v. a
  second-order sequel to godunov's method, Journal of computational Physics
  32~(1) (1979) 101--136, https://doi.org/10.1016/0021-9991(79)90145-1.

\bibitem{jameson1991time}
A.~Jameson, Time dependent calculations using multigrid, with applications to
  unsteady flows past airfoils and wings, in: {10th Computational Fluid
  Dynamics conference}, Honolulu, Hawaii, U.S.A., p. 1596.

\bibitem{hall2002computation}
K.~C. Hall, J.~P. Thomas, W.~S. Clark, Computation of unsteady nonlinear flows
  in cascades using a harmonic balance technique, AIAA journal 40~(5) (2002)
  879--886, https://doi.org/10.2514/2.1754.

\bibitem{rubino2018adjoint}
A.~Rubino, M.~Pini, P.~Colonna, T.~Albring, S.~Nimmagadda, T.~Economon,
  J.~Alonso, Adjoint-based fluid dynamic design optimization in quasi-periodic
  unsteady flow problems using a harmonic balance method, Journal of
  Computational Physics 372 (2018) 220--235,
  https://doi.org/10.1016/j.jcp.2018.06.023.

\bibitem{rinaldi2015flux}
E.~Rinaldi, P.~Colonna, R.~Pecnik, Flux-conserving treatment of non-conformal
  interfaces for finite-volume discretization of conservation laws, Computers
  \& Fluids 120 (2015) 126--139,
  https://doi.org/10.1016/j.compfluid.2015.07.017.

\bibitem{rendall2010parallel}
T.~Rendall, C.~Allen, Parallel efficient mesh motion using radial basis
  functions with application to multi-bladed rotors, International journal for
  numerical methods in engineering 81~(1) (2010) 89--105,
  https://doi.org/10.1002/nme.2678.

\bibitem{wang2015improved}
G.~Wang, H.~H. Mian, Z.-Y. Ye, J.-D. Lee, Improved point selection method for
  hybrid-unstructured mesh deformation using radial basis functions, AIAA
  Journal 53~(4) (2015) 1016--1025, https://doi.org/10.2514/1.J053304.

\bibitem{xie2017efficient}
L.~Xie, H.~Liu, Efficient mesh motion using radial basis functions with volume
  grid points reduction algorithm, Journal of Computational Physics 348 (2017)
  401--415, https://doi.org/10.1016/j.jcp.2017.07.042.

\bibitem{morelli2020radial}
M.~Morelli, T.~Bellosta, A.~Guardone, Efficient radial basis function mesh
  deformation methods for aircraft icing, in: {In the 7th European Seminar on
  Computing}, Pilsen, Czech Republic, June 8-12, 2020.

\bibitem{caradonna1981experimental}
F.~X. Caradonna, C.~Tung, Experimental and analytical studies of a model
  helicopter rotor in hover, in: {Presented at the 6th European Rotorcraft and
  Powered Lift Aircraft Forum}, Bristol, England, September 16-19, 1980.

\bibitem{palacios2014stanford}
F.~Palacios, T.~D. Economon, A.~Aranake, S.~R. Copeland, A.~K. Lonkar, T.~W.
  Lukaczyk, D.~E. Manosalvas, K.~R. Naik, S.~Padron, B.~Tracey, et~al.,
  {Stanford university unstructured (SU2): Analysis and design technology for
  turbulent flows}, in: {In 52nd Aerospace Sciences Meeting}, National Harbor,
  Maryland, January 13-17, 2014, p. 0243, https://doi.org/10.2514/6.2014-0243.

\bibitem{cross1988tip}
J.~L. Cross, M.~E. Watts, Tip aerodynamics and acoustics test: a report and
  data survey, National Aeronautics and Space Administration, Scientific and
  Technical Information Division, Reference Publication No. 1179, Ames Research
  Center, Moffett Field, California, December, 1988.

\end{thebibliography}

\end{document}